 \documentclass{ajour}
 \usepackage{amssymb}

\newtheorem{lemmab}[theorem]{Lemma}
\newtheorem{propb}[theorem]{Proposition}
\newtheorem{corb}[theorem]{Corollary}
\newtheorem{defb}[theorem]{Definition}

\newtheorem{openquestion}[theorem]{Question}

\begin{document}

%


\authorrunninghead{Arturo Magidin}
\titlerunninghead{Amalgams of nilpotent groups of class two}





\title{Amalgams of nilpotent groups of class two%
\thanks*{2000 MSC: Primary 20E06,
20F18, 20D15, 08B25; Secondary 20E10.}%
}


\author{Arturo Magidin}
\affil{Instituto de Matem\'aticas,
       Universidad Nacional Aut\'onoma de M\'exico\\
       Circuito Exterior, Ciudad Universitaria,
       04510 Mexico City, Mexico}

\email{magidin@matem.unam.mx}

\abstract{We give necessary and sufficient conditions for weak and
strong embeddability of amalgams in each subvariety of ${\cal N}_2$,
the class of all nilpotent groups of class at most two; this
generalizes B.~Maier's result for ${\cal N}_2$. We also discuss
dominions, and characterize the weak, strong, and special amalgamation
bases for each subvariety, contrasting the resulting classes with one
another.}


\begin{article}


\zerosection{Introduction}

In this paper we will prove analogues of B.~Maier's characterization of
weak and strong embeddability of amalgams in~$\mathcal{N}_2$
\cite{amalgone,amalgtwo} for the subvarieties of $\mathcal{N}_2$. We
will also give analogues of D.~Saracino's characterization of weak and
strong amalgamation bases \cite{saracino}, the author's work on
dominions \cite{nildoms} and on amalgamation bases in some varieties
of nil-2 groups \cite{closures,absclosed}. Definitions will be
recalled below. The results are obtained by extending the methods used
by the cited authors.

Groups will be written multiplicatively, unless otherwise
specified. We will use $Z$ to denote the infinite cyclic group, which
we also write multiplicatively.  All maps are assumed to be group
morphisms unless we explicitly note otherwise. The multiplicative
identity of a group $G$ will be denoted by $e$, and we will use $e_G$
if there is danger of ambiguity. For a group $G$ and elements $x,y\in
G$, $x^y$ represents $y^{-1}xy$, and the commutator of $x$ and $y$ is
$[x,y]=x^{-1}y^{-1}xy$; note that $[x,y]^{-1}=[y,x]$, and
$x^y=x[x,y]$. Given subsets $A,B$ of $G$, not necessarily subgroups,
$[A,B]$ denotes the subgroup of $G$ generated by all commutators
$[a,b]$ with $a\in A$ and $b\in B$. The commutator subgroup of $G$ is
the subgroup $[G,G]$, which is also denoted by~$G'$. The center
of $G$ is denoted by $Z(G)$. All commutators will be written
left-normed, so $[x,y,z]=[[x,y],z]$.  Given a positive integer $k>0$,
$G^k$ is the subgroup of $G$ generated by all $k$-th powers of
elements of $G$, that is $G^k = \left\langle g^k \,|\, g\in
G\right\rangle$. We also set $G^0=\{e\}$.

So that our statements can be written generally, we will sometimes say
that an element is of ``order zero'' to mean it is of infinite order,
and of ``exponent zero'' to mean it has arbitrary (possibly, but not
necessarily, infinite) order. We will say that $0$ divides~$0$, and no
other integer, so $\alpha\equiv\beta \pmod{0}$ means $\alpha=\beta$.
Also, recall that for any nonnegative integer $b$, ${\rm gcd}(0,b)=b$
and ${\rm lcm}(0,b)=0$.

Recall that given a class ${\cal C}$ of groups, a ${\cal C}$-amalgam
of two groups $A$ and $B$ in ${\cal C}$ with core $D$ consists of
groups $A$, $B$, and $D$, all in ${\cal C}$, equipped with one to one
group morphisms
\begin{eqnarray*}
\Phi_A\colon D &\to& A\\
\Phi_B\colon D&\to& B.
\end{eqnarray*}
To simplify notation, we denote this situation by $(A,B;D)$, and we
think of $D$ as a subgroup of both $A$ and $B$; i.e.~$D\subseteq A$
and $D\subseteq B$, with $\Phi_A$ and $\Phi_B$ the subgroup
inclusions. To say that the amalgam $(A,B;D)$ is \textit{(weakly)
embeddable (in ${\cal C}$)} means that there exists a group $M$ in
${\cal C}$ and one to one group morphisms
\[\lambda_A\colon A\to M,\qquad \lambda_B\colon B\to M,\qquad
\lambda\colon D\to M,\]
such that
\[\lambda_A\circ\Phi_A = \lambda \qquad{\rm and}\qquad
\lambda_B\circ\Phi_B=\lambda.\]
We also say that the group $M$ \textit
{realizes the amalgam.} Note that there is nothing in the
definition to prevent an element $x\in A\setminus D$ to be identified
in $M$ with an element $y\in B\setminus D$.  If, in addition, we have
\[\lambda_A(A)\cap\lambda_B(B)=\lambda(D),\] then we say that the
embedding is a \textit{strong} embedding, and that the amalgam is
\textit{strongly embeddable (in ${\cal C}$).} 

A classical theorem of Schreier states that if ${\cal C}$ is the class
of all groups, or the class of all finite groups, then every ${\cal
C}$-amalgam is strongly embeddable. It is not hard to verify that if
${\cal C}$ is the class of all abelian groups, or the class of all
abelian groups of a given exponent, then we again have that all ${\cal
C}$-amalgams are strongly embeddable. However, this situation is
actually fairly rare. For example, not every amalgam of finite
$p$-groups is embeddable into a finite $p$-group \cite{higmanpgroups}.

So, given a class ${\cal C}$, it is natural to look for necessary and
sufficient conditions for the embeddability of an amalgam in
${\cal C}$. For example, the problem was studied by Higman for finite
$p$-groups \cite{higmanpgroups}, by Leinen for solvable groups
\cite{leinen}, and by Maier \cite{amalgone,amalgtwo} for the class of
nilpotent groups of class two, among others. Ideally, we want to be
able to characterize embeddability by conditions ``internal'' to the
amalgam, that is conditions which involve only the groups $A$, $B$,
and~$D$; examples of this kind of conditions are Maier's conditions,
and Higman's conditions for $p$-groups. This in opposition to what one
might call conditions ``external'' to the amalgam, which are
conditions which invoke the existence of other groups or morphisms
which are not inherent in $A$, $B$, or $D$. For example, see condition
(A)(c) of Theorem~2 in \cite{leinen}.

Closely related to amalgams is the notion of \textit{amalgamation
bases}. Given a class ${\cal C}$, a group $G\in{\cal C}$ is a
\textit{weak amalgamation base} (for ${\cal C}$) if every ${\cal
C}$-amalgam $(A,B;G)$ with core $G$ is weakly embeddable in ${\cal
C}$; and it is a \textit{strong amalgamation base} (for ${\cal C}$) if
every ${\cal C}$-amalgam $(A,B;G)$ with core $G$ is strongly
embeddable in ${\cal C}$. Another natural question is to characterize
the amalgamation bases for a given ${\cal C}$.

Clearly, any strong base in ${\cal C}$ is necessarily a weak base in
${\cal C}$, but the converse need not hold. The gap between the two is
``measured'' by the \textit{special amalgams}. An amalgam $(A,B;D)$ is
\textit{special} if there exists an isomorphism $\varphi\colon A\to
B$, such that $\varphi\circ\Phi_A=\Phi_B$. In that case, we usually
write $(A,A;D)$ with $\varphi={\rm id}_A$ being understood. Clearly, a
special amalgam is always weakly embeddable, so we are interested in
strong embeddability of special amalgams. A group $G\in{\cal C}$ is a
\textit{special amalgamation base} (for ${\cal C}$) if and only if
every special ${\cal C}$-amalgam $(A,A;G)$ is \textit{strongly}
embeddable. We then have:

\begin{theorem}[see for example \cite{episandamalgs}]
\label{weakstrongdifference}
Let ${\cal C}$ be a class of groups. A group $G\in{\cal C}$ is a strong
amalgamation base in ${\cal C}$ if and only if it is both a weak and a special
amalgamation base in ${\cal C}$.
\end{theorem}

Among the most natural classes to work in are the \textit{varieties}
of groups. A class ${\cal C}$ is a variety if and only if it is closed
under taking homomorphic images, subgroups, and
arbitrary direct products of members of ${\cal C}$. Alternatively, a
variety is the collection of all groups which satisfy a given set of
laws or identities. Thus, the variety of all abelian groups is
determined as the collection of all groups that satisfy the identity
$xy=yx$; all groups of exponent $n$ are determined by the identity
$x^n=e$; etc. We direct the reader to Hanna Neumann's excellent book
\cite{hneumann} for more information on varieties of groups.

Let ${\cal N}_2$ denote the variety of all nilpotent groups of class
at most two; that is, groups $G$ such that $[G,G]\subseteq Z(G)$, or
equivalently, for which the identity $[x,y,z]=e$ holds. It is easy
to verify that for any nil-2 group (i{.}e{.} any nilpotent group of
class at most two) the following identities hold, so we will use them
without comment throughout the paper:

\begin{propb}
Let $G\in {\cal N}_2$. For all $x,y,z\in G$ and all integers $n$:
\begin{itemize}
\item[(a)]~$[xy,z]=[x,z][y,z]$; $[x,yz]=[x,y][x,z]$.
\item[(b)]~$[x^n,y]=[x,y]^n=[x,y^n]$.
\item[(c)]~$(xy)^n = x^n y^n[y,x]^{n(n-1)/2}$.
\item[(d)]~The value of $[x,y]$ depends only on the congruence
classes of $x$ and~$y$ modulo $G'$ (in fact, modulo $Z(G)$). 
\end{itemize}
\end{propb}

To simplify notation, we set ${n \choose 2} = \frac{n(n-1)}{2}$
for every integer $n$.

It is worth noting that by (c) above, any element of $G^n$ can be
written as $g^ng'$, with $g\in G$ and $g'\in G'$.

In \cite{amalgone} and \cite{amalgtwo}, B.~Maier gave necessary and
sufficient conditions for an ${\cal N}_2$-amalgam $(A,B;D)$ to be
weakly or strongly embeddable in ${\cal N}_2$. D.~Saracino
characterized the weak and strong amalgamation bases for ${\cal N}_2$
in~\cite{saracino}. The author characterized the special amalgamation
bases in~\cite{absclosed}, and characterized the weak, strong, and
special amalgamation bases in the varieties of all ${\cal N}_2$-groups
of a given odd exponent~$m$ in~\cite{closures}.

A subvariety of a variety ${\cal V}$ is a class of groups ${\cal W}$
with ${\cal W}\subseteq{\cal V}$, and ${\cal W}$ a variety of groups
in its own right. In the case of ${\cal N}_2$, a full description of
all subvarieties is given by:

\begin{theorem}[see for example \cite{jonsson}, \cite{classifthree}]
Every subvariety of ${\cal N}_2$ may be defined by the identities
\begin{equation}
x^m = [x_1,x_2]^n = [x_1,x_2,x_3] = e
\label{definingeqs}
\end{equation}
for unique nonnegative integers $m$ and $n$ satisfying $n|m/{\rm
gcd}(2,m)$, yielding a bijection between pairs of nonnegative integers
$(m,n)$ satisfying this condition, and subvarieties of ${\cal N}_2$.
\end{theorem}

The case $m=n=0$ corresponds to the full variety ${\cal N}_2$. When
$n=1$, we obtain the varieties of abelian groups of exponent $m$
(where ``abelian groups of exponent zero'' just means ``abelian
groups''). For odd $m>0$, the variety of all nil-2 groups of exponent
$m$ is given by $(m,m)$, and if $m=2k>0$, then the variety of all nil-2
groups of exponent $m$ is given by $(m,k)=(2k,k)$.

To see why the possibly mysterious factor ${\rm gcd}(2,m)$ appears,
note that if all elements are of exponent $m$, then
\[e  =  (xy)^m 
 =  x^m y^m [y,x]^{{m\choose 2}}
 =  [y,x]^{{m\choose 2}}.\]
If $m$ is odd, ${m\choose 2}$ is a multiple of $m$, so we
get that $[y,x]^{{m\choose 2}}$ is trivial. But if $m$ is even, then
$[y,x]$ must be of exponent dividing $\frac{m}{2}$.

We will denote the subvariety corresponding to $(m,n)$ simply by the
pair $(m,n)$. Thus, if we write $G\in (m,n)$ we mean that $G$
satisfies the identities in~(\ref{definingeqs}). Note that if $G\in
(m,n)$, then all $n$-th powers are necessarily central, that is
$G^n\subseteq Z(G)$. Also note that $(m,n)\subseteq (m',n')$ if and
only if $n|m/{\rm gcd}(m,2)$, $n'|m'/{\rm gcd}(m',2)$, $m|m'$, and
$n|n'$. If we write $(m,n)$, we will assume that $n|m/{\rm gcd}(m,2)$.

The subvarieties of ${\cal N}_2$ form a $01$-lattice, with minimum
$(1,1)$, the variety consisting only of the trivial group, and maximum
$(0,0)={\cal N}_2$. The meet and join are defined, respectively, by,
\begin{eqnarray*}
(m,n) \wedge (m',n') & = & \left({\rm gcd}(m,m'),{\rm
gcd}(n,n')\right),\\
(m,n) \vee (m',n') & = & \left({\rm lcm}(m,m'),{\rm
lcm}(n,n')\right).
\end{eqnarray*}

Let ${\cal V}$ be a variety, and ${\cal W}$ a subvariety of ${\cal
V}$. If we have an amalgam $(A,B;D)$ of ${\cal W}$-groups, it need not
be embeddable in ${\cal W}$ even if it is embeddable in ${\cal V}$
(the converse, of course, always holds). If a group $G\in{\cal W}$ is
an amalgamation base in ${\cal V}$, it need not be an amalgamation
base in ${\cal W}$ (for, perhaps, a given amalgam of ${\cal W}$-groups
cannot be embedded into a ${\cal W}$-group, even if it can be
embedded into a ${\cal V}$-group). Nor is it the case that if $G$ is
an amalgamation base in ${\cal W}$, it need be one in ${\cal V}$.
Thus, we want to give necessary and sufficient conditions for
embeddability in $(m,n)$, contrast them with those in $(m',n')$,
and do the same with amalgamation bases.

Here is one example to keep in mind:

\begin{example}
\label{guidingex}
Consider the following situation: let $M$ be the nil-2
group presented by:
\[M = \Bigl\langle x,y \,\Bigm|\,
x^4=y^4=[x,y]^4=[x,y,x]=[x,y,y]=e\Bigr\rangle\]
and consider the subgroups $A=\langle x\rangle$ and $B=\langle
x^2,y\rangle$. Their intersection is a cyclic group of order 2,
generated by $x^2$. Call it $D$. The amalgam $(A,B;D)$ is thus
strongly embeddable in ${\cal N}_2$. Both $A$ and $B$ are of exponent
$4$, so they lie in $(4,2)$, while $M\in(8,4)$. We claim, however,
that $(A,B;D)$ is not even weakly embeddable in $(4,2)$. To see this,
note that $x^2$ is a square in $A$, but not central in $B$. If we
could embed the amalgam $(A,B;D)$ into a $(4,2)$-group $K$, then in
$K$ all squares would be central, so $x^2\in D$ would necessarily be
central in $K$, but it does not even centralize $B$. Thus
$(A,B;D)$ cannot be even weakly embeddable in $(4,2)$, even though it
is strongly embeddable in ${\cal N}_2$ (and in $(8,4)$).
\end{example}

For the remainder of this paper, all groups will be assumed to lie in
${\cal N}_2$ unless otherwise specified. If we say a group $G$ lies in
$(m,n)$, and we give a presentation for the group, it will be understood
that it is a presentation in $(m,n)$; that is, the identities of
$(m,n)$ will be imposed on the group, as well as all the relations
specified in the presentation.

In Section~1 we will give preliminary results. Among them we will
prove some results on adjunction of roots. In Section~2 we will give
necessary and sufficient conditions for weak and strong embeddability
of amalgams in $(m,n)$; the result with $m=n=0$ gives Maier's
conditions for ${\cal N}_2$, and all conditions are internal. The
proofs are based on Maier's methods. In Section~3 we use the
embeddability result to give a description of dominions in these
varieties (definitions will be recalled there). With this information,
we will return in Section~4 to the question of embeddability, and we
will contrast embeddability in $(m,n)$ with that in~$(m',n')$. In
Section~5 we turn to the characterization the weak and strong
amalgamation bases in $(m,n)$, and prove the two classes are actually
equal. We will give the description, and then make some reductions in
special cases and give examples. Finally, in Section~6 we do the same
for the special amalgamation bases. Again, the characterizations of
$m=n=0$ will yield Saracino's and the author's previous results. We
will close by proving that a group in $(m,n)$ has an absolute closure
if and only if it is already absolutely closed (definitions are
recalled below).

\section{Preliminary results}

One important property of nilpotent groups is that a torsion
nilpotent group can be decomposed into the direct product of its
$p$-parts. The following lemma tells us we can study amalgamation
properties in this situation by dealing with the $p$-parts separately.

\begin{propb}
\label{prop:pparts} Let $(G,K;H)$ be an amalgam of $(m,n)$-groups,
with $m>0$. Then $(G,K;H)$ is weakly (resp.~strongly) embeddable into
an $(m,n)$ group if and only if for every prime $p$ dividing $m$, the
amalgam $(G_p,K_p;H_p)$ is weakly  (resp.~strongly) embeddable in
$(p^{{\rm ord}_p(m)},p^{{\rm ord}_p(n)})$, where $G_p$ is the $p$-part
of $G$, and likewise for $K_p$ and $H_p$.
\end{propb}
\begin{proof}
Since $G$, $K$, and $H$ are of finite exponent $m$, they are the
direct product of their $p$-parts. Thus $G=\prod G_p$, $K=\prod K_p$
and $H=\prod H_p$. Clearly, an embedding of $(G,K;H)$ into $M\in
(m,n)$ provides embeddings for $(G_p,K_p;H_p)$ into $M_p$ for each
prime $p$. Conversely, if we have embeddings into $M_p$ for each
$(G_p,K_p;H_p)$, then $\prod M_p$ gives an embedding for $(G,K;H)$.
\end{proof}

Thus, in the case when $m>0$, we may restrict ourselves to varieties
of the form $(p^{a+b},p^a)$ with $p$ a prime, $a,b\geq 0$ ($b>0$
if $p=2$). 

In addition, another important reduction is given by the following
result, which follows from the usual compactness results of logic:

\begin{propb}[Lemma 5 in \cite{amalgone}]
\label{fgamalgams}
Let $A$ and $B$ be groups in a variety ${\cal V}$, and let $(A,B;D)$
be a ${\cal V}$-amalgam. The amalgam is weakly (resp. strongly)
embeddable in ${\cal V}$ if and only if for each finite collection of elements
$a_1,\ldots,a_m\in A$ and $b_1,\ldots,b_n\in B$, the amalgam 
\[\bigl(\langle a_1,\ldots,a_m\rangle,\langle b_1,\ldots,b_n\rangle;
\langle a_1,\ldots,a_m\rangle\cap\langle b_1,\ldots,b_n\rangle\bigr)\]
is weakly (resp.~strongly) embeddable in ${\cal V}$.
\end{propb}

\begin{remark}
The fact that ${\cal V}$ is a variety is important. The result may not
hold even for closely related classes, such as pseudovarieties
(classes closed under subgroups, quotients, and \textit{finite} direct
products). For example, compare Example~8.87 in \cite{nildoms}, with
Theorem~3.11 in \cite{fgnil}. 
\end{remark}

\subsection*{Coproducts}

Another advantage of working in varieties is the existence of
coproducts. Given a variety of groups ${\cal V}$, and two groups $G$
and $K$ in ${\cal V}$, their \textit{$\mathcal{V}$-coproduct},
$G\amalg^{{\cal V}} K$ is the unique $\mathcal{V}$-group (up to
isomorphism) equipped with embeddings $i_G\colon G\to G\amalg^{{\cal
V}}K$ and $i_K\colon K\to G\amalg^{{\cal V}}K$, and the following
universal property: given any pair of morphisms $f\colon G\to M$ and
$g\colon K\to M$ with $M\in {\cal V}$, there exists a unique morphism
$\varphi\colon G\amalg^{{\cal V}}K\to M$ with $\varphi\circ i_G=f$ and
$\varphi\circ i_K=g$. One way to construct the coproduct is by taking
the free product $G*K$ of $G$ and $K$, and then moding out by the
verbal subgroup corresponding to $\mathcal{V}$, i.e.~the least normal
subgroup $N\triangleleft G*K$ such that $(G*K)/N\in\mathcal{V}$.

From the coproduct, one can construct \textit{amalgamated
coproducts}. Given an amalgam of ${\cal V}$-groups $(A,B;D)$, the
amalgamated coproduct of $A$ and $B$ over $D$ is the unique (up to
isomorphism) ${\cal V}$-group $A\amalg_D^{{\cal V}} B$, equipped with
morphisms $\psi_A\colon A\to A\amalg_D^{{\cal V}} B$ and $\psi_B\colon
B\to A\amalg_D^{{\cal V}} B$ with $\psi_A(d)=\psi_B(d)$ for each $d\in
D$, and the following universal property: given any pair of morphisms
$f\colon A\to M$ and $g\colon B\to M$, with $M\in{\cal V}$, and such
that for all $d\in D$, $f(d)=g(d)$, there exists a unique morphism
$\varphi\colon A\amalg_D^{{\cal V}} B\to M$ such that
$\varphi\circ\psi_A=f$ and $\varphi\circ \psi_B=g$. The amalgamated
coproduct can be constructed as the quotient of the coproduct
$A\amalg^{{\cal V}} B$ modulo the least normal subgroup which contains
all elements $i_A(d)(i_B(d))^{-1}$ for each $d\in D$.  Note that it may
occur that the morphisms into the amalgamated coproduct are not
injections.

If we are working in a variety $\mathcal{V}$, and we examine whether or not the
amalgam $(A,B;D)$ is embeddable, the obvious candidate for $M$ is the
$\mathcal{V}$-coproduct with amalgamation of $A$ and $B$ over $D$,
$A\amalg^{{\cal V}}_{D} B$. Because of the universal
property, the amalgam is weakly (resp.~strongly) embeddable in~${\cal
V}$ if and only if it is weakly (resp.~strongly) embeddable into
$A\amalg_{D}^{{\cal V}}B$. 

Given $G,K\in{\cal N}_2$, every element of their coproduct
$G\amalg^{{\cal N}_2} K$ has a unique expression of the form
$\alpha\beta\gamma$, where $\alpha\in G$, $\beta\in K$, and $\gamma\in
[G,K]$, the `cartesian.' A theorem of T.~MacHenry \cite{machenry}
states that the cartesian subgroup $[G,K]$ of $G\amalg^{{\cal N}_2}K$
is isomorphic to the tensor product $G^{\rm ab}\otimes K^{\rm ab}$, by
the mapping that sends $[g,k]$ to $\overline{g}\otimes \overline{k}$.

We can also provide a similar description of the coproduct in $(m,n)$:

\begin{propb}
\label{prop:coprod}
Let $G,K\in (m,n)$. Then
\[ G\amalg^{(m,n)}K \cong \left(G\amalg^{{\cal N}_2} K\right)\Bigm/ [G,K]^n.\]
In particular, every element of $G\amalg^{(m,n)} K$ can be written
uniquely as $\alpha\beta\gamma$, with $\alpha\in G$, $\beta\in K$, and
$\gamma$ in the `cartesian', which is isomorphic to $[G,K]/[G,K]^n$.
\end{propb}

\begin{proof} It suffices to show that $(G\amalg^{{\cal N}_2}K)/[G,K]^n$
lies in $(m,n)$ and has the corresponding universal property. That it
lies in $(m,n)$ follows because both $G$ and $K$ are of exponent $m$,
and all commutators are either in $G'$, $K'$, or the cartesian. The
former two are already of exponent $n$, and the latter is of exponent
$n$ since we have killed the $n$-th powers explicitly. That $G$ and
$K$ are embedded into the quotient follows from the description given
above.

Given maps $f\colon G\to M$ and $g\colon K\to M$, with $M\in(m,n)$,
the universal property of $G\amalg^{{\cal N}_2} K$ gives a unique map
to $M$, which will factor through the quotient $G\amalg^{{\cal
N}_2}K/[G,K]^n$, since in $M$ all $n$-th powers of commutators are
central. It is now easy to verify that this map is indeed unique from
the quotient, giving that this is in fact isomorphic to the coproduct
in $(m,n)$. The final statement now follows.
\end{proof}

And we can generalize T.~MacHenry's result to give a
description of the cartesian:

\begin{propb}
\label{prop:machenrygen}
Let $G,K\in(m,n)$. The cartesian $[G,K]$ of the coproduct
$G\amalg^{(m,n)}K$ is isomorphic to
\[\frac{G}{G^nG'}\otimes \frac{K}{K^nK'}.\]
\end{propb}
\begin{proof}
By Proposition~\ref{prop:coprod}, the cartesian of $G\amalg^{(m,n)}K$
is isomorphic to the cartesian of $G\amalg^{{\cal N}_2}K$ modulo its
$n$-th power. Using T.~MacHenry's isomorphism, this gives, in $G\amalg^{(m,n)}K$:
\begin{eqnarray*}
[G,K] & \cong & \frac{G^{\rm ab}\otimes K^{\rm ab}}{\left(G^{\rm
ab}\otimes K^{\rm ab}\right)^n}\\
& \cong & \frac{G^{\rm ab}}{\left(G^{\rm ab}\right)^n}\otimes
\frac{K^{\rm ab}}{\left(K^{\rm ab}\right)^n}\\
& \cong & \frac{G}{G^nG'}\otimes \frac{K}{K^nK'}.
\end{eqnarray*}
\end{proof}

\begin{corb}[cf. Lemma~3 in \cite{amalgone}]
\label{whatcommutes}
Let $G\in(m,n)$. For each element $g\in G$, the following are equivalent:
\begin{itemize}
\item[(a)]~$g\in G^nG'$.
\item[(b)]~$[g,c]=e$ in $G\amalg^{(m,n)}\langle c\rangle$, where
$\langle c\rangle$ is a cyclic group of order $m$ (infinite cyclic if
$m=0$).
\item[(c)]~$g\in Z(H)$ for any $(m,n)$-overgroup $H$ of $G$.
\end{itemize}
\end{corb}

\begin{proof}
$(a)\Rightarrow(c)$ For any $(m,n)$-overgroup $H$ of $G$, we have
\[g\in G^nG'\subseteq H^nH'\subseteq Z(H).\]

$(c)\Rightarrow(b)$ Note that $G\amalg^{(m,n)}\langle c\rangle$ is an
 $(m,n)$-overgroup of $G$, so $g$ is central there, and thus $[g,c]=e$.

$(b)\Rightarrow (a)$ Assume that $g\notin G^nG'$. Since the
cartesian of $G\amalg^{(m,n)}\langle c\rangle$ is isomorphic to
\[\frac{G}{G^nG'}\otimes \frac{\langle
c\rangle}{\langle c^n\rangle} \cong \frac{G}{G^nG'}\]
it follows that $[g,c]$ is nontrivial in $G\amalg^{(m,n)}\langle
c\rangle$, as desired.
\end{proof}

\subsection*{Central amalgams}

There is one class of amalgams for which strong embeddability is
easy. These are the amalgams where the core is central in both
factors:

\begin{propb}
\label{centralamalgs}
Let ${\cal C}$ be a class of groups closed under quotients and finite
direct products. If $(A,B;D)$ is an amalgam in ${\cal C}$, and $D$ is
central in both $A$ and $B$, then the amalgam is strongly embeddable
into $A\times_D B$, where:
\[A\times_D B = \frac{A\times B}{\Bigl\{(d,d^{-1})\in A\times B\Bigm|
d\in D\Bigr\}}.\]
\end{propb}

If $G$ is a group, and $H$ a subgroup of $G$, we let $C_G(H)$ be the
centralizer of $H$ in $G$; that is, the set of all $g\in G$ which
commute with all elements of $H$. The following results are quoted
here for convenience, and we direct the reader to Maier's paper
\cite{amalgone} for their proofs:

\begin{lemmab}[Lemma~1 in \cite{amalgone}]
\label{lemmaone}
Let $U\leq G$ and $X\leq C_G(U)$. Then the subgroup $\langle
U,X\rangle$ of $G$ is equal to $UX$, and isomorphic to $U\times_{U\cap
X} X$.
\end{lemmab}

\begin{lemmab}[Korollar~1 in~\cite{amalgone}]
\label{korollarone}
Let $X$ be central in both $A$ and $B$, and $U\leq A$. The subgroup
$\langle U,B\rangle$ in $A\times_X B$ is isomorphic to $U\times_{U\cap
X} B$.
\end{lemmab}

We say a subgroup $U$ of $G$ is \textit{co-central} if $G=\langle
U,X\rangle$, with $X$ central in~$G$.

\begin{lemmab}[Lemma~2 in~\cite{amalgone}]
\label{lemmatwo}
Let $U$ be co-central in both $G$ and~$H$, i.e.~$G=\langle U,X\rangle$
and $H=\langle U,Y\rangle$, with $X\leq Z(G)$, $Y\leq Z(H)$. Then
\[P = \frac{U\times X\times Y}{\Bigl\{ (xy,x^{-1},y^{-1}) \Bigm| x\in
U\cap X, y\in U\cap Y\Bigr\}}\]
is isomorphic to $G_Y=G\times_{U\cap Y} Y$ and to $H_X=H\times_{U\cap
X} X$, and the isomorphism identifies $G\leq G_Y$ (resp.~$H\leq H_X$)
with $\langle U,X\rangle\subseteq H_X$ (resp.~$\langle
U,Y\rangle\subseteq G_Y$).
\end{lemmab}

\begin{lemmab}[cf. Korollar~2 in~\cite{amalgone}]
\label{korollartwo}
Let $A,B\in(m,n)$, and let $D$ be a subgroup of both $A$ and $B$, with
$A^nA'\cap D\subseteq Z(B)$, $B^nB'\cap D\subseteq Z(A)$. Then the subgroups
$\langle D,A^nA',B^nB'\rangle$ of $A_0=A\times_{D\cap B^nB'} B^nB'$
and of $B_0=A^nA'\times_{A^nA'\cap D} B$ are isomorphic, and
$A_0,B_0\in (m,n)$.
\end{lemmab}
\begin{proof}
Set $D=U$, $X=A^nA'$, $Y=B^nB'$, and apply Lemma~\ref{lemmatwo}. That
they both lie in $(m,n)$ follows because they are obtained as
quotients of products of groups in $(m,n)$.
\end{proof}

\subsection*{Root adjunction}

To characterize amalgamation bases, we will need to know when it is
possible to adjoin roots to elements of $G$; i.e.~give an overgroup
$K\in(m,n)$ such that a given set of elements of $G$ lie in
$K^q$. Mostly, we will need to know if we can adjoin a $q$-th root modulo
the commutator subgroup, so we will ask that the given
set of elements lie in $K^qK'$.

The basic result on root adjunction in ${\cal N}_2$ is due to
Saracino, and we quote it here for reference::

\begin{lemmab}[Theorem 2.1 in \cite{saracino}]
\label{lemma:roots}
Let $G$ be a nil-2 group, let $r>0$, let $\mathbf{n}$ be an
$r$-tuple of  positive integers, and let $\mathbf{g}$ be an $r$-tuple
of elements of $G$. Then there exists an ${\cal N}_2$-overgroup $K$ of
$G$ containing an $n_j$-th root for $g_j$ ($1\leq j\leq r$) if and
only for every $r\times r$ array $\{c_{ij}\}$ of integers such that
$n_i c_{ij} = n_j c_{ji}$ for all $i$, $j$, and for all
$y_1,\ldots,y_r\in G$,
\[ \mbox{if}\quad
\forall j\left(y_j\equiv \prod_{i=1}^r g_i^{c_{ij}}
\pmod{G'}\right)\qquad
\mbox{then}\qquad
\prod_{j=1}^r [y_j,g_j]=e.\]
\end{lemmab}

In particular note that we can always adjoin roots to a finite family
of central elements. In fact, we can always adjoin \textit{central}
roots to a finite family of central elements, although the argument is
much simpler: to adjoin a central $n_i$-th root to the central elements
$g_i$, $i=1,\ldots,r$, just let $q_i$ be the order of $g_i$ ($q_i=0$
if $g_i$ is of infinite order), take
\[G\times (Z/n_1q_1Z)\times\cdots\times (Z/n_r q_rZ)\]
and identify the $n_i$-th power of the generator of the $i$-th cyclic
group with~$g_i$.

The following result will be used several times:
\begin{lemmab}
\label{basiceqs}
Let $G\in(m,n)$, $x,y\in G$. Let $K$ be an $(m,n)$-overgroup of $G$,
and assume that for some $q\in\mathbb{Z}$, $r,s\in K$, and $r',s'\in
K^nK'$, we have $r^qr'=x$ and $s^qs'=y$. For any integers $a$, $b$,
$c$, and~$d$, and any $g_1,g_2\in K$,
\[\mbox{if}\qquad\begin{array}{rclc}
g_1^q & \equiv & x^a y^b &\pmod{K^nK'}\\
g_2^q & \equiv & x^{c}y^d &\pmod{K^nK'}
\end{array}\qquad\mbox{ then }[r,s]^{q(c-b)}=[g_1,x][g_2,y].\]
In particular, if such a congruence has a solution with $g_1,g_2\in
G$, then $[r,s]^{q(c-b)}\in G$.
\end{lemmab}
\begin{proof}
Since $xr^{-q}$ and $ys^{-q}$ lie in $K^nK'$, they are central in
$K$. So,
\begin{eqnarray*}
e & = & [g_1r^{-a}s^{-c},xr^{-q}][g_2r^{-b}s^{-d},ys^{-q}]\\
& = & [g_1,x][g_2,y][r^{-a}s^{-c},x][r^{-b}s^{-d},y]\\
&&\qquad [g_1,r^{-q}][g_2,s^{-q}][s^{-c},r^{-q}][r^{-b},s^{-q}]\\
& = &
[g_1,x][g_2,y][r,x^{-a}y^{-b}][s,x^{-c}y^{-d}][g_1^{-q},r][g_2^{-q},s][r,s]^{qb-qc}\\
& = &
[g_1,x][g_2,y][r,g_1^{-q}][s,g_2^{-q}][g_1^{-q},r][g_2^{-q},s][r,s]^{q(b-c)}\\
& = & [g_1,x][g_2,y][r,s]^{q(b-c)}.
\end{eqnarray*}
Therefore, $[r,s]^{q(c-b)}=[g_1,x][g_2,y]$, as desired. If $g_1,g_2\in
G$, then $[g_1,x][g_2,y]\in G$, so $[r,s]^{q(c-b)}\in G$.
\end{proof}

\begin{remark}
Note that the conditions depend only on the equivalence classes of
$x$ and $y$ modulo $K^nK'$; that the conditions are symmetric on
$x$ and~$y$; and that we may restrict $a,b,c,d$ to $0\leq a,b,c,d\leq n-1$.
\end{remark}

In our applications, we will want to adjoin $q$-th roots to either one
or two elements of $G\in (m,n)$, with $q|n$, and staying in $(m,n)$,
so we restrict our statements to that situation. First we need a
lemma:

\begin{lemmab}
\label{lemma:comms}
Let $G\in(m,n)$, and let $g\in Z(G)$ be an element of exponent
$n$. Then there exists an $(m,n)$-overgroup $K$ of~$G$, and elements
$r_1,r_2\in C_K(G)$ such that $g=[r_1,r_2]$. 
\end{lemmab}
\begin{proof}
Let $q$ be the order of $g$ ($q=0$ if $n=0$ and $g$ is of infinite
order). In particular, $q|n$, and also $q|m$. Let $H$ be the group:
\[H = \Bigl\langle r_1, r_2 \Bigm|
r_1^q=r_2^q=[r_1,r_2]^q=[r_1,r_2,r_1]=[r_1,r_2,r_2]=e\Bigr\rangle.\]
Note that $H\in(m,n)$ necessarily, and that $\langle g\rangle\cong
\langle[r_1,r_2]\rangle$.  Now let \[K=G\times_{\langle g\rangle\simeq
\langle[r_1,r_2]\rangle} H.\]
\end{proof}

We will consider the case when $m$ and $n$ are prime powers first, and
deal with the general case later.

\begin{theorem}
\label{addtworootsp}
Let $p$ be a prime, and let $G\in(p^{a+b},p^a)$, with $a\geq 0$ and $b\geq
0$ ($b>0$ if $p=2$). Let $x,y\in G$, and $i$ be an integer with $1\leq
i\leq a$. There exists a group $K\in(p^{a+b},p^a)$, overgroup of $G$, with
$x,y\in K^{p^i}K'$ if and only if both of the following two conditions
hold:
\begin{itemize}
\item[(a)]~$x^{\zeta}=y^{\zeta} = e$, where $\zeta =
{\rm lcm}(p^{a+b-i}, p^a)$.
\item[(b)]~For all $g_1,g_2\in G$,
$\alpha,\beta,\gamma,\delta\in\mathbb{Z}$ with
$\beta\equiv\gamma\pmod{p^{a-i}}$,
\[\mbox{if}\quad\begin{array}{rcl}
g_1^{p^i} & \equiv & x^{\alpha} y^{\beta}\\
g_2^{p^i} & \equiv & x^{\gamma} y^{\delta}
\end{array}
\pmod{G^{p^a}G'},\quad\mbox{then}\quad [g_1,x][g_2,y]=e.\]
\end{itemize}
Moreover, if $i\leq b$, then we may in fact choose $K$ so that $x$ and
$y$ have $p^i$-th roots in~$K$.
\end{theorem}
\begin{remark}
We do the result one prime at a time because the case $i>b$
introduces complications which would be hard to keep track of if we
try to deal with the general case directly. 
\end{remark}

\begin{proof}
For necessity, assume that $K$ is an overgroup of $G$ with
$r^{p^i}r'=x$, $s^{p^i}s'=y$, and $r',s'\in K^{p^a}K'$. If $i>b$, then
$x^{p^a}=r^{p^{a+i}}r'^{p^a}=e$,
since $a+i>a+b$, and similarly with
$y^{p^a}$. If $i\leq b$, then 
\[x^{p^{a+b-i}}=r^{p^{a+b}}r'^{p^{a+b-i}}=e,\]
since $a+b-i\geq a$, and again similarly for $y^{p^{a+b-i}}$, yielding~(a).

For (b), note that since the congruences hold modulo $G^{p^a}G'$, they
also hold modulo $K^{p^a}K'$. We have $p^{a-i}|(\gamma-\beta)$, so 
by Lemma~\ref{basiceqs}
\[[r,s]^{p^i(\gamma-\beta)}=[g_1,x][g_2,y].\]
Since $K\in(p^{a+b},p^a)$, and $p^a|p^i(\gamma-\beta)$, the left hand side must be
trivial, so $[g_1,x][g_2,y]=e$, as desired. 

Moreover, if $i\leq b$, we can adjoin central $p^i$-th roots to $r'$
and $s'$ to get $p^i$-th roots for $x$ and $y$. This gives necessity.

For sufficiency, we will construct $K$, along the lines of Saracino's
proof of Lemma~\ref{lemma:roots}. Let
\[K_0 = G \amalg^{{\cal N}_2}\left(\langle r\rangle \amalg^{(0,p^a)}
\langle s\rangle\right),\]
where $\langle r\rangle$ and $\langle s\rangle$ are infinite cyclic
groups. Let $N$ be the smallest normal subgroup of $K_0$ containing
$xr^{-p^i}$ and $ys^{-p^i}$.

We claim that $N\cap G=\{e\}$. Indeed, a general element of $N$ may be
written as
\[
\prod_{k=1}^{s_1}\left(b_{1k}z_{1k}\right) \left(xr^{-p^i}\right)^{\epsilon_{1k}}
\left(b_{1k}z_{1k}\right)^{-1}\cdot
\prod_{k=1}^{s_2}\left(b_{2k}z_{2k}\right)
\left(ys^{-p^i}\right)^{\epsilon_{2k}}
\left(b_{2k}z_{2k}\right)^{-1}
\]
with $s_i\geq 0$, $\epsilon_{jk}=\pm 1$, $b_{jk}\in G$, and
$z_{jk}=r^{a_{jk1}}s^{a_{jk2}}$. Rewriting,
\begin{eqnarray*}
\lefteqn{\prod_{k=1}^{s_1}
[b_{1k},x]^{\epsilon_{1k}}[b_{1k},r^{-p^i}]^{\epsilon_{1k}}
[z_{1k},x]^{\epsilon_{1k}}[z_{1k},r^{-p^i}]^{\epsilon_{1k}} \cdot
\left(xr^{-p^i}\right)^{t_1}}\\
& \quad & \cdot
\prod_{k=1}^{s_2}
[b_{2k},y]^{\epsilon_{2k}}[b_{2k},s^{-p^i}]^{\epsilon_{2k}}
[z_{2k},y]^{\epsilon_{2k}}[z_{2k},s^{-p^i}]^{\epsilon_{2k}} \cdot
\left(ys^{-p^i}\right)^{t_2}
\end{eqnarray*}
\noindent where $t_j=\sum \epsilon_{jk}$. We do this by first
replacing the conjugate $u^v$ by $[v,u^{-1}]u$, and then expanding the
commutator brackets bilinearly.
We write it as $\alpha\beta\gamma$,
with $\alpha\in G$, $\beta\in \langle r\rangle\amalg^{(0,p^a)}\langle
s \rangle$, and $\gamma\in [G,\langle r\rangle\amalg^{(0,p^a)}\langle
s\rangle]$. Assume that this element lies in $G$, and equals $g\in
G$. Then by uniqueness we have $\alpha=g$, $\beta=\gamma=e$.

On the other hand, $\beta=r^{-p^it_1}s^{-p^it_2}[r,s]^u$ for some
integer $u$, so we must have $t_1=t_2=0$.  Therefore, we have that:
\begin{eqnarray*}
g & = & \prod_{k=1}^{s_1}[b_{1k},x]^{\epsilon_{1k}}\cdot
\prod_{k=2}^{s_2}[b_{2k},y]^{\epsilon_{2k}}.\\
\beta & = & \prod_{k=1}^{s_1}[z_{1k},r^{-p^i}]^{\epsilon_{1k}}\cdot
\prod_{k=1}^{s_2} [z_{2k},s^{-p^i}]^{\epsilon_{2k}}
=  e.\\
\gamma & = & \prod_{k=1}^{s_1}
[b_{1k},r^{-p^i}]^{\epsilon_{1k}}[z_{1k},x]^{\epsilon_{1k}}
\prod_{k=1}^{s_2}
[b_{2k},s^{-p^i}]^{\epsilon_{2k}}[z_{2k},y]^{\epsilon_{2k}} =  e.
\end{eqnarray*}
Since $z_{jk}=r^{a_{jk1}}s^{a_{jk2}}$, we get
\begin{eqnarray}
g & = & \left[\prod_{k=1}^{s_1}b_{1k}^{\epsilon_{1k}},x\right]
\left[\prod_{k=1}^{s_2}
b_{2k}^{\epsilon_{2k}},y\right]. \label{valueg}\\
\beta & = & [r,s]^{-p^i\sum \epsilon_{2k}a_{2k1} + p^i\sum
\epsilon_{1k}a_{1k2}} = e.\label{valuebeta}\\
\gamma & = & \left[\left(\prod_{k=1}^{s_1}
b_{1k}^{\epsilon_{1k}}\right)^{-p^i} x^{-\sum \epsilon_{1k}a_{1k1}}
y^{-\sum \epsilon_{2k}a_{2k1}}, r\right]\nonumber\\
& &\quad \cdot\left[\left(\prod_{k=1}^{s_2} b_{2k}^{\epsilon_{2k}}\right)^{-p^i}
x^{-\sum\epsilon_{1k}a_{1k2}} y^{-\sum\epsilon_{2k}a_{2k2}},
s\right]\label{valuegamma}\\
& = & e.\nonumber
\end{eqnarray}
Let $g_1=\prod b_{1k}^{\epsilon_{1k}}$, $g_2=\prod
b_{2k}^{\epsilon_{2k}}$, and let $c_{mn}=-\sum
\epsilon_{mk}a_{mkn}$. We rewrite (\ref{valuegamma}), and we get
\[ e = \left[g_1^{-p^i}x^{c_{11}}y^{c_{21}}, r\right]
\left[g_2^{-p^i}x^{c_{12}}y^{c_{22}}, s\right].\]
Since the cartesian equals $G^{\rm ab}\otimes(Z\times Z)$, this means
that
\[g_j^{-p^i}x^{c_{1j}}y^{c_{2j}}\in G'\]
so we must have
\[ \begin{array}{rcl}
g_1^{p^i} & \equiv & x^{c_{11}}y^{c_{21}}\\
g_2^{p^i} & \equiv & x^{c_{12}}y^{c_{22}}
\end{array}\pmod{G'}.\]
The congruences also hold modulo $G^{p^a}G'$, and
from~(\ref{valuebeta}), we have that
\[e = \beta = [r,s]^{p^i(c_{21}-c_{12})}\]
so $p^a\,|\,p^i(c_{21}-c_{12})$. Thus, $c_{12}\equiv
c_{21}\pmod{p^{a-i}}$, so by condition~(b) we must have
$[g_1,x][g_2,y]=e$; but since 
$g = [g_1,x][g_2,y]$ by (\ref{valueg}), 
this means that $g=e$. So $G\cap N=\{e\}$, as claimed.

Let $K_1= K_0/N$. Note that $G$ is a subgroup of $K_1$, and that, in
$K_1$, $r^{p^i}=x$, $s^{p^i}=y$. Note as well that $[r,s]^{p^a}=e$ in
$K_0$, and hence also in $K_1$.

If $i\leq b$, consider the subgroup of $K_1$
generated by $G$, $r$, and $s$. We already have that $G\in
(p^{a+b},p^a)$. Note also that:
\begin{eqnarray*}
r^{p^{a+b}} = \left(r^{p^i}\right)^{p^{a+b-i}} & = & x^{p^{a+b-i}} = e
\\
s^{p^{a+b}} = \left(s^{p^i}\right)^{p^{a+b-i}} & = & y^{p^{a+b-i}} = e
\end{eqnarray*}
\noindent by condition~(a); to show that the subgroup $\langle
G,r,s\rangle$ lies in $(p^{a+b},p^a)$, it will now suffice to show
that both $r^{p^a}$ and $s^{p^a}$ are central. First, note that both
$x^{p^{a-i}}$ and $y^{p^{a-i}}$ are central in $G$: for given any
$g\in G$,
\begin{eqnarray*}
\left( g^{p^{a-i}}\right)^{p^i} \equiv x^0y^0\pmod{G^{p^a}G'}\\
e^{p^i} \equiv x^0y^0\pmod{G^{p^a}G'}
\end{eqnarray*}
so by (b), $e=[g^{p^{a-i}},x]=[g,x^{p^{a-i}}]$. A similar calculation
holds for $y^{p^{a-i}}$.  Therefore, $r^{p^a}=x^{p^{a-i}}$ centralizes
$G$; it trivially centralizes $r$, so we just need to know that
$r^{p^a}$ centralizes $s$ to get that it is central in the subgroup in
question. But $[r^{p^a},s]=[r,s]^{p^a}=e$ in $K_0$; so $r^{p^a}$ is central
in $\langle G,r,s\rangle$, and an analogous calculation holds for
$s^{p^a}$. Thus, if we let $K=\langle G,r,s\rangle$, we have $x,y$
have $p^i$-th roots in~$K$, and $K\in(p^{a+b},p^a)$.

Assume now that $i>b$. We still have that both $r^{p^a}$ and
$s^{p^a}$ are central in $\langle G,r,s\rangle$, by the same argument
as above; however, the orders of $r$ and $s$ may be greater than
$p^{a+b}$. 

Let $K_2$ be the result of adjoining to $K_1$ central $p^a$-th roots $t$ and
$v$ to $x^{p^{a-i}}$ and $y^{p^{a-i}}$, respectively. Note that both
$t^{p^i}$ and $v^{p^i}$ are of exponent~$p^a$: 
\[(t^{p^i})^{p^a} = (t^{p^a})^{p^i} = (x^{p^{a-i}})^{p^i} = x^{p^a} =
e\]
by (a), and analogously with $v^{p^i}$. Let $K_3$ be the result of
adjoining to $K_2$ elements $q_1,q_2,q_3,q_4$, of exponent $p^a$, with
$[q_1,q_2]=t^{p^i}$, $[q_3,q_4]=v^{p^i}$, as in
Lemma~\ref{lemma:comms}, so that $q_1$ commutes with all elements of
$K_3$, except $q_2$, etc.

Let $K$ be the subgroup of $K_3$ generated by $G$, $rt^{-1}$,
$sv^{-1}$, and the $q_i$. We claim $K\in(p^{a+b},p^{a})$. Note that
\begin{eqnarray*}
(rt^{-1})^{p^{a+b}} & = & r^{p^{a+b}}t^{-p^{a+b}}\\
& = & x^{p^{a+b-i}} (x^{-p^{a-i}})^{p^b}\\
& = & x^{p^{a+b-i}}x^{-p^{a+b-i}}\\
& = & e,
\end{eqnarray*}
and analogously with $sv^{-1}$. And $(rt^{-1})^{p^a}$ is central, since
$r^{p^a}$ is central here, and so is $t$; the same is true of
$(sv^{-1})^{p^a}$, so
$K\in(p^{a+b},p^a)$. Finally, note that in $K$, we have
\[\left(rt^{-1}\right)^{p^i}[q_1,q_2] = r^{p^i}t^{-p^i}[q_1,q_2] =
r^{p^i} t^{-p^i}t^{p^i} = r^{p^i} = x.\]
and similarly $(sv^{-1})^{p^i}[q_3,q_4]= y$; thus $x,y\in K^{p^i}K'$,
as desired.
\end{proof}

\begin{theorem}
\label{addtworootsgen}
Let $G\in(m,n)$. Let $x,y\in G$, and $q>0$ with
$q|n$. There exists $K\in (m,n)$, overgroup of $G$, with $x,y\in
K^qK'$ if and only if both of the following conditions hold:
\begin{itemize}
\item[(i)]~$x^\zeta=y^{\zeta}=e$, where $\zeta={\rm
lcm}(\frac{m}{q},n)$; and
\item[(ii)]~for all $g_1,g_2\in G$,
$\alpha,\beta,\gamma,\delta\in\mathbb{Z}$ with $\beta\equiv
\gamma\pmod{\frac{n}{q}}$, 
\[\mbox{if}\quad\begin{array}{rcl}
g_1^q & \equiv & x^{\alpha} y^{\beta}\\
g_2^q & \equiv & x^{\gamma} y^{\delta}
\end{array}
\pmod{G^nG'},\qquad\mbox{ then } [g_1,x][g_2,y]=e.\]
\end{itemize}
If $qn|m$, then we may choose $K$ so that $x$ and $y$ both have $q$-th
roots in~$K$.
\end{theorem}
\begin{proof} 
For the case $m>0$, we can proceed one prime at a time by working in
the $p$-parts of $G$. If we simply use the analogous construction
directly when $m=0$, we will always be in the situation analogous to
the case $i\leq b$ in the previous theorem, so everything works
out. Note that condition (i) is trivially true if $m=0$.
\end{proof}

By setting $y=e$, we may obtain the conditions for adjoining a root to
a single element:

\begin{corb}
\label{addonerootgen}
Let $G\in(m,n)$. Let $x\in G$ and $q>0$ be an integer with $q|n$.
There exists a group $K\in(m,n)$, overgroup of $G$, with $x\in K^qK'$
if and only if $x^{\zeta}=e$, where $\zeta={\rm lcm}(\frac{m}{q},n)$, and
for all $g\in G$,
$\alpha\in\mathbb{Z}$ with
$g^q\equiv x^{\alpha}\pmod{G^nG'}$, we have $[g,x]=e$.
If $qn|m$, then we may choose $K$ so that $x$ has a $q$-th root in $K$.
\end{corb}

\subsection*{Controlling the center}

One key ingredient in our embedding result will be constructing
a group $G\in(m,n)$ which has a specified structure
for $G^nG'$, the elements that \textit{must} be central in $G$ and any
$(m,n)$-overgroup. Because we can restrict to finitely generated
groups, we only need to consider finitely generated abelian groups as
the specified structure of $G^nG'$. We begin with the easy case of
$m=0$:

\begin{lemmab}
\label{abtocomms}
Let $A$ be a finitely generated abelian group. Then there exists
$G\in(0,n)$ such that $G^nG' \cong A$.
\end{lemmab}
\begin{proof}
Decompose $A$ into a sum of cyclic groups, $A=\oplus(Z/a_iZ)$. If
$n=0$, then for each $i$ let:
\[G_i = \Bigl\langle u_i,v_i \Bigm|
u_i^{a_i}=v_i^{a_i}=[u_i,v_i]^{a_i}=[u_i,v_i,u_i] = [u_i,v_i,v_i] =
e\Bigr\rangle.\]
Then let $G=\oplus G_i$. It is easy to verify $G$ satisfies the given
condition. If $n>0$, then simply adjoin a central $n$-th root to the
generator of each cyclic subgroup, to obtain an abelian group
$G\in(0,n)$ with $G^n\cong A$.
\end{proof}

The case $m>0$ is slightly more complicated, so we deal with it one
prime at a time:

\begin{lemmab}
\label{cyclictocomm}
Let $p$ be a prime, and let $a,i\geq 0$. There is a group $G\in (p^{a+b},p^a)$ with 
\[G^{p^a}G' = Z(G) \cong Z/p^iZ\]
provided that:
\begin{itemize}
\item[(a)]~if $i\leq a$ and $p$ is odd, then $b\geq 0$.
\item[(b)]~if $i\leq a$ and $p=2$, then $b>0$.
\item[(c)]~if $i>a$, then $b\geq i$.
\end{itemize}
\end{lemmab}
\begin{proof}
If $i\leq a$, let
\[G = \Bigl\langle u,v \Bigm|
u^{p^i}=v^{p^i}=[u,v]^{p^i}=[u,v,u]=[u,v,v]=e\Bigr\rangle.\]
If $i>a$, we let $G$ be the split metacyclic $p$-group:
\[G = \Bigl\langle u, v \Bigm| u^{p^{a+i}}=v^{p^a}=e;\quad
uv=vu^{1+p^i}\Bigr\rangle.\]
It is not hard to verify these definitions will work.
\end{proof}
\begin{corb}
\label{abtopcomm}
Let $A$ be an abelian group of exponent $p^i$, with $i\geq 0$, and $p$ a
prime. Then there exists $G\in(p^{a+b},p^a)$ with
\[G^{p^a}G'=Z(G)\cong A\]
provided $i$, $a$, and $b$ satisfy conditions (a)--(c) of
Lemma~\ref{cyclictocomm}.
\end{corb}
\begin{proof} 
Since $A$ is of bounded exponent, we can decompose it as a sum of
cyclic groups, each of order $p^j$ for some $j$, with $j\leq i$. For each direct
summand, use Lemma~\ref{cyclictocomm}, and then take the direct sum of
the groups so obtained.
\end{proof}

\begin{corb}
\label{boundabtocomm}
Let $A$ be an abelian group of bounded exponent $k>0$. If $m>0$ and $k|{\rm
lcm}(\frac{m}{n},n)$, then there exists $G\in(m,n)$ with
\[Z(G)=G^nG'\cong A.\] 
\end{corb}
\begin{proof}
Decompose $A$ into a sum of cyclic groups of prime power
order and deal with the $p$-parts separately. The condition that
$k|{\rm lcm}(\frac{m}{n},n)$ guarantees that we can apply
Corollary~\ref{abtopcomm} in each case, and we are done.
\end{proof}

\section{Weak and strong amalgamation}

In this section, we will give necessary and sufficient conditions for
weak and strong embeddability of amalgams in $(m,n)$.

We begin by quoting the following result:
\begin{lemmab}[Lemma 4 in~\cite{amalgone}]
\label{constructionlemma}
Let $A,B\in{\cal N}_2$, and let $U$ be a subgroup, $A'\leq U\leq
A$. Assume that $A/U=\langle a_1U\rangle\times\cdots\times \langle
a_nU\rangle$, for $a_i\in A$ and let $n_i$ be the order of $a_i$
modulo $U$ (with $0<n_i<\infty$ for $1\leq i\leq t$, and $n_i=0$
for $t<i\leq n$).  Given morphisms $\varphi\colon U\to B$ and
$\varphi_i\colon\langle a_i\rangle \to B$, with $\varphi_i(a_i)=b_i$,
there exists a homomorphism $\Phi\colon A\to B$ with $\Phi|_U=\varphi$
and $\Phi(a_i)=b_i$ for each $i$, provided that:
\begin{itemize}
\item[(a)]~$\varphi\bigl([g,a_i]\bigr) = \bigl[\varphi(g),b_i\bigr]$
for all $g\in U$, $1\leq i\leq n$.
\item[(b)]~$\varphi\bigl([a_i,a_j]\bigr) = [b_i,b_j]$ for all $1\leq
i,j\leq n$.
\item[(c)]~$\varphi(a_i^{n_i}) = b_i^{n_i}$ for $1\leq i\leq t$.
\end{itemize}
Moreover, $\Phi$ is injective provided that both $\varphi$ and the
map induced by the $\varphi_i$ on
$A/U\to \Phi(A)/\Phi(U)$ are injective.
\end{lemmab}

We will prove our characterization of embeddability by solving a
special case and then reducing.  Most of the ``heavy lifting'' will be
done by the following proposition, whose proof is patterned after the
final step in the Hauptsatz in~\cite{amalgone}:

\begin{propb}
\label{heavylifting}
Let $(A,B;D)$ be an $(m,n)$-amalgam, and assume that $A$ is
finitely generated over $D$. If
\begin{itemize}
\item[(a)]~$A^nA'=D^nD' = B^nB'$; and
\item[(b)]~For all $q>0$ with $q|n$, $a\in A$ and $b\in B$ with
$a^q,b^q\in D$, we have
\[[a^q,b] = [a,b^q];\]
\end{itemize}
then the amalgam $(A,B;D)$ is strongly embeddable in $(m,n)$.
\end{propb}
\begin{remark} It is not hard to verify that we may restrict $q$ to
prime powers.
\end{remark}

\begin{proof}
Say that $A=\langle D,a_1,\ldots,a_k\rangle$, with $a_iD$ of order
$q_i$ modulo $D$ ($q_i=0$ if $n=0$ and $a_i$ is of infinite order
modulo $D$), $q_i|n$. Write:
\begin{eqnarray*}
a_i^{q_i} & = & g_i^{-1} \in D.\\
{}[a_i,a_j] & = & z_{ij}^{-1}\in D^nD'.\\
{}[a_i,d] & = & z_{id}^{-1}\in D^nD'\qquad\mbox{ for each }d\in D.
\end{eqnarray*}

Let $H = B\amalg^{(m,n)}\langle c_1\rangle \amalg^{(m,n)} \langle
c_2\rangle \amalg^{(m,n)}\cdots\amalg^{(m,n)} \langle c_k\rangle$,
where $\langle c_i\rangle$ is a cyclic group of order $m$ (infinite
cyclic if $m=0$). Note that
$H\in(m,n)$.

Let $N$ be the normal subgroup of $H$ generated by the elements
$x_i=c_i^{q_i}g_i$, $[c_i,c_j]z_{ij}$, and $[c_i,d]z_{id}$ for $d\in
D$.

First we claim that $N\cap B=\{e\}$. Note that $[c_i,c_j]z_{ij}$,
$[c_i,d]z_{id}$ are central for each $i,j,d$, and that if $q_i>0$,
$g_i^{n/q_i}$ is central and $g_i^{m/q_i}=e$. If $q_i=0$, then $g_i=e$.

A general element of $N$ can be written as:
\begin{equation}
 h = \prod_{\mu=1}^r
\left(x_{\sigma_{\mu}}^{\delta_{\mu}}\right)^{h_{\mu}} \cdot \!\!\!\!
\prod_{1\leq j<i\leq k}\!\!\!\!\left([c_i,c_j]z_{ij}\right)^{n_{ij}} \cdot
\prod_{i=1}^{k}\prod_{d\in J_i}\left([c_i,d]z_{id}\right)^{n_{id}}
\label{hexpression}
\end{equation}
where $r\geq 0$, $\sigma_{\mu}\in\{1,\ldots,k\}$, $\delta_{\mu}=\pm
1$, $n_{ij},n_{id}\in\mathbb{Z}$, $h_{\mu}\in H$, and each $J_i$ is a
finite set of elements of $D$. We can use a single finite set $J$ by
setting the necessary $n_{id}$ to zero.

We first look at the factor
$\prod_{\mu=1}^{r}\left(x_{\sigma_{\mu}}^{\delta_{\mu}}\right)^{h_{\mu}}$.
Ordering the factors lexicographically, we obtain
\begin{eqnarray}
\prod_{\mu=1}^{r}\left(x_{\sigma_{\mu}}^{\delta_{\mu}}\right)^{h_{\mu}}
& = & \prod_{\mu=1}^{r} x_{\sigma_{\mu}}^{\delta_{\mu}} \left[
x_{\sigma_{\mu}}, h_{\mu}^{\delta_{\mu}}\right] \nonumber\\
& = & \prod_{i=1}^k x_i^{m_i} \cdot \!\!\!\!\prod_{1\leq j<i\leq k}
\!\!\!\![x_i,x_j]^{m_{ij}} \cdot \prod_{i=1}^k [x_i,h'_i]\label{firstfactor}
\end{eqnarray}
where 
\[m_i=\sum_{\sigma_{\mu}=i} \delta_{\mu},\qquad m_{ij}=\!\!\!\!
\mathop{\mathop{\sum_{\sigma_{\mu}=i}}\limits_{\sigma_{\nu}=j}}\limits_{1\leq\mu<\nu\leq
r}\!\!\!\!\delta_{\mu}\delta_{\nu},\qquad h'_i=\prod_{\sigma_{\mu}=i} h_{\mu}^{\delta_{\mu}}.\]
Since $x_i=c_i^{q_i}g_i$, we have
$x_i^{m_i} = \left(c_i^{q_i}g_i\right)^{m_i} =
c_i^{m_iq_i}g_i^{m_i}[g_i,c_i^{q_i}]^{m_i\choose 2}$.
So we have that
\begin{eqnarray}
\lefteqn{\prod_{\mu=1}^{r} \left(x_{\sigma_{\mu}}^{\delta_{mu}}\right)^{h_{\mu}}
 =}\nonumber\\
& &\prod_{i=1}^k c_i^{m_iq_i}g_i^{m_i}[g_i,c_i^{q_i}]^{m_i\choose 2}
\cdot\!\!\!\! \prod_{1\leq j<i\leq k}\!\!\!\![x_i,x_j]^{m_{ij}} \cdot
\prod_{i=1}^{k}[x_i,h'_i].\label{expandfirst}
\end{eqnarray}

Expanding the second factor in (\ref{expandfirst}), we have
\begin{eqnarray*}
\prod_{1\leq j<i\leq k}\!\!\!\! [x_i,x_j]^{m_{ij}} & = & \!\!\!\!\prod_{1\leq j<i\leq
k}\!\!\!\! [c_i^{q_i}g_i,c_j^{q_j},g_j]^{m_{ij}}\\
& = &\!\!\!\! \prod_{1\leq j<i\leq k}\!\!\!\![c_i,c_j]^{q_iq_jm_{ij}}
[c_i,g_j]^{q_im_{ij}}[c_j,g_i]^{-q_jm_{ij}} [g_i,g_j]^{m_{ij}}.
\end{eqnarray*}
If we expand the third factor in (\ref{expandfirst}), and writing 
\[h'_i = \left(\prod_{j=1}^{k}c_j^{\ell_{ij}}\right)b_ih''_i\]
with $\ell_{ij}\in\mathbb{Z}$, $b_i\in B$, $h''_i\in H'$, we have
\begin{eqnarray*}
{}[x_i,h'_i] & = &
\left[c_i^{q_i}g_i,\left(\prod_{j=1}^kc_j^{\ell_{ij}}\right)
b_ih''_i\right] \\
& = &
\left(\prod_{j=1}^k[c_i,c_j]^{q_i\ell_{ij}}\cdot[c_j,g_i]^{-\ell_{ij}}\right)
\cdot [c_i,b_i]^{q_i}[g_i,b_i].
\end{eqnarray*}

Looking back to the second factor in~(\ref{hexpression}), we have
\[\prod_{1\leq j<i\leq k}\!\!\!\!\left([c_i,c_j]z_{ij}\right)^{n_{ij}} = \!\!\!\!\prod_{1\leq
j<i\leq k}\!\!\!\![c_i,c_j]^{n_{ij}}z_{ij}^{n_{ij}}.\]
And the third factor in~(\ref{hexpression}) yields:
\begin{eqnarray*}
\prod_{i=1}^k \prod_{d\in J} \left([c_i,d]z_{id}\right)^{n_{id}} & = &
\prod_{i=1}^k \left[ c_i,\prod_{d\in J} d\,^{n_{id}}\right]\cdot
\prod_{d\in J} z_{id}^{n_{id}}\\
& = & \prod_{i=1}^k[c_i,d_i]z_i,
\end{eqnarray*}
by setting $d_i=\prod d\,^{n_{id}}$ and $z_i=\prod z_{id}^{n_{id}}$. 

Making the substitutions into the expression for $h$ given
in~(\ref{hexpression}), we obtain:
\begin{eqnarray*}
h & = & \left(\prod_{i=1}^{r}
c_i^{m_iq_i}g_i^{m_i}[g_i,c_i^{q_i}]^{m_i\choose 2}\right)\\
& &\qquad \cdot \!\!\!\!\prod_{1\leq j<i\leq k}\!\!\!\![c_i,c_j]^{q_iq_jm_{ij}}
[c_i,g_j]^{q_im_{ij}} [c_j,g_i]^{-q_jm_{ij}} [g_i,g_j]^{m_{ij}}\\
& & \qquad \cdot \left(\prod_{i,j=1}^{k}[c_i,c_j]^{q_i\ell_{ij}}
[c_j,g_i]^{-\ell_{ij}}\right) \cdot
\left( \prod_{i=1}^k [c_i,b_i]^{q_i}[g_i,b_i]\right)\\
& & \qquad \cdot \left(\prod_{1\leq j<i\leq k}\!\!\!\!
[c_i,c_j]^{n_{ij}}z_{ij}^{n_{ij}}\right) \cdot \prod_{i=1}^k
[c_i,d_i]z_i.
\end{eqnarray*}
Assume that $h\in N\cap B$, and we will prove that $h=e$. The maps
${\rm id}\colon B\to B$ and $c_i\mapsto e$ induce a unique map
$\psi\colon H\to B$ by the universal property of the coproduct. Since
$h\in B$, $\psi(h)=h$; on the other hand, all $c_i$ map to $e$, so
\begin{eqnarray*}
h & = & \psi(h)\\
& = & \prod_{i=1}^r g_{i}^{m_i} \cdot \!\!\!\!\prod_{1\leq j<i\leq k}
\!\!\!\![g_i,g_j]^{m_{ij}} \cdot \prod_{i=1}^k[g_i,b_i]\cdot\!\!\!\! \prod_{1\leq
j<i\leq k}\!\!\!\!z_{ij}^{n_{ij}}\cdot \prod_{i=1}^k z_i\\
& = & h_1.
\end{eqnarray*}
Let $X$ be the product of $k$ copies of cyclic groups of order $m$
(infinite cyclic if $m=0$),
and denote the generators of the factors by $u_i$. Mapping $B\to
B\times X$ by mapping identically to the first coordinate, and mapping
$c_i\mapsto u_i$, we get another induced map $\psi_2\colon H\to
B\times X$. On the one hand, $h$ maps to itself, but on the other hand
we have:
\[h = \psi_2(h) = \left(\prod_{i=1}^r u_i^{m_iq_i}\right) \cdot h_1.\]
Since we already know that $h=h_1$, this yields that the first factor
is trivial, so $u_i^{m_i q_i}=e$ for each $i$; thus, $m\,|\,m_i q_i$ for
each $i$. Therefore, we also have that $n\,|\,q_i{m_i\choose 2}$; hence,
since $H\in(m,n)$, 
\[ [g_i,c_i^{q_i}]^{m_i\choose 2} =  c_i^{m_iq_i} = e.\]
In addition, since $m|m_i q_i$, and $q_i|m$, we have that
either
$\frac{m}{q_i}|m_i$, or $m=0$ and $q_im_i=0$. If $q_i=0$, then $x_i$
is trivial anyway, and if $q_i\not=0$, then $m_i=0$. So
several factors in the original expression are actually trivial, and we have:
\[ h = h_1 = \!\!\!\!\prod_{1\leq j<i\leq k}\!\!\!\![g_i,g_j]^{m_{ij}} \cdot
\prod_{i=1}^k[g_i,b_i]\cdot\!\!\!\! \prod_{1\leq j<i\leq k}\!\!\!\! z_{ij}^{n_{ij}}
\cdot \prod_{i=1}^k z_i.\]
In particular, all the other terms in the original expression for $h$
are actually trivial, so
\begin{eqnarray*} 
e & = & h_2\\
& = &\!\!\!\! \prod_{1\leq j<i\leq k}\!\!\!\![c_i,c_j]^{q_iq_jm_{ij}}
[c_i,g_j]^{q_im_j} [c_j,g_i]^{-q_jm_{ij}}\\
& & \qquad\cdot \prod_{i,j=1}^k[c_i,c_j]^{q_i\ell_{ij}}
[c_j,g_i]^{-\ell_{ij}} \cdot \prod_{i=1}^k[c_i,b_i]^{q_i}\\
& & \qquad \cdot \!\!\!\!\prod_{1\leq j<i\leq k}\!\!\!\![c_i,c_j]^{n_{ij}} \cdot
\prod_{i=1}^k [c_i,d_i]\\
& = &\!\!\!\!\prod_{1\leq j<i\leq k}\!\!\!\! [c_i,c_j]^{(q_iq_jm_{ij} +
q_i\ell_{ij}-q_j\ell_{ji}+n_{ij})}\\
& & \qquad\cdot\!\!\!\! \prod_{1\leq j<i\leq k}\!\!\!\! [c_i,g_j]^{q_im_{ij}}
[c_j,g_i]^{-q_jm_{ij}}\\
& & \qquad\cdot \prod_{i,j=1}^k[c_j,g_i]^{-\ell_{ij}}\cdot
\prod_{i=1}^k\left([c_i,b_i]^{q_i}[c_i,d_i]\right).\\
\end{eqnarray*}
Fix $j<i$. Pick $X\in(m,n)$ with a commutator of order $n$ (of
infinite order if $n=0$), for
example the relatively free group of rank 2 generated by $u$ and
$v$. Map $c_i\mapsto u$, $c_j\mapsto v$, and $c_{\ell}\mapsto e$ for
$\ell\not=i,j$, and map $B\mapsto e$; this induces a map $\sigma\colon
H\to X$, which necessarily maps $h_2$ to $e$; it also maps $h_2$ to
\[ e = \sigma(h_2) = [u,v]^{(q_iq_jm_{ij} + q_i\ell_{ij} - q_j\ell_{ji} +
n_{ij})}.\]
Therefore,
\begin{equation}
n \,|\, q_iq_jm_{ij} + q_i\ell_{ij}-q_j\ell_{ji}+n_{ij}
\label{longexpdiv}
\end{equation}
for each choice of $j<i$, so each of these commutators in $h_2$ are
trivial. 

Let $\langle c\rangle$ be a cyclic group of order $n$ (infinite cyclic
if $n=0$) and consider the group $B\amalg^{(m,n)}\langle
c\rangle$. Fix $i$, and map $B$ to itself via the canonical inclusion,
send $c_i$ to $c$, and $c_j\mapsto e$ for $j\not=i$. This induces a map
$H\to B\amalg^{(m,n)}\langle c\rangle$. The image of $h_2$ is again
trivial, and also equal to
\begin{eqnarray*}
e & = & \left(\prod_{1\leq j<i} [c,g_j]^{q_im_{ij}}\cdot\!\! \prod_{i<j\leq k}
\!\!\![c,g_j]^{-q_im_{ji}} \cdot \prod_{j=1}^k [c,g_j]^{-\ell_{ji}}\right) \cdot
[c,b_i]^{q_i}\cdot [c,d_i]\\
& = & \left[ c, \left(\prod_{1\leq j<i} g_j^{q_im_{ij}} \cdot \prod_{i<j\leq
k} g_j^{-q_im_{ji}} \cdot \prod_{i=1}^k g_{j}^{-\ell_{ji}}\right)
b_i^{q_i}d_i\right].
\end{eqnarray*}

By Corollary~\ref{whatcommutes}, we must have
\[ \left(\prod_{1\leq j<i} g_j^{q_im_{ij}} \cdot \prod_{i<j\leq
k} g_j^{-q_im_{ji}} \cdot \prod_{i=1}^k g_{j}^{-\ell_{ji}}\right)
b_i^{q_i}d_i \in B^nB'=D^nD'.\]
The $g_i$ all lie in $D$, as does $d_i$; so we must have $b_i^{q_i}\in
D$. We know that $a_i^{q_i}\in D$, and moreover, that $q_i|n$, so
by condition (b) in the statement of the Proposition,
$[a_i^{q_i},b_i]=[a_i,b_i^{q_i}]\in D$. 
Therefore, 
\[ [a_i,b_i^{q_i}] = [a_i^{q_i},b_i]=[g_i^{-1},b_i].\]

Looking at $A$, we must have (since the second factor lies in $D^nD'$,
it makes sense to calculate this in $A$):
\begin{eqnarray*}
e & = & \left[ a_i, \left(\prod_{1\leq j<i} g_j^{q_im_{ij}} \cdot \prod_{i<j\leq
k} g_j^{-q_im_{ji}} \cdot \prod_{i=1}^k g_{j}^{-\ell_{ji}}\right)
b_i^{q_i}d_i\right]\\
& = &\left[ a_i, \left(\prod_{1\leq j<i} a_j^{-q_jq_im_{ij}} \cdot \prod_{i<j\leq
k} a_j^{q_jq_im_{ji}} \cdot \prod_{i=1}^k a_{j}^{q_j\ell_{ji}}\right)
b_i^{q_i}d_i\right]\\
& = & \prod_{1\leq j<i} [a_i,a_j]^{-q_iq_jm_{ij}} \cdot \prod_{i<j\leq
k} [a_i,a_j]^{q_iq_jm_{ji}} \cdot
\prod_{j=1}^k[a_i,a_j]^{q_j\ell_{ji}}\\
& & \qquad \cdot [a_i,b_i^{q_i}][a_i,d_i].
\end{eqnarray*}
Taking the product over all $i=1,\ldots,k$, we have
\begin{eqnarray*}
e & = & \!\!\!\!\prod_{1\leq j<i\leq k}\!\!\!\! [a_i,a_j]^{(-q_iq_jm_{ij}-q_iq_jm_{ij}
+ q_j\ell_{ji} -q_i\ell_{ij})}
\cdot \prod_{i=1}^k [a_i,b_i^{q_i}][a_i,d_i]\\
& = & \!\!\!\!\prod_{1\leq j<i\leq
k}\!\!\!\![a_i,a_j]^{(-2q_iq_jm_{ij}+q_j\ell_{ji}-q_i\ell_{ij})}
\cdot \prod_{i=1}^k[g_i^{-1},b_i][a_i,d_i].
\end{eqnarray*}
From (\ref{longexpdiv}), we have that
\[n_{ij}\equiv -q_iq_jm_{ij} - q_i\ell_{ij} + q_j\ell_{ji} \pmod{n}\]
so $n_{ij}-q_iq_jm_{ij} \equiv -2q_iq_jm_{ij} - q_i\ell_{ij} +
q_j\ell_{ji} \pmod{n}$.
Using this in the expression above, we have
\[ e = \!\!\!\!\prod_{1\leq j<i\leq k}\!\!\!\! [a_i,a_j]^{n_{ij}-q_iq_jm_{ij}} \cdot
\prod_{i=1}^k[g_i^{-1},b_i][a_i,d_i].\]
We also have that:
\begin{eqnarray*}
{}[a_i,a_j]^{n_{ij}-q_iq_jm_{ij}} & = &
[a_i,a_j]^{-q_iq_jm_{ij}}[a_i,a_j]^{n_{ij}}\\
& = & [a_i^{q_i},a_j^{q_j}]^{-m_{ij}} z_{ij}^{-n_{ij}}\\
& = & [g_i^{-1},g_j^{-1}]^{-m_{ij}}z_{ij}^{-n_{ij}}\\
& = & [g_i,g_j]^{-m_{ij}}z_{ij}^{-n_{ij}}
\end{eqnarray*}
and that
\begin{eqnarray*}
{}[a_i,d_i] & = & \left[a_i,\prod_{d\in J} d\,{}^{n_{id}}\right]\\
& = & \prod_{d\in J} [a_i,d]^{n_{id}}
 =  \prod_{d\in J} z_{id}^{-n_{id}}
 =  z_i^{-1}.
\end{eqnarray*}

Therefore,
\begin{eqnarray*}
e & = & \!\!\!\!\prod_{1\leq j<i\leq k}\!\!\!\!
[g_i,g_j]^{-m_{ij}}z_{ij}^{-n_{ij}}\cdot
\prod_{i=1}^k\left([g_i,b_i]z_{i}\right)^{-1}\\
& = & \left( \prod_{1\leq j<i\leq k}\!\!\!\! [g_i,g_j]^{m_{ij}}z_{ij}^{n_{ij}}
\cdot \prod_{i=1}^k[g_i,b_i]z_i\right)^{-1}.
\end{eqnarray*}
Since we already have that
\[ h = \!\!\!\!\prod_{1\leq j<i\leq k}\!\!\!\! [g_i,g_j]^{m_{ij}}z_{ij}^{n_{ij}} \cdot
\prod_{i=1}^k [g_i,b_i]z_i\]
this gives that $h=e$, as claimed. Therefore, $N\cap B=\{e\}$. 

Let $G=H/N\in (m,n)$. Note that, since $N\cap B=\{e\}$, $B$ embeds
into~$G$. What we want to do now is to embed $A$ into $G$, so that
$D\subseteq A\cap B$ in~$G$. Since $A'\subseteq A^nA'=D^nD'\subseteq
D\subseteq A$, and $A$ is generated by $a_1,\ldots,a_k$ modulo $D$,
we may use Lemma~\ref{constructionlemma}. Map $\varphi\colon D\to G$
via the inclusion into $B$; define $\varphi_i\colon \langle
a_i\rangle\to G$ by $\varphi_i(a_i)=\overline{c_i}$. Note that since
$\overline{c_i}^{q_i} = g_i^{-1} = a_i^{q_i}$, the $\varphi_i$ are
well defined.  We need to check
the conditions of Lemma~\ref{constructionlemma} to ensure a map 
$A\to G$ is induced, and that this map is injective.

First, we need to verify that $\varphi([d,a_i])=[\varphi(d),\overline{c_i}]$. 
As $[a_i,d]=z_{id}^{-1}$, and $[c_i,d]z_{id}\in N$, the
image of $[d,c_i]$ in $G$ is equal to $z_{id}$, as desired. Second, we
need to check that $\varphi([a_i,a_j])=\overline{[c_i,c_j]}$. But again,
$[c_i,c_j]z_{ij}$ lies in $N$, and $[a_i,a_j]=z_{ij}^{-1}$. Finally,
to verify that $\varphi(a_i^{q_i})=\overline{c_i}^{q_i}$, which we
already noted. Thus, $\varphi$ and the $\varphi_i$ induce a morphism
$\Phi\colon A\to G$. Since $\Phi|_D=\varphi$, $D$ is contained in $B\cap
\Phi(A)$. Thus, in order to prove that $G$ realizes a (weak) embedding
of $(A,B;D)$, we only need to show that $\Phi$ is injective.

Since $\varphi$ is the identity, it is injective. The maps $c_i\mapsto
a_iD$ and $B\mapsto e$ induce a morphism $H\to A/D$. The kernel
contains $B$, and also contains~$N$, so it factors through
$H/BN$. This gives an induced map
\[\frac{G}{B}\to \frac{A}{D}\]
by $\overline{c_i}B\mapsto a_iD$. Note that $G'\subseteq B$, so we may
look at $\Phi(A)B/B$, which is a subgroup of $G/B$. So we have
\[\frac{\Phi(A)}{\Phi(A)\cap B} \cong \frac{\Phi(A)B}{B}\to
\frac{G}{B}\to\frac{A}{D}\to \frac{\Phi(A)}{\Phi(D)}.\]
Also, $\Phi(A)/\Phi(D)$ maps to $G/B$ by
\[\Phi(a_i)\Phi(D) \to \Phi(a_i)\left(\Phi(A)\cap B\right) \to
\Phi(a_i)B=\overline{c_i}B.\]
Since $G/B$ maps again to $A/D$, we get an endomorphism 
\[\frac{A}{D}\to\frac{\Phi(A)}{\Phi(D)}\to\frac{A}{D}\]
which is the identity on each $\overline{a_i}$. So the first map, which is the
map induced by $\Phi$, must be injective, and together with the
injectivity of $\varphi$ we conclude from
Lemma~\ref{constructionlemma} that $\Phi$ is injective. Thus, $G$
realizes a (weak) embedding of the amalgam.

Finally, we strengthen this embedding into a strong embedding. Since
$A'\subseteq D^nD'\subseteq D$, we have $D\triangleleft A$. Let $K=G\times
(A/D)$. We construct embeddings $A\to K$ by mapping to $G$ via $\Phi$
and to $A/D$ via the canonical projection; and $B\to K$ by mapping to
$G$ via the inclusion and to $A/D$ via the zero map. Now assume that
$a\in A$, $b\in B$ map to the same element of $K$. Then, $aD= D$,
so $a\in D$, proving that their intersection in $K$ is contained in
$D$. Since the embeddings into both coordinates agree on $D$, their
intersection contains at least $D$; thus $K$ realizes the
\textit{strong} embedding of $(A,B;D)$, giving the result.
\end{proof}

\begin{remark}
\label{normalmeansstrong}
Note that the argument in the final paragraph only uses
the fact that $D\triangleleft A$ and that $(A,B;D)$ is weakly
embeddable. In this situation, we can always find a strong embedding
using the argument given.
\end{remark}

We can now give the strong amalgamation criterion:

\begin{theorem}[cf. Satz~3 \cite{amalgtwo}]
\label{strongamalgs}
Let $(A,B;D)$ be an $(m,n)$-amalgam. The amalgam is strongly
embeddable in $(m,n)$ if and only if the following two conditions are
satisfied:
\begin{itemize}
\item[(a)]~$D\cap A^nA' \subseteq Z(B)$ and $D\cap B^nB'\subseteq
Z(A)$.
\item[(b)]~For every $q>0$ with $q|n$, $a\in A$, $a'\in A^nA'$, $b\in
B$, $b'\in B^nB'$ with $a^qa',b^qb'\in D$, we have
\[ \left[a^qa',b\right] = \left[a,b^qb'\right] \in D.\]
\end{itemize}
\end{theorem}

\begin{remark} Again, it is not hard to verify that we may restrict
$q$ to prime powers.
\end{remark}

\begin{proof}
First we prove necessity. Suppose that $K\in(m,n)$ contains $A$ and
$B$, and $A\cap B=D$ in $K$. Any element in $A^nA'$ is central in $K$,
so $A^nA'\cap D\subseteq Z(K)\subseteq C_K(B)$, giving (a) by
symmetry. For (b), note that both $a'$ and $b'$ are central in $K$ and
commutator brackets act bilinearly, so we have
\[ [a^qa',b] = [a^q,b] = [a,b]^q\]
and this equals $[a,b^qb']$ by symmetry. Since $[a^qa',b]\in B$ and
$[a,b^qb']\in A$, their common value lies in $A\cap B$, hence in $D$.

For sufficiency, note that both conditions are inherited to any
finitely generated subamalgam 
$\left(\langle a_1,\ldots,a_n\rangle,\langle b_1,\ldots,b_m\rangle;
\langle a_1,\ldots,a_n\rangle\cap\langle
b_1,\ldots,b_m\rangle\right)$,
so we may assume that both $A$ and $B$ are finitely generated (and
finitely generated over $D$) by Proposition~\ref{fgamalgams}. 

Let $A_0 = A\times_{D\cap B^nB'} B^nB'$ and
$B_0=A^nA'\times_{A^nA'\cap D} B$ be central amalgams. By
Lemma~\ref{korollartwo}, the subgroups $\langle D,A^nA',B^nB'\rangle$
of $A_0$ and $B_0$ are isomorphic. 

Note that if $m>0$, then $A^n$ is of exponent $\frac{m}{n}$, so
$A^nA'$ is of exponent ${\rm lcm}(\frac{m}{n},n)$, and the same with
$B^nB'$. By Lemma~\ref{abtocomms} (for $m=0$) or
Corollary~\ref{boundabtocomm} (for $m>0$), there exists $U\in(m,n)$
such that $U^nU'\cong \langle A^nA',B^nB'\rangle$.

Let $X=A_0\times_{\langle A^nA',B^nB'\rangle} U$ and
$Y=B_0\times_{\langle A^nA',B^nB'\rangle} U$. Then 
\[\langle
A^nA',B^nB'\rangle = X^nX' = Y^nY'\] since $U^nU'\subseteq X^nX'$,
$A_0^nA'\subseteq\langle
A^nA',B^nB'\rangle$, and $[A_0,U]=\{e\}$; and analogous for
$Y^nY'$. Let $E=\langle D,U\rangle$ in both $X$ and $Y$. We want to
apply Proposition~\ref{heavylifting} to $(X,Y;E)$. Considering $E$ as
a subgroup of either $X$ or $Y$, by Lemma~\ref{korollarone} we have
that:
\[ E = \langle D,U\rangle \cong D\times_{D\cap U} U.\]
Again, it is straightforward to verify that
$E^nE'=X^nX'=Y^nY'=\langle A^nA',B^nB'\rangle$.
This gives condition (a) in Proposition~\ref{heavylifting}. To test
condition (b), let $q>0$, $q|n$, $x\in X$, $y\in Y$ with $x^q,y^q\in
E$. We want to show that $[x,y^q]=[x^q,y]$.

By construction of $X$, we may write $x=ab'u$ with $a\in A$, $b'\in B^nB'$,
and $u\in U$. Since all three commute with one another,
$x^q=a^qb'^qu^q\in E$. Since $U,B^nB'\subseteq E$, this says that
$a^q\in E$. We know $E=D\times_{D\cap U} U$, and $a\in A$, $a^q\in
\langle D,U\rangle\cap A=\langle D,A^nA'\rangle$, so we may write
\[a^q = da_2;\qquad d\in D,\quad a_2\in A^nA'.\]
Therefore, $x^q=da_2b'^qu^q$. Note that $a^qa_2^{-1}=d\in D$, with
$a_2^{-1}\in A^nA'$. Symmetrically, we may write $y=ba'v$ with $b\in
B$, $a'\in A^nA'$, $v\in U$, and $y^q=d'b_2a'^qv^q$ with $d'\in D$,
$b_2\in B^nB'$, and $b^qb_2^{-1}=d'\in D$. 

Since $A$, $B$, and $D$ satisfy (b) in the statement, we must have that
\[{}[a^qa_2^{-1},b] = [a,b^qb_2^{-1}] \in D.\]

Now, in $X$ we have:
\begin{eqnarray*}
\left[x,y^q\right] & = & \left[ ab'u,b^qa'^qv^q\right]\\
& = & [ab'u,b^qa'^q][ab'u,v^q]\\
& = & [ab'u,b^qa'^q][u,v^q]\\
& = & [ab',b^qa'^q][u,b^qa'^q][u,v]^q\\
& = & [ab',b^qa'^q][u,v]^q\\
& = & [a,b^q][u,v]^q\\
& = & [a,b^qb_2^{-1}][u,v]^q.
\end{eqnarray*}
Symmetrically, in $Y$ we have
\[\left[x^q,y\right] = \left[a^qa_2^{-1},b\right][u,v]^q.\]
Since $[a^qa_2^{-1},b]=[a,b^2b_2^{-1}]$, we conclude that
$[x,y^q]=[x^q,y]$, as claimed.

Thus $(X,Y;E)$ satisfies the conditions of
Proposition~\ref{heavylifting}, so there exists a $K\in (m,n)$ which
realizes a strong embedding of $(X,Y;E)$. We claim that this
embedding, via the inclusions of $A$ and $B$ into $X$ and $Y$,
respectively, yields a strong embedding of $(A,B;D)$. Since $X\cap
Y=E$ in $K$, we must have that $A\cap B = (A\cap E)\cap (B\cap
E)$. Since $E$ is itself a strong amalgam of $D$ and $U$, $A\cap E=D$
(since $U\cap A_0 = \langle A^nA',B^nB'\rangle$ in $X$, and $A\cap
B^nB'=D\cap B^nB'$ in $A_0$). Symmetrically, $B\cap E=D$. So $A\cap
B=D$ in $K$, and therefore, $K$ provides a strong embedding for
$(A,B;D)$, as desired.
\end{proof}

\begin{remark}
It is also worth noting that the conditions given are independent
of~$m$. So if $(A,B;D)$ is an amalgam of $(m,n)$ groups, and
$(m,n)\subseteq(m',n)$ for some $m'$, then the amalgam is strongly embeddable
in $(m',n)$ if and only if it is also strongly embeddable in
$(m,n)$. Note that the analogous remark with $(m,n)\subseteq (m,n')$
is not true, as shown in Example~\ref{guidingex}.
\end{remark}

The conditions for weak embeddability are more complicated. We give
them in the next result:

\begin{theorem}[cf. Hauptsatz in~\cite{amalgone}]
\label{weakamalgs}
Let $(A,B;D)$ be an amalgam of $(m,n)$-groups. The amalgam is weakly
embeddable in $(m,n)$ if and only if the following two conditions are
satisfied:
\begin{itemize}
\item[(1)]~$A^nA'\cap D\subseteq Z(B)$ and $B^nB'\cap D\subseteq
Z(A)$.
\item[(2)]~ For every $k>0$, $q_i>0$ with $q_i|n$, $a_i\in A$,
$a'_i\in A^nA'$ with $a_i^{q_i}a'_i\in D$, $b_i\in B$, $b'_i\in B^nB'$
with $b_i^{q_i}b'_i\in D$ for $i=1,\ldots,k$, and for each $d\in D$,
\[ \prod_{i=1}^k \left[a_i,b_i^{q_i}b'_i\right] = d
\Longleftrightarrow \prod_{i=1}^k \left[a_i^{q_i}a'_i,b_i\right] =
d.\]
\end{itemize}
\end{theorem}
\begin{proof}
Necessity follows along the same lines as the necessity in
Theorem~\ref{strongamalgs}. Simply note that for each $i$, in any
$(m,n)$-group realizing the amalgam, we have
$\left[a_i^{q_i}a'_i,b_i\right] = \left[a_i,b_i\right]^{q_i} = 
\left[ a_i,b_i^{q_i}b'_i\right]$.

For sufficiency, we will construct an overgroup $D_0$ of $D$, which is
a subgroup of both
$A$ and $B$, and which satisfies conditions (a) and (b) in
Theorem~\ref{strongamalgs}. That will yield that $(A,B;D_0)$ can be
strongly embedded, and thus that $(A,B;D)$ can be weakly embedded
using the same group that realizes $(A,B;D_0)$.

Let $D^A$ be the normal closure of $D$ in $A$, and analogous for
$D^B$. Let
\[\mathcal{S} = \Bigl\{ E\leq A,B\Bigm| E\subseteq D^A,D^B;\ 
D\subseteq E;\ (A,B;E)\mbox{ satisfies (1) and (2)}\Bigr\}.\]
Trivially, $D\in\mathcal{S}$. Order $\mathcal{S}$ by inclusion; the
union of any chain in $\mathcal{S}$ again lies in $\mathcal{S}$, as
there will only be finitely many indices involved in checking
condition~(2). By Zorn's Lemma, $\mathcal{S}$ has maximal
elements. Let $D_0$ be a maximal element of $\mathcal{S}$. We claim
that $D_0$ satisfies conditions (a) and (b) in
Theorem~\ref{strongamalgs}.

Condition (a) there is the same as condition (1) here, so $D_0$
satisfies it by virtue of lying in $\mathcal{S}$.
To show that $D_0$ satisfies condition (b) in
Theorem~\ref{strongamalgs}, let $q>0$ with $q|n$, $a\in A$, $a'\in
A^nA'$, $b\in B$, $b'\in B^nB'$ with $a^qa',b^qb'\in D_0$. We want to
prove that $[a^qa',b]=[a,b^qb']\in D_0$. 

Let $b_0=[a^qa',b]$ and $a_0=[a,b^qb']$. Then $a_0^r\in D_0$ if and
only if $b_0^r\in D_0$ for each $r$, and if so $a_0^r=b_0^r$. Indeed, this
follows from the fact that $D_0$ satisfies condition (2): set $k=r$,
and all terms equal to $a_0$ or $b_0$. So we conclude that
$D_0\cap\langle a_0\rangle = \langle d \rangle = \langle a_0^r\rangle$
and
$D_0\cap\langle b_0\rangle = \langle d \rangle = \langle b_0^r\rangle$
with the same $d$ and the same $r$. 

Since $a_0\in Z(A)$, and $b_0\in Z(B)$, we have that
\[\langle D_0,a_0\rangle \cong D_0\times_{\langle d\rangle} \langle
a_0\rangle; \qquad \langle D_0,b_0\rangle \cong D_0\times_{\langle
d\rangle} \langle b_0\rangle.\]
Since $\langle a_0\rangle \cong \langle b_0\rangle$, $\langle
D_0,a_0\rangle\cong \langle D_0,b_0\rangle$; we denote the common
subgroup of $A$ and $B$ by $D_1$. Then $D_1\leq D_0^A=D^A$,
$D_1\leq D_0^B=D^B$, and we claim that in fact $D_1\in\mathcal{S}$.

If $x\in D_1\cap A^nA'$, then $x=da_0^r$ for some $d\in D_0$,
$r\in\mathbb{Z}$. Since $a_0\in A^nA'$, we also have $d\in A^nA'$, and
by (1) for $D_0$, $d\in Z(B)$; so $x\simeq db_0^r$ lies in
$Z(B)$. Symmetrically for $y\in D_1\cap B^nB'$, yielding that $D_1$
satisfies (1). For (2), let $k>0$, $q_i|n$, $q_i>0$, $a_i\in A$,
$a'_i\in A^nA'$, $b_i\in B$, $b'_i\in B^nB'$ with
$a_i^{q_i}a'_i,b_i^{q_i}b'_i\in D_1$, and let $d\in D_1$. We want to
prove that
\[ \prod_{i=1}^k \left[a_i^{q_i}a'_i,b_i\right] = d
\Longleftrightarrow \prod_{i=1}^k\left[a_i,b_i^{q_i}b'_i\right] = d.\]
Write $d=d'a_0^{\ell} = d'b_0^{\ell}$, for some $d'\in D_0$,
$\ell\in\mathbb{Z}$. In addition, we may write
$a_i^{q_i}a'_i=d_ib_0^{\ell_i}$, so
$a_i^{q_i}\left(a'_ia_0^{-\ell_i}\right)=d_i\in D_0$; so we have
\[\prod_{i=1}^k [a_i^{q_i}a'_i,b_i] = d \Longleftrightarrow
\left(\prod_{i=1}^k [a_i^{q_i}a'_i,b_i]\right)b_0^{-\ell_i}=d'.\]
Since $b_0=[a^qa',b]$, we may apply (2) to $D_0$ to obtain that
\[\left(\prod_{i=1}^k [a_i^{q_i}a'_i,b_i]\right)b_0^{-\ell_i} = d'
\Longleftrightarrow \left(\prod_{i=1}^k [a_i,b_i^{q_i}b'_i]\right)a_0^{-\ell_i} =
d',\]
which means that
\[\prod_{i=1}^k
[a_i^{q_i}a'_i,b_i]=d'b_0^{\ell_i}=d\Longleftrightarrow
\prod_{i=1}^k[a_i,b_i^{q_i}b'_i]=d'a_0^{\ell_i}=d.\]
Therefore, $D_1$ satisfies (2). 

Thus, $D_1\in\mathcal{S}$, $D_0\subseteq D_1$, so by maximality of
$D_0$, $D_1=D_0$. Therefore, $a_0\simeq b_0\in D_0$. Which in turn
shows that $D_0$ satisfies (b) from Theorem~\ref{strongamalgs}. Thus
$(A,B;D_0)$ is strongly embeddable, and therefore $(A,B;D)$ is weakly
embeddable.
\end{proof}
\begin{remark}
Taking a maximal $D_0$ may not be the most efficient amalgam available
(i.e.~the amalgam with the smallest core containing $D$ which is
strongly embeddable). But note that if two common subgroups $D_1$ and
$D_2$ containing $D$ satisfy both (a) and (b) from
Theorem~\ref{strongamalgs}, then so does $D_1\cap D_2$; so in order to
obtain the most efficient weak embedding for $(A,B;D)$, we could take
the intersection of all subgroups which contain $D$ and satisfy (a)
and (b) from Theorem~\ref{strongamalgs}; this is the
smallest common subgroup containing $D$ over which the amalgam is
strongly embeddable.
\end{remark}

\begin{remark} Again, if $(A,B;D)$ is an amalgam of $(m,n)$ groups,
and $(m,n)\subseteq(m',n)$ for some $m'$, then the amalgam is weakly
embeddable in $(m,n)$ if and only if it is weakly embeddable in
$(m',n)$, since conditions (1) and (2) in Theorem~\ref{weakamalgs}
do not depend on~$m$.
\end{remark}

\begin{example} We now see what happened in
Example~\ref{guidingex}. The amalgam is strongly embeddable in
$(0,0)=\mathcal{N}_2$; but for $(4,2)$, it fails to satisfy conditions
(a) from Theorem~\ref{strongamalgs} and (1) from
Theorem~\ref{weakamalgs}: $A^2A'\cap D=D$, but $D$ is not central in
$B$. So the amalgam cannot be embedded in $(4,2)$, nor in $(m,2)$ for
any $m$ multiple of $4$. However, if $4|n$, then $A^nA'\cap D=\{e\}$,
so condition (1) is satisfied if we work in $(m,n)$ for any $n$ with
$4|n$.  More about this in Section~4.
\end{example}

\section{Dominions and Special Amalgams}

Aside from measuring the gap between weak and strong embeddability,
special amalgams also provide a link between the study of amalgams and
the study of dominions. Recall that Isbell \cite{isbellone} defines
for a category ${\cal C}$ of algebras (in the sense of Universal
Algebra) of a fixed type $\Omega$, an algebra $A\in {\cal C}$ and
subalgebra $B$ of $A$, the \textit{dominion of $B$ in $A$} (in the
category $\cal C$) to be the intersection of all equalizers containing
$B$. Explicitly,
\[{\rm dom}_A^{{\cal C}}(B) = \Bigl\{a\in A \,\Bigm|\, \forall f,g\colon A\to C,
\ {\rm if}\ f|_B=g|_B\ {\rm then}\ f(a)=g(a)\Bigr\}\] where $C$
ranges over all algebras $C\in {\cal C}$, and $f,g$ are morphisms. The
connection between amalgams and dominions when working in a variety
${\cal V}$ is
given precisely by the special amalgams: if we let
$A'$ be an isomorphic copy of $A$, and $M=A\amalg^{{\cal V}}_B A'$,
we have that
\[{\rm dom}_A^{{\cal V}}(B) = A\cap A'\subseteq M,\]
where we have identified $B$ with its common image in $A$ and $A'$. In
other words, ${\rm dom}_A^{{\cal V}}(B)$ is the smallest subalgebra
$D$ of~$A$ such that $B\subseteq D$ and the special amalgam $(A,A;D)$ is
strongly embeddable. If ${\rm dom}_A^{{\cal V}}(B)=B$, we say that the
dominion of~$B$ is ``trivial'' (meaning it is as small as possible),
and we say it is ``nontrivial'' otherwise.

In general, ${\rm dom}_A^{{\cal C}}(-)$ is a closure operator on the
lattice of subalgebras of~$A$. If we are working in a variety of
groups, then normal subgroups have trivial dominions, the dominion
construction respects finite direct products (that is, if $H_1<G_1$
and $H_2<G_2$, then the dominion of $H_1\times H_2$ in $G_1\times G_2$
is the product of the dominions of $H_1$ in $G_1$ and of $H_2$ in
$G_2$), and also respects quotients: if $N\triangleleft G$ is
contained in $H$, then
\[{\rm dom}_{G/N}^{\cal V}(H/N) = \left({\rm dom}_G^{\cal
V}(H)\right)\bigm/ N.\]
In addition, if the class $\mathcal{C}$ is contained in the class
$\mathcal{C}'$, then 
\[{\rm dom}_A^{{\cal C}'}(B) \subseteq {\rm dom}_A^{{\cal C}}(B)\]
for every $A\in\mathcal{C}$ and $B\leq A$.  For a proof of these
assertions we direct the reader to \cite{nildoms}.  Note that if
$G\in(m,n)$ and $H<G$, then we must have ${\rm dom}_G^{{\cal
N}_2}(H)\subseteq {\rm dom}_G^{(m,n)}(H)$, but equality need not hold
in principle. So our goal is to give a description of dominions in
$(m,n)$, and contrast it with the same dominions, but taken in~${\cal
N}_2$.

A group $G\in{\cal V}$ is said to be \textit{absolutely closed} (in
${\cal V}$) if and only if for all $K\in{\cal V}$ with $G$ a subgroup
of~$K$, we have ${\rm dom}_K^{{\cal V}}(G)=G$. The connection with
special amalgams shows that a group $G\in{\cal V}$ is absolutely
closed in~${\cal V}$ if and only if it is a special amalgamation base
in ${\cal V}$. We direct the reader to the survey paper by
Higgins~\cite{episandamalgs} for a more complete discussion of
amalgams and their connection with dominions. We will use the terms
``special amalgamation base'' and ``absolutely closed''
interchangeably. 

Intuitively dominions, in the context of amalgams, measure how hard it
is to ``upgrade'' a weak embedding into a strong embedding. That is,
when we know that a ${\cal V}$-amalgam $(A,B;D)$ is weakly embeddable,
and we wonder whether it is also \textit{strongly} embeddable
(possibly into a different group $M$), or at least to give the
smallest possible $D_0$ containing $D$ over which the amalgam is
strongly embeddable. If all dominions are trivial (i{.}e{.}, if for
all $G\in{\cal V}$ and $H<G$, ${\rm dom}_G^{{\cal V}}(H)=H$) then any
weak embedding can be upgraded into a strong embedding for
$(A,B;D)$. If dominions are ``large'', then it tends to be hard to
upgrade weak embeddings into strong ones, and the minimal $D_0$ tends
to be large compared to $D$. If they are ``small'' (meaning they
differ very little from the subgroup itself), the weak embeddings tend
to be easy to upgrade, and $D_0$ tends to be close to~$D$.

We now turn to dominions in the varieties $(m,n)$. We can obtain a
description of them as a consequence of Theorem~\ref{strongamalgs},
since the dominion of $H$ in $G$ is the smallest subgroup $D$
containing $H$ such that $(G,G;D)$ is strongly embeddable.

\begin{theorem}[cf. Theorem~3.31 in~\cite{nildoms}]
\label{domsmn}
Let $G\in(m,n)$, and let $H$ be a subgroup of
$G$. Let $D$ be the subgroup generated by $H$ and all elements of the
form $[a,b]^q$, where $q>0$, $q|n$, and $a^q,b^q\in H(G^nG')$. Then
\[D = {\rm dom}_G^{(m,n)}(H).\] 
\end{theorem}

\begin{proof}
First we prove that all such $[a,b]^q$ must lie in the
dominion. Indeed, say $K\in(m,n)$, and we have morphisms $f,g\colon
G\to K$ which agree on $H$. Then
\begin{eqnarray*}
f\left([a,b]^q\right) & = & \left[ f(a),f(b) \right]^q\\
& = & \left[ f(a^q),f(b)\right]\\
& = & \left[ f(a^q)f(a'),f(b)\right]\qquad\;\mbox{(for some $a'\in
G^nG'$ with $a^qa'\in H$)}\\
& = & \left[ f(a^qa'),f(b)\right]\\
& = & \left[ g(a^qa'),f(b)\right]\qquad\qquad\mbox{(since $a^qa'\in H$ and $f|_H=g|_H$)}\\
& = & \left[ g(a^q)g(a'),f(b)\right]\\
& = & \left[ g(a^q),f(b)\right]\qquad\qquad\quad \mbox{(since $g(a')\in
K^nK'\subseteq Z(K)$)}\\
& = & \left[ g(a),f(b)\right]^q,
\end{eqnarray*}
and by symmetry this equals $g\left([a,b]^q\right)$. So
$D\subseteq{\rm dom}_G^{(m,n)}(H)$. 

To prove equality, it suffices to show that $(G,G;D)$ is strongly
embeddable in $(m,n)$. We check conditions (a) and (b) from
Theorem~\ref{strongamalgs}.

Condition (a) is trivial, since $G^nG'\cap D\subseteq G^nG'\subseteq
Z(G)$. For condition (b), assume that $a,b\in G$ satisfy
$a^qa',b^qb'\in D$ for some $q>0$, $q|n$, $a',b'\in G^nG'$. We want to
prove that $[a^qa',b]=[a,b^qb']\in D$. The equality of the two follows
by bilinearity of the commutator bracket, which we can now apply since
we are working inside a single group $G$:
\[
[a^qa',b]=[a^q,b][a',b]=[a^q,b]=[a,b]^q=[a,b^q]=[a,b^q][a,b']=[a,b^qb'].\]
To prove that it lies in $D$, note that every element of $D$ can be
written as an element of $H$ times some commutators. Thus, if $x\in
D$, then there exist $x'\in G'$ such that $xx'\in H$. Since $a^qa'\in
D$, there is some $a''\in G'\subseteq G^nG'$ with
$a^qa'a''=a^q(a'a'')\in H$, and analogously there exists an element
$b''\in G^nG'$ with $b^q(b'b'')\in H$. By construction of~$D$, we must
have $[a,b]^q\in D$; but this is equal to $[a^qa',b]$, giving
condition (b). Thus, $(G,G;D)$ is strongly embeddable, so ${\rm
dom}_G^{(m,n)}(H)=D$.
\end{proof}

\begin{remark}
\label{ppowersfordoms}
It is straightforward to check that we may restrict $q$
to prime powers. Indeed, if $q=p_1^{a_1}\cdots p_r^{a_r}$ is a prime
factorization for $q$, and we have that all elements
\[ \left[ a^{q/p_i^{a_i}},b^{q/p_i^{a_i}}\right]^{p_i^{a_i}} =
[a,b]^{q^2/p_i^{a_i}}\]
lie in the dominion, then we have $[a,b]^k$ in the dominion, with
\[ k = {\rm gcd}
\left(\frac{q^2}{p_1^{a_1}},\ldots,\frac{q^2}{p_r^{a_r}}\right) = q.\]
\end{remark}

From Theorem~\ref{domsmn} we deduce the following corollary, which answers
Question 9.93 in~\cite{nildoms}:
\begin{corb}
Let $G\in(m,n)$, with $n$ squarefree. Then for all subgroups $H$ of~$G$, 
${\rm dom}_G^{(m,n)} (H) = H$.
\end{corb}
\begin{proof}
We may assume that $G$ is finitely generated.

Suppose first that $m>0$; by looking at the $p$-parts separately, which
we can do since dominions respect finite direct products, 
we may assume that $m=p^a$ is a prime power, and $n=p$; for if $n=1$,
then $G$ is abelian, and $H\triangleleft G$. 

Suppose that $a^q,b^q\in H(G^pG')$. Looking at Theorem~\ref{domsmn}, we
need only consider $q=1$ and $q=p$. If $q=p$, then $[a,b]^q=e$, so it
lies in $H$. If $q=1$, then, choosing $a',b'\in G^pG'$ so that
$aa',bb'\in H$, we have
\[ [a,b] = [aa',bb']\in H\]
so again $[a,b]^q\in H$. Therefore, $D=H$, and the dominion is
trivial.

Now suppose that $m=0$, and $G\in (0,n)$.
Let $x\in G\setminus
H$. We will prove that $x\notin {\rm dom}_G^{(0,n)}(H)$. Let $N$ be a
normal subgroup of $G$, with $G/N$ finite and such that $xN\notin
HN$ (the existence of such an $N$ follows from the fact that in a
finitely generated nilpotent group, every subgroup is closed in the
profinite topology).
Let $m$ be such that $G/N$ is of exponent $m$, and $n| m/{\rm
gcd}(m,2)$. Since $n$ is squarefree, the dominion of $HN/N$ in $G/N$
in $(m,n)$ is trivial by the previous case, so
\begin{eqnarray*}
 xN \notin HN/N & \subseteq & \left({\rm dom}_G^{(0,n)}(HN)\right)/N\\
& = & {\rm dom}_{G/N}^{(0,n)} \left(HN/N\right)\\
&\subseteq & {\rm dom}_{G/N}^{(m,n)}\left(HN/N\right)\\
& = & HN/N
\end{eqnarray*}
\noindent so $HN={\rm dom}_G^{(0,n)}(HN)$, and since $x$ is not in
that dominion, we conclude that
$x\notin {\rm dom}_G^{(0,n)}(H)$
either. This gives the result.
\end{proof}

\begin{corb}
\label{corgenmaierone}
Let $m\geq 0$, $n|m/{\rm gcd}(m,2)$ and $n$ squarefree. If an amalgam of
$(m,n)$ groups is weakly embeddable in $(m,n)$, then it is also
strongly embeddable in $(m,n)$.
\end{corb}
\begin{proof}
Since all dominions are trivial in $(m,n)$, every group is a special
amalgamation base. Say $K$ gives a weak embedding of $(A,B;D)$ in
$(m,n)$. Then we may strongly embed $(K,K;D)$ in $(m,n)$. By embedding $A$ into
the first copy of $K$ and $B$ into the second copy of $K$, we get a
strong embedding for $(A,B;D)$.
\end{proof}
\begin{remark} This result partially generalizes Corollary~1.3
in~\cite{nilexpp}. There, Maier proves that in the class of all
nilpotent groups of class at most $c$ and exponent a prime $p$, with
$p>c$, a weakly embeddable amalgam is always strongly embeddable. For
$c=2$, the result would cover only the varieties $(p,p)$ with $p$ an
odd prime. 
\end{remark}

Given a group $G\in (m,n)$, we must have for
all $H<G$ that
\[{\rm dom}_G^{{\cal N}_2}(H) \subseteq {\rm dom}_G^{(m,n)}(H).\]
Since there are fewer groups to map to in $(m,n)$, it would be
possible, at least in principle, that the dominion in $(m,n)$ is
strictly larger than the dominion in ${\cal N}_2$ (or in some
intermediate variety). This turns out not to be the
case:

\begin{theorem}
\label{equalitydoms}
Let $G\in(m,n)$, and let $H$ be a subgroup of~$G$. Then
${\rm dom}_G^{{\cal N}_2}(H) = {\rm dom}_G^{(m,n)}(H)$.
\end{theorem}
\begin{proof}
The inclusion of the left hand side into the right hand side is immediate,
since $(m,n)\subseteq {\cal N}_2$. To prove the reverse inclusion, say
that we have $x,y\in G$, $q>0$ with $q|n$, and $x',y'\in G^nG'$ such
that $x^qx',y^qy'\in H$. We want to prove that $[x,y]^q\in {\rm
dom}_G^{{\cal N}_2}(H)$. Note that ${\rm dom}_G^{{\cal N}_2}(H)$
contains all $[a,b]^q$ with $q>0$ and $a^q,b^q\in H[G,G]$: no
restriction on $q$ dividing $n$, but also a smaller group $HG'$
rather than $H(G^nG')$. 

Write $x'=r^nr'$ with $r\in G$, $r'\in G'$, and analogously
$y=s^ns'$. Since $q|n$, we may write $n=qk$ for some integer $k$. Then
\[\left(xr^k\right)^q = x^q r^n [r^k,x]^{q\choose 2}.\]
So $(xr^k)^qx''\in H$, where $x''=[x,r]^{k{q\choose 2}}r'\in G'$. In a
similar manner, $(ys^k)^qy''\in H$, with $y''\in G'$. Therefore,
$\left[ xr^k,ys^k \right]^q \in {\rm dom}_G^{{\cal N}_2}(H)$. 

But since $G\in(m,n)$, $n$-th powers of commutators are trivial. So:
\begin{eqnarray*}
\left[ xr^k,ys^k\right]^q & = &
\left[\left(xr^k\right)^q,ys^k\right]\\
& = & \left[x^qr^n,ys^k\right]\\
& = & \left[x^q,ys^k\right]\left[r,ys^k\right]^n\\
& = & \left[x,(ys^k)^q\right]\\
& = & \left[x,y^qs^n\right]\\
& = & \left[x,y^q\right]\left[x,s\right]^n\\
& = & [x,y]^q.
\end{eqnarray*}
Therefore, $[x,y]^q\in {\rm dom}_G^{{\cal N}_2}(H)$, as desired. This
proves the equality.
\end{proof}

In particular, ${\rm dom}_G^{{\cal N}_2}(H) = {\rm dom}_G^{{\cal
W}}(H)$, where ${\cal W}$ is the variety generated by~$G$.

\section{Return to embeddability}

We mentioned in Section~3 that, intuitively, dominions measure
how hard it is to strengthen a weak embedding into a strong one. In
view of Theorem~\ref{equalitydoms}, we might guess that if we have an
$(m,n)$-amalgam which is weakly embeddable in $(m,n)$, 
that it should be just as hard to strengthen that
embedding in ${\cal N}_2$ (i.e.~to embed the amalgam strongly in
${\cal N}_2$) as in~$(m,n)$.

In other words: we have seen that it is possible for an
$(m,n)$-amalgam to be strongly embeddable in ${\cal N}_2$ but not even
weakly embeddable in $(m,n)$. But Theorem~\ref{equalitydoms} suggests
that it may be impossible to have an $(m,n)$-amalgam which is strongly
embeddable in~${\cal N}_2$, and weakly, but \textit{not} strongly,
embeddable in~$(m,n)$. The intuition does indeed pay off:

\begin{propb} 
Let $(A,B;D)$ be an amalgam of $(m,n)$-groups.
If the amalgam is weakly embeddable in $(m,n)$ and strongly
embeddable in ${\cal N}_2$, then the amalgam is strongly embeddable in
$(m,n)$. 
\end{propb}
\begin{proof}
We check conditions (a) and (b) from Theorem~\ref{strongamalgs}. Since
$(A,B;D)$ is weakly embeddable in $(m,n)$, by condition (1) of
Theorem~\ref{weakamalgs} we have that $A^nA'\cap D\subseteq Z(B)$ and
$B^nB'\cap D\subseteq Z(A)$, giving (a). Let $q>0$, $q|n$, $a\in
A$, $r^nr'\in A^nA'$, $b\in B$, $s^ns'\in B^nB'$, with
$a^qr^nr',b^qs^ns'\in D$. We want to prove that
\[\left[ a^qr^nr',b\right] = \left[a,b^qs^ns'\right]\in D.\]

Let $k$ be such that $qk=n$. Then
\begin{eqnarray*}
\left(ar^k\right)^q[a,r^k]^{q\choose 2}r' & = & a^qr^nr'\\
\left(bs^k\right)^q[b,s^k]^{q\choose 2}r' & = & b^qs^ns'.
\end{eqnarray*}
Since $[a,r^k]^{q\choose 2}r'\in A'$, $[b,s^k]^{q\choose 2}s'\in B'$,
and the amalgam is strongly embeddable in ${\cal N}_2$, we know that
\[
\left[ \left(ar^k\right)^q[a,r^k]^{q\choose 2}r', bs^k\right] =
  \left[ ar^k, \left(bs^k\right)^q[b,s^k]^{q\choose 2}s'\right]\in
  D.\]

Since $a^qr^nr'\in D$, its $k$-th power also lies in~$D$. This is equal to:
\[\left(a^qr^nr'\right)^k  =  \left(a^qr^n\right)^k r'^k
 =  a^n r^{kn}r'^k \in A^nA',\]
\noindent so $\left(a^qr^nr'\right)^k\in A^nA'\cap D\subseteq Z(B)$, and
symmetrically we also have that $\left(b^qs^ns'\right)^k\in B^nB'\cap
D\subseteq Z(A)$.

Therefore, in $B$, we have
\begin{eqnarray*}
\left[\left(ar^k\right)^q[a,r^k]^{q\choose 2}r',bs^k\right] & = &
\left[a^qr^nr',bs^k\right]\\
& = & [a^qr^nr',b][a^qr^nr',s^k]\\
& = & [a^qr^nr',b][(a^qr^nr')^k,s]\\
& = & [a^qr^nr',b]
\end{eqnarray*}
\noindent since $(a^qr^nr')^k\in Z(B)$. Symmetrically, in $A$ we have
\[\left[ar^k,\left(bs^k\right)^q[b,s^k]^{q\choose 2}s'\right] =
[a,b^qs^ns'].\]
Therefore, $[a^qr^nr',b]=[a,b^qs^ns']\in D$, as needed.

So $(A,B;D)$ satisfies conditions (a) and (b) in
Theorem~\ref{strongamalgs}, and thus is strongly embeddable in~${\cal
N}_2$.
\end{proof}

A closer examination of the argument shows that the only property
derived from the weak embeddability of $(A,B;D)$ in $(m,n)$ which was
used was the fact that $A^nA'\cap D\subseteq Z(B)$ and $B^nB'\cap
D\subseteq Z(A)$, and we never used the fact that the amalgam also
satisfies condition (2) of Theorem~\ref{weakamalgs}. This yields the
following:

\begin{theorem}
\label{contraststrong}
Let $(A,B;D)$ be an amalgam of $(m,n)$ groups. The amalgam
is strongly embeddable in $(m,n)$ if and only if it is strongly
embeddable in ${\cal N}_2$, and in addition we have
\[A^nA'\cap D\subseteq Z(B)\qquad\mbox{and}\qquad B^nB'\cap D\subseteq
Z(A).\]
\end{theorem}

Similarly, we have:
\begin{theorem}
\label{contrastweak}
Let $(A,B;D)$ be an amalgam of $(m,n)$ groups. The amalgam
is weakly embeddable in $(m,n)$ if and only if it is weakly embeddable
in ${\cal N}_2$, and in addition we have
\[A^nA'\cap D\subseteq Z(B)\qquad\mbox{and}\qquad B^nB'\cap D\subseteq
Z(A).\]
\end{theorem}
\begin{proof}
If $(A,B;D)$ is weakly embeddable in $(m,n)$, then
trivially it is so in ${\cal N}_2$, and the condition on $A^nA'\cap D$
and $B^nB'\cap D$ follows from Theorem~\ref{weakamalgs}. For the
converse, assume that $(A,B;D)$ is weakly embeddable in ${\cal
N}_2$. Let $k>0$, $q_i>0$, $q_i|n$, $a_i\in A$, $a_i'\in A^nA'$, with
$a_i^{q_i}a_i'\in D$, $b_i\in B$, $b_i'\in B^nB'$ with
$b_i^{q_i}b_i'\in D$ as in the statement of
Theorem~\ref{weakamalgs}. We want to prove that for each $d\in D$,
\[ \prod[a_i,b_i^{q_i}b_i'] = d \Longleftrightarrow
\prod[a_i^{q_i}a_i',b_i]=d.\]

Write $a_i'=r_i^nr_i'$, $b_i'=s_i^ns_i'$. Let $n=q_it_i$. Note that,
as before,
\[(a_ir_i^{t_i})^{q_i}[a_i,r^{t_i}]^{{q_i\choose
2}}r'=a_i^{q_i}r_i^nr_i'=a_i^{q_i}a_i'.\]
Also, since $a_i^{q_i}r_i^nr_i'\in D$, so is
\[
(a_i^{q_i}r_i^nr_i')^{t_i} =
a^nr_i^{nt_i} r_i'^{t_i}\in
A^nA'\cap D\subseteq Z(B)\]
and symmetrically for $b_i^{q_i}s_i^ms_i'$. Since the amalgam is
weakly embeddable in ${\cal N}_2$, we know that for each $d\in D$,
\begin{eqnarray*}
\lefteqn{ \prod \Bigl[
(a_ir_i^{t_i}),(b_is_i^{t_i})^{q_i}[b_i,s_i^{t_i}]^{q_i\choose 2}s_i'\Bigr] = d 
\Longleftrightarrow}\\
& \Longleftrightarrow & \prod \Bigl[
(a_ir_i^{t_i})^{q_i}[a_i,r_i^{t_i}]^{q_i\choose
2}r_i',(b_is_i^{t_i})\Bigr] = d.
\end{eqnarray*}

The $i$-th factor on the left hand side of this equation yields
\begin{eqnarray*}
\Bigl[
(a_ir_i^{t_i}),(b_is_i^{t_i})^{q_i}[b_i,s_i]^{t_i{q_i\choose 2}}s_i'\Bigr]  & = &
[a_i,b_i^{q_i}b_i'][r_i^{t_i},b_i^{q_i}b_i']\\
& = & [a_i,b_i^{q_i}b_i'][r_i,(b_i^{q_i}b_i')^{t_i}]\\
& = & [a_i,b_i^{q_i}b_i']
\end{eqnarray*}
since $(b_i^{q_i}b_i')^{t_i}\in Z(A)$, and symmetrically for the right
hand side. Therefore, for each $d\in D$:
\begin{eqnarray*}
\prod [a_i,b_i^{q_i}b_i'] = d & \Longleftrightarrow & \prod
[a_ir_i^{t_i},b_i^{q_i}b_i'] = d\\
& \Longleftrightarrow & \prod [a_i^{q_i}a_i',b_is_i^{t_i}] = d\\
& \Longleftrightarrow & \prod [a_i^{q_i}a_i',b_i] = d.
\end{eqnarray*}
giving the conditions for weak embeddability in $(m,n)$, as desired.
\end{proof}

Now we see that Example~\ref{guidingex} captures the difference
between embedding in $(m,n)$, with $n>0$, and embedding in ${\cal
N}_2$: in $(m,n)$ we need, in addition, that $n$-th powers and
commutators which lie in the core be central in both factors of the
amalgam.

Let $(A,B;D)$ be an amalgam of nil-2 groups. Let $\mathcal{L}$ be the
$01$-lattice of all varieties $(m,n)$ such that that $A,B\in(m,n)$,
and let $\mathcal{F}_s$ (resp.~$\mathcal{F}_w$) be the subset of all
$(m,n)\in\mathcal{L}$ such that $(A,B;D)$ is strongly (resp.~weakly)
embeddable in $(m,n)$. Trivially, $\mathcal{F}_s$ and $\mathcal{F}_w$
are upward closed: if we have $(m,n)\subseteq (m',n')$, and
$(m,n)\in\mathcal{F}$, then $(m',n')\in\mathcal{F}$ as well. In fact,
we have a bit more:

\begin{theorem}
\label{filter}
Let $(A,B;D)$, $\mathcal{L}$, $\mathcal{F}_s$, and $\mathcal{F}_w$ be
as in the previous paragraph. Then $\mathcal{F}_s$
(resp.~$\mathcal{F}_w$) is either empty, or is a filter in
$\mathcal{L}$; explicitly, if $(m,n),(m',n')\in \mathcal{F}_s$
(resp.~$\mathcal{F}_w)$, then 
\[ (m,n)\wedge (m',n') = \left({\rm gcd}(m,m'),{\rm
gcd}(n,n')\right)\]
also lies in $\mathcal{F}_s$ (resp.~$\mathcal{F}_w$). In particular,
if $\mathcal{F}_s$ (resp.~$\mathcal{F}_w$) is non-empty, then it is a
principal filter.
\end{theorem}

\begin{proof} By Theorems~\ref{contraststrong}
and~\ref{contrastweak}, it suffices to verify that
\[ A^{\ell}A'\cap D \subseteq Z(B)\mbox{ and }B^{\ell}B'\cap
D\subseteq Z(A)\] where $\ell={\rm gcd}(n,n')$, since the amalgam is
embeddable in $(0,0)$. Write $n=\ell k$, $n'=\ell k'$, with ${\rm
gcd}(k,k')=1$. Let $a^{\ell}a'\in A^{\ell}A'\cap D$. We want to show
that it lies in $Z(B)$. Since $(a^{\ell}a')^k = a^na'^{k}\in A^nA'\cap
D$, it lies in $Z(B)$, because $(A,B;D)$ is embeddable in $(m,n)$;
likewise, we also have that $(a^{\ell}a')^{k'}\in A^{n'}A'\cap D$, and
therefore lies in $Z(B)$, since the amalgam is also embeddable in
$(m',n')$.

Therefore, $(a^{\ell}a')$ has exponent ${\rm gcd}(k,k')$ modulo
$Z(B)$ in $B$. Since ${\rm gcd}(k,k')=1$, this means that it is of
exponent $1$ modulo $Z(B)$, so it lies in $Z(B)$, as desired. A
symmetric argument shows that $B^{\ell}B'\cap D$ is contained
in~$Z(A)$, giving the result. The final statement follows because
$\mathcal{L}$ has the descending chain condition, so any filter is
principal.
\end{proof}

\subsection*{Two special cases}

We will now give analogues of two results of Maier, which deal with
the special situation in which the core $D$ is either central or
cocentral in one of the two groups in the amalgams. Recall that we say
a subgroup $D$ is cocentral in $B$ if there is a central subgroup $H$
of $B$ such that $B=\langle D,H\rangle$. Equivalently, if $B=\langle
D,Z(B)\rangle$.

\begin{propb}[cf. Satz~2, part~1 in~\cite{amalgone}]
\label{cocentralconds}
Suppose that $(A,B;D)$ is an amalgam of $(m,n)$-groups, and $D$ is
cocentral in $B$. Then the amalgam is strongly
embeddable in $(m,n)$ if and only if
\begin{itemize}
\item[(1)]~$B^n\cap D\subseteq Z(A)$; and
\item[(2)]~$\forall q>0$, $q|n$, if $a\in A$ satisfies
$a^q\in D(A^nA')$, then for all $b\in Z(B)$, $b'\in B^n$ with $b^qb'\in D$,
we have $[a,b^qb']=e$. 
\end{itemize}
\end{propb}

\begin{remark} It is worth noting that since $D$ is cocentral in $B$,
it is normal; and therefore, if the amalgam is embeddable at all, it
is strongly embeddable, as noted in Remark~\ref{normalmeansstrong}.
\end{remark}

\begin{proof}
If the amalgam is strongly embeddable, then
\[B^n\cap D \subseteq B^nB'\cap D\subseteq Z(A)\]
yielding condition (1). For condition (2), note that
$[a,b^qb'] = [a^qa',b]=e$, since $b$ is central in $B$.

Conversely, we prove that if the amalgam satisfies (1) and~(2), then
it is strongly embeddable. Note that $A^nA'\cap D\subseteq
Z(D)\subseteq Z(B)$, since $D$ is cocentral in $B$. On the other hand,
\[ B^nB'\cap D = B^nD'\cap D = \langle B^n\cap D, D'\rangle \subseteq
Z(A)\]
since $D'\subseteq Z(A)$, and $B^n\cap D\subseteq Z(A)$ by (1). 

Now let $a\in A$, $a'\in A^nA'$, $b\in B$, $b'\in B^nB'$ with $a^qa',
b^qb'\in D$ for some $q>0$, $q|n$. We want to prove that 
$[a^qa',b]=[a,b^qb']\in D$. 

Write $b=dz$, $b'=f^nx^ny'$, with $d,f\in D$, $x,z\in Z(B)$, $y'\in
B'=D'$. Then
\[ b^qb' = d^qz^qf^nx^ny' = (d^qf^n)z^qx^ny'.\]
Therefore, if $b^qb'\in D$, then so is $z^qx^n$. By condition (2),
we must have $[a,z^qx^n]=e$, and $y'\in D'\subseteq A'$, so $[a,y']=e$
as well. Therefore,
\begin{eqnarray*}
[a,b^qb'] & = & [a,d^qz^qf^nx^ny']\\
& = & [a,d^qf^n][a,z^qx^n][a,y']\\
& = & [a,d^qf^n]\\
& = & [a,d^q][a,f^n]\\
& = & [a,d^q]\\
& = & [a^q,d]\\
& = & [a^qa',d].
\end{eqnarray*}
On the other hand,
\[ [a^qa',b] = [a^qa',dz] = [a^qa',d][a^qa',z]=[a^qa',d]\]
so we have that $[a^qa',b]=[a,b^qb']$. To see that it lies in $D$,
note that it lies in $B'$, and $B'=D'$.
\end{proof}

\begin{remark} It should be noted that Maier's statement of Satz~2
part~1 in~\cite{amalgone} is in fact incorrect. There he states that for an amalgam in
${\cal N}_2$, it is enough to consider in (2) $a\in A$ with $a^q\in DA'$ and
$b\in B\setminus D$ with $b^q\in D$. The condition, as given there, is
sufficient but not necessary. Here is a
counterexample to necessity:

Let $D$ be the relatively free ${\cal N}_2$-group on two generators
$x$ and~$y$. The elements of the group can be written uniquely in the
form
\[x^ay^b[y,x]^c;\qquad a,b,c\in\mathbb{Z}.\]
Let $Z/qZ$ be the cyclic group of order $q$, generated by $z$; let
$A=B=(Z/qZ)\times D$. Then $(A,B;D)$ is strongly embeddable, for
example into $(Z/qZ)\times D\times (Z/qZ)$, and $D$ is cocentral in
$B$. Let $a=(e,x)\in A$, and let $b=(z,y)\in B\setminus D$. Then
$a^q=(e,x^q)\in D$, and $b^q=(e,y^q)\in D$. However,
$[a,b^q] = [x,y^q] = [x,y]^q \not= e$, so the condition given is not
necessary.  
\end{remark}

When $D$ is central, we get the following:

\begin{propb}[cf. Satz~2 part 2 in~\cite{amalgone}]
Suppose that $(A,B;D)$ is an amalgam of $(m,n)$ groups, and
$D\subseteq Z(B)$. Then the  amalgam is
strongly embeddable in $(m,n)$ if and only if
\begin{itemize}
\item[(1)]~$B^nB'\cap D\subseteq Z(A)$; and
\item[(2)]~For all $q>0$, $q|n$, $a\in A$, $a'\in A^nA'$ with
$a^qa'\in D$, and all $b\in B$, $b'\in B^nB'$ with $b^qb'\in D$, we
have $[a,b^qb']=e$.
\end{itemize}
\end{propb}

\begin{proof}
Again note that if the amalgam is embeddable at all, then it is
strongly embeddable, since $D$ is normal in~$B$.

If the amalgam is strongly embeddable, then (1) certainly holds, and
we have $[a,b^qb']=[a^qa',b]$. Since $D$ is central, the latter
commutator is trivial, yielding (2).

Conversely, assume the amalgam satisfies (1) and (2). We have that
\[ A^nA'\cap D\subseteq D \subseteq Z(B),\]
and together with (1), we get condition (a) from
Theorem~\ref{strongamalgs}. For condition~(b), let $x\in A$, $x'\in
A^nA'$, $y\in B$, $y'\in B^nB'$ with $x^qx',y^qy'\in D$ for some
$q>0$, $q|n$. We want to prove that $[x^qx',y]=[x,y^qy']\in D$.  But
$[x^qx',y]=e$ because $D$ is central in $B$; and condition (1) says
that $[x,y^qy']=e$. So the amalgam satisfies the conditions of
Theorem~\ref{strongamalgs}, and we are done.
\end{proof}

\section{Weak and Strong Amalgamation bases}

We pass now to the study of amalgamation bases in the varieties $(m,n)$.
We will look at the weak and strong bases first, and then deal with
the special bases separately in the next section.

\subsection*{Characterization}

The conditions for weak and strong embeddability involve two parts:
one relating to central elements (conditions (a) in
Theorem~\ref{strongamalgs}, and (1) in Theorem~\ref{weakamalgs}), and
one that relates to how certain roots of elements of the core are to
interact (conditions (b) in Theorem~\ref{strongamalgs} and (2) in
Theorem~\ref{weakamalgs}).  Similarly, the conditions for a group to
be a weak and strong base will also consist of two parts: one which is
used to ensure that condition (a) of Theorem~\ref{strongamalgs} holds,
and another which is used to ensure that condition (b) will hold for
any amalgam with the given core.

We give the following definition:
\begin{definition} Given a group $G$ and an integer $n>0$, we define
\[\Omega^{n} = \Bigl\{g\in G\,\Bigm|\, g^n=e\Bigr\}.\]
By convention, we set $\Omega^0(G)=G$.
\end{definition}

\begin{theorem}[cf. Theorem~3.3 in \cite{saracino}, Satz~5
in~\cite{amalgtwo}]
\label{weakstrongbases}
For a group $G$ in $(m,n)$, the following are
equivalent:
\begin{itemize}
\item[(i)]~$G$ is a weak $(m,n)$-amalgamation base.
\item[(ii)]~(a) $\Omega^{\beta}(Z(G))=G^nG'$, where $\beta=0$ if
$m=n=0$, and $\beta={\rm
lcm}(\frac{m}{n},n)$ otherwise; and (b) for each $q>0$ with $q|n$, and each
$g\in G$, either $g\in G^qG'$, or else no $(m,n)$-overgroup $K$ of $G$
satisfies $g\in K^qK'$.
\item[(iii)]~(a) $\Omega^{\beta}(Z(G))=G^nG'$, where $\beta=0$ if
$m=n=0$ and $\beta={\rm lcm}(\frac{m}{n},n)$ otherwise; and (b) for
each $q>0$ with $q|n$, and each $g\in G$, either $g\in G^qG'$, or
$g^{\zeta}\not=e$ where $\zeta={\rm lcm}(\frac{m}{q},n)$, or there
exists $h\in G$ such that $h^q\equiv g^c\pmod{G^nG'}$ for some $c\in
\mathbb{Z}$, and $[h,c]\not=e$.
\item[(iv)]~$G$ is a strong $(m,n)$-amalgamation base.
\end{itemize}
\end{theorem}

\begin{proof}
The equivalence of $(ii)$ and $(iii)$ follows from
Corollary~\ref{addonerootgen}. 

Clearly, $(iv)\Rightarrow(i)$. To see that $(i)\Rightarrow(ii)$,
assume that $(ii)$ does not hold. First, we note that
$G^nG'\subseteq\Omega^{\beta}(Z(G))$. For, given $g\in G$, $g'\in G'$,
then $g^ng'\in Z(G)$, and also
\[(g^ng')^{\beta} = g^{n\beta}g'^{\beta}=e\]
because $m|n\beta$, and $n|\beta$.

Therefore, if $\Omega^{\beta}(Z(G))\not=G^nG'$, then there exists
$g\in\Omega^{\beta}(Z(G))\setminus G^nG'$. Let $K_1$ be an
$(m,n)$-overgroup of $G$ such that $g\in K_1^nK_1'$. That $K_1$ exists
follows by Corollary~\ref{addonerootgen}.

Let $K_2= G\amalg^{(m,n)}\langle c\rangle$, where $\langle c\rangle$
is cyclic of order $m$, infinite cyclic if $m=0$. Then $[g,c]\not=e$
in $K_2$ by Corollary~\ref{whatcommutes}. Now consider the amalgam
$(K_1,K_2;G)$. It cannot be weakly embeddable in $(m,n)$, because
$g\in K_1^nK_1'\cap G$, but $g\notin Z(K_2)$.

Now assume instead that there exist $q>0$, with $q|n$, and $g\in G$
with $g\notin G^qG'$, but for which there exists an $(m,n)$-overgroup
$K_1$ of $G$ with $g\in K_1^qK_1'$. Let $K_2 = G\amalg^{(m,n)}\langle
c\rangle$, with $\langle c\rangle$ cyclic of order $q$. Note that
$[g,c]\not=e\in K_2$, because the cartesian is isomorphic to
\[ \frac{G}{G^nG'}\otimes \left(Z/qZ\right) \cong \frac{G}{G^qG'}.\] 
If the amalgam $(K_1,K_2;G)$ were
embeddable into $M\in(m,n)$, then we would have $g\in M^qM'$, since
$g\in K_1^qK_1'$, so we would have $[g,c]=e$ in $M$, which is a
contradiction.

Therefore, if $G$ fails (ii), it is not a weak amalgamation base.

Finally, we prove (ii)$\Rightarrow$(iv). Let $(A,B;G)$ be an amalgam
of $(m,n)$ groups, where $G$ satisfies the conditions given in
(ii). We prove $(A,B;G)$ satisfies the conditions in
Theorem~\ref{strongamalgs}.

Let $a\in A^nA'\cap G$. Since $a\in A^nA'$, $a\in Z(G)$. Also,
$a=r^nr'$ for some $r\in A$, $r'\in A'$, so
$a^{\beta} = (r^nr')^{\beta} = r^{n\beta}r'^{\beta} = e$.
Thus, we have that $a\in \Omega^{\beta}(Z(G)) = G^nG'$. Therefore, 
$a\in G^nG' \subseteq B^nB'\subseteq Z(B)$,
and symmetrically for elements of $B^nB'\cap G$. So $(A,B;G)$
satisfies condition~(a) of Theorem~\ref{strongamalgs}.

For condition (b), let $q|n$, $a\in A$, $a'\in A^nA'$, $b\in B$,
$b'\in B^nB'$ such that $a^qa',b^qb'\in G$.  By condition (ii)(b) here,
there must exist $g_1,g_2\in G$, $g_1',g_2'\in G^nG'$ such that
$a^qa'=g_1^qg_1'$, $b^qb'=g_2^qg_2'$, since it is clear that $A$ is an
overgroup of $G$ where $a^qa'\in A^qA'$, and $B$ is an overgroup of
$G$ where $b^qb'\in B^qB'$. Thus, in $B$ we have that:
\begin{eqnarray*}
[a^qa',b] & = & [g_1^qg_1',b]\\
& = & [g_1^q,b]\\
& = & [g_1,b^q]\\
& = & [g_1,b^qb']\\
& = & [g_1,g_2^qg_2']\\
& = & [g_1,g_2]^q,
\end{eqnarray*}
\noindent which clearly lies in $G$, and is equal to $[a,b^qb']$ by
symmetry. Therefore $(A,B;G)$ is strongly embeddable in $(m,n)$,
proving $(iv)$.
\end{proof}

\begin{remark} It is again straightforward to verify that we may
restrict $q$ to prime powers.
\end{remark}

\begin{remark}
\label{twobisenough}
It is worth noting that if $n>0$, then condition
(iii)(a) follows from condition (iii)(b) (and hence (ii)(a) follows
from (ii)(b)) by setting $q=n$. Let $z\in \Omega^{\beta}(Z(G))$; since $z$ cannot
satisfy the last clause of (iii)(b), it must either have
$z^{\zeta}\not=e$, where $\zeta={\rm lcm}(\frac{m}{n},n)=\beta$, or else $z\in
G^nG'$. Since $z^{\beta}=e$, we must have 
$z\in G^nG'$; the other inclusion is always true, yielding
(iii)(a). If we were to set $q=0$, then (iii)(b) is equivalent to the
statement that $G'=Z(G)$, which is not necessary for the cases
$n>0$. But if we allow $q=0$ only when $n=0$, then again we would have
that (iii)(b) implies (iii)(a). We have kept condition (iii)(a) explicitly,
however, so the statement parallels more closely Theorem~\ref{strongamalgs} and
Saracino's Theorem~3.3 in~\cite{saracino}.\end{remark}

\subsection*{Reductions and Examples}

We pause now to explore some of the
differences between strong bases in the different varieties $(m,n)$,
to give some simplifications and classifications for special cases,
and to give examples. The examples illustrate how one uses the
conditions of Theorem~\ref{weakstrongbases} in practice.

In principle, it could be that a group $G$ lying in $(m,n)$ is a strong
base in ${\cal N}_2$, but not in $(m,n)$; for perhaps some amalgam of
$(m,n)$ groups with core is embeddable in ${\cal N}_2$, but not in
$(m,n)$. However, as it turns out, this cannot occur, and we show
something stronger:

\begin{theorem}
\label{downisokayforstrong}
Let $G\in (m,n)\subseteq(m',n')$. If $G$ is a strong amalgamation base
in $(m',n')$, then it is also a strong amalgamation base in $(m,n)$.
\end{theorem}
\begin{proof}
Let $G$ be a strong $(m',n')$ base. We may assume that $n>0$, for
otherwise we would have $n=n'=m=m'=0$ and the statement is trivial.
Thus, it will suffice to check condition (ii)(b) for $(m,n)$, as per
Remark~\ref{twobisenough}. Let $q>0$, $q|n$, and let $g\in G$. Since
$G$ is a strong $(m',n')$ base, either $g$ lies in $G^qG'$, or there is no
$(m',n')$-overgroup $K$ of $G$ with $g$ in $K^qK'$. If $g\in G^qG'$, we
are done. Otherwise, there can be no $(m,n)$-overgroup $K$ of $G$ with
$g$ in $K^qK'$, since there is no $(m',n')$-overgroup with the property,
and $(m,n)\subseteq(m',n')$. This proves (ii)(b) for $(m,n)$, and
since $n>0$, this implies (ii)(a), giving the result.
\end{proof}

We also know that if $m$ is odd, then the strong amalgamation bases in
$(m,m)$ are exactly those that are strong amalgamation bases in ${\cal
N}_2$ (Theorem~4.4 in~\cite{closures}). To quickly see this, note that
no element in an $(m,m)$ group can be the ``wrong order'' to have a
$q$-th root modulo a commutator, and that $G^mG'=G'$, so the condition
for adjunction of a $q$-th root modulo the commutator while staying in
$(m,m)$ turns out to be exactly the same as for adjunction of a $q$-th
root in some ${\cal N}_2$-overgroup.

There is a very simple characterization of the abelian groups which
are strong amalgamation bases in the different $(m,n)$. We do it in
two parts: $m>0$ and $m=0$. For the first case, we may assume $m$ is a
prime power.

\begin{theorem}
\label{abelianforp}
Let $G$ be an abelian group, $G\in (p^{a+b},p^a)$, with $p$ a prime,
$a>0$, $b\geq 0$ ($b>0$ if $p=2$). 
\begin{itemize}
\item[(1)]~If $b\geq a$, then $G$ is a strong amalgamation base if
and only if \[G=\mathop{\oplus}\limits_{j\in J} (Z/p^{a+b}Z).\]
\item[(2)]~If $b<a$, then $G$ is a strong amalgamation base if and
only if $G$ is trivial.
\end{itemize}
\end{theorem}
\begin{proof}
Since $G$ is bounded, we may write $G$ as a sum of cyclic groups,
$G=\oplus (Z/p^{i_j}Z)$, with $1\leq i_j\leq a+b$.

Assume first that $b\geq a$. If $G$ is a sum of cyclic groups of order
$p^{a+b}$, then it is easy to verify that $\Omega^{p^b}(G) =
G^{p^a}=G^{p^a}G'$. Also, fix $i$ with $1\leq i\leq a$, and let $g\in
G$. There is a $(p^{a+b},p^a)$ overgroup K of $G$ with $g\in
K^{p^i}K'$ if and only if $g^{p^{a+b-i}}=e$. But this occurs if and
only if $g\in G^{p^i}$. Thus, $G$ satisfies the conditions for being a
strong $(p^{a+b},p^a)$ base.

Conversely, note that if some $i_j\leq b$, then $\Omega^{p^b}(G)$
includes the entire cyclic summand, which cannot all be $p^a$-th
powers unless the summand is trivial. And if for some $j$ we have
$b<i_j<a+b$, then the $p^{i_j-b}$-th powers are $p^a$-th powers only
if $p^a|p^{i_j-b}$, which occurs only if $a+b\leq i_j$, which is
impossible. So all cyclic summands must be of order $p^{a+b}$ for
condition (ii)(a) in Theorem~\ref{weakstrongbases} to be satisfied,
yielding the converse.

For (2), assume that $b<a$. Here we need $\Omega^{p^a}(G)=G^{p^a}$. If
some cyclic summand has $i_j<a$, then the entire summand is of
exponent $p^a$, which yields a contradiction unless the summand is
trivial. And if $i_j\geq a$, the elements annihilated by $p^a$ are the
$p^{i_j-a}$-th powers; they are $p^a$-th powers only if $2a\leq i_j$.
But $a+b<2a$, so this is impossible. Thus, an abelian group cannot
be a strong amalgamation base when $b<a$, unless it is
trivial.
\end{proof}

\begin{remark} In particular, the converse in Theorem~\ref{downisokayforstrong} is
false in general: for example, the cyclic group of order $p^5$ is a
strong base in $(p^5,p)$, but not in $(p^6,p)$. We will explore the
extent to which the converse holds later in the section.
\end{remark}

\begin{remark} The case $a=0$ is easy: the class $(p^b,1)$ is a
class of abelian groups, so all groups are strong amalgamation bases.
\end{remark}

\begin{theorem}
\label{absformzero}
Let $G$ be an abelian group, $G\in (0,n)$. Then $G$ is a strong
amalgamation base in $(0,n)$ if and only if it is $n$-divisible
(i.e.~$G^n=G$).
\end{theorem}
\begin{proof}
Note that for $n=0$, the claim is that an abelian group is a strong
$(0,0)$-amalgamation base if and only if $G=G^0=\{e\}$, which also
follows from the fact that we also need $Z(G)=G'$. So we may assume $n>0$.

Let $q|n$, $g\in G$; the third clause of condition (iii)(b) from
Theorem~\ref{weakstrongbases} is always false in $G$, since we cannot
have $[h,g]\not=e$ for any $h\in G$; and the condition on the order of
$g$ is trivially false, since $\zeta=0$. So, if $G$ is a strong
amalgamation base, then we must have $g\in G^qG'=G^q$. Setting $q=n$
yields that $G$ is $n$-divisible. Conversely, if $G$ is $n$-divisible,
then it is $q$-divisible for all $q|n$, giving condition (iii)(b), and
from there we get condition (iii)(a) as stated in
Remark~\ref{twobisenough}.
\end{proof}

\begin{corb}
\label{corabsformzero}
Let $G$ be an abelian group, $G\in (0,n)$ with $n>0$. Then $G$ is a
strong $(0,n)$ amalgamation base if and only if it is $p$-divisible
for every prime $p$ with $p|n$.
\end{corb}

We have seen that condition (ii)(b) in Theorem~\ref{weakstrongbases}
is necessary and sufficient for a group to be a strong $(m,n)$
amalgamation base, provided that $n>0$. There is one situation in which (ii)(a) is both
necessary and sufficient. To describe it we need a couple of
lemmas. 

Let $\pi$ be a set of primes. We say that a group is $\pi$-divisible
if every element has a $p$-th root for each prime $p\in\pi$; we say it
is $\pi'$-divisible if every element has a $q$-th root for each prime
$q\notin\pi$.

\begin{lemmab}
\label{ppartsall}
Let $G=A\oplus B\in (m,n)$, and let $\pi$ be a set of primes. If $A$
is $\pi$-divisible and $B$ is $\pi'$-divisible, then $G$ is a strong
amalgamation base in $(m,n)$ if and only if $A$ and $B$ both are.
\end{lemmab}
\begin{proof}
It is easy to verify that, in general, if $A\oplus B$ is an
amalgamation base (weak, strong, or special), then so are both $A$ and
$B$. Conversely, suppose that both $A$ and $B$ are strong amalgamation
bases in $(m,n)$.

Since $\Omega^{\beta}(Z(G))=\Omega^{\beta}(Z(A))\oplus\Omega^{\beta}(Z(B))$,
and $G^nG'=A^nA'\oplus B^nB'$, it follows that $G$ satisfies condition
(ii)(a) from Theorem~\ref{weakstrongbases}. For condition (ii)(b),
assume $q=p^{\alpha}$ is a prime power, $q|n$, and let $(a,b)\in G$. 

If $p\notin\pi$, then $b\in B^qB'$, since $B$ is
$\pi'$-divisible. Suppose that there exist an
$(m,n)$-overgroup $K$ of~$G$ with $(a,b)\in K^qK'$. Since $b\in
B^qB'\subseteq K^qK'$, we would have $(a,b)(e,b^{-1})\in K^qK'$, which
shows that $(a,e)\in K^qK'$. Since $A$ is a strong base, and $K$ is
also an overgroup of~$A$, we must have $a\in A^qA'$; but if both $a\in
A^qA'$ and $b\in B^qB'$, then we have $(a,b)\in G^qG'$, yielding
(ii)(b), and we are done.

If, on the other hand, $p\in\pi$, then $a\in A^qA'$ since $A$ is
$\pi$-divisible, and the symmetric argument holds.
\end{proof}

\begin{lemmab}[cf. Theorem 3.5 in~\cite{saracino}]
\label{ppartsanywhere}
Let $A,B\in (m,n)$, and assume that $A$ and $B$ are of relatively
prime exponents. Then $A\oplus B$ is a strong $(m,n)$-amalgamation
base if and only if both $A$ and $B$ are.
\end{lemmab}
\begin{proof} Use Lemma~\ref{ppartsall}, with $\pi$ the set of all
primes not dividing the exponent of $A$.
\end{proof}

\begin{corb}
\label{cor:pparts}
Let $G\in (m,n)$ be a torsion group. Then $G$ is a strong
$(m,n)$-amalgamation base if and only if the $p$-parts of $G$ are.
\end{corb}

We may now give the promised class of groups for which condition
(ii)(a) is necessary and sufficient:

\begin{lemmab}[cf. Theorem~3.6 in~\cite{saracino}]
Let $G\in (m,n)$ be of exponent $k$, where $k$ is either squarefree or
twice a squarefree number. Then $G$ is a strong amalgamation base for
$(m,n)$ if and only if $\Omega^{\beta}(Z(G))=G^nG'$, where $\beta=0$
if $n=m=0$, and $\beta={\rm lcm}(\frac{m}{n},n)$ otherwise.
\end{lemmab}

\begin{proof}
By Corollary~\ref{cor:pparts}, we may assume that $k=p$ a prime,
or $k=4$. Necessity follows from Theorem~\ref{weakstrongbases}, so we
only need to prove sufficiency.

Assume first that $k=p$. Let $q$ be a prime power dividing $n$. Since
$G$ is $p'$-divisible, for any prime $p'\not=p$, we may assume that
$q=p^a$, with $a\geq 1$. Let $g\in G$, and we want to prove that either $g\in
G^{p^a}G'$, or else no $(m,n)$-overgroup $K$ of $G$ has a $p^a$-th
root for $g$. Note that $p^a|\beta$, so
$\Omega^{\beta}(Z(G))=Z(G)=G^nG'$. So, if $g\in Z(G)$, then $g\in
G^nG'\subseteq G^{p^a}G'$, and we are done. Otherwise, there exists
$h\in G$ with $[h,g]\not=e$, but $h^{p^a}\equiv g^0\pmod{G^nG'}$, thus
showing that $g$ satisfies the last clause of condition (iii)(b). Thus
$G$ is a strong amalgamation base for $(m,n)$.

Next, suppose that $G$ is of exponent 4. Then $G^2\subseteq Z(G)$. If
$n$ is odd, then $G$ is $q$-divisible for any $q|n$, so (iii)(b)
always holds and we are done. If $4|n$, then $G^n$ is trivial and
$4|\beta$, so the condition on the center stated in the theorem is
equivalent to $Z(G)=G'$. In particular, $G^2\subseteq G'$. Let $q|n$
be a power of $2$ (any other prime power will yield that (iii)(b)
holds), and let $g\in G$. If $g$ is central, then $g\in G'$ and the
first clause of (iii)(b) holds. Otherwise, there is some $h\in G$ such
that $[h,g]\not=e$; and since $G^2\subseteq G'$, $h^q\equiv
g\pmod{G'}$, so the last clause of (iii)(b) holds instead.

Finally, assume that ${\rm ord}_2(n)=1$. Note that we may then assume
that $q=2$. If $4|\beta$, we may proceed
as above. Otherwise, we must have that ${\rm ord}_2(m)=2$, so that
${\rm ord}_2(\beta)=1$. Then $G^2G'$ is exactly the central elements
of exponent $2$. Again, if $g$ is not central, we get an $h$ which
does not commute with $g$, to satisfy the last clause of (iii)(b); if
$g$ is central of exponent $2$, then it lies in $G^2G'$ and the first
clause of (iii)(b) is satisfied. And if $g$ is of order 4, then
$g^{\zeta}\not=e$, so the second clause is satisfied, and we are done.
\end{proof}

Note in addition that if $k|n$, then $k|\beta$ and $G^n=\{e\}$, so the
condition becomes $Z(G)=G'$, which does not depend on $m$ and
$n$. 

\begin{corb}
Let $G\in(m,n)$ be a group of exponent $k$, $k|n$. If $k$ is
squarefree or twice a squarefree number, then the following are
equivalent:
\begin{itemize}
\item[(1)]~$G'=Z(G)$.
\item[(2)]~$G$ is a strong $(m,n)$-amalgamation base.
\item[(3)]~$G$ is a strong $(m',n')$-amalgamation base for all $n'$
with $k|n'$.
\end{itemize}
\end{corb}

\begin{example}
In general, condition (ii)(a) is not enough for groups of exponent
a multiple of $p^2$, with $p$ an odd prime and $p^2|n$, or multiple of
$8$ with $4|n$. For let $p$ be a prime, and consider the ${\cal N}_2$-group
\[G = \langle x,y,z\,|\,
x^{p^2}=y^{p^2}=z^p=[x,y]^{p^2}=[x,z]^p=[y,z]^p=e\rangle.\]
If $G\in(m,n)$, then $p^2|n$, and $z\notin G^{p^2}G'$, but there is an
$(m,n)$-overgroup $K$ of $G$ with $z\in K^{p^2}K'$, so $G$ cannot be
an $(m,n)$ base. 
\end{example}

We would like to know, in general, whether there are $(m,n)$ bases that
are not $(m',n')$-bases for $(m,n)$ properly contained in
$(m',n')$. It is not hard to see that there are some trivial cases
where this will be impossible, for example when $n=n'=1$. To
fully explore the question, we begin with a series of examples:

\begin{example}
\label{bsmall}
 Let $G\in (p^{a+b},p^a)$, with $p$ a prime, $a,b>0$,
be the relatively free group of rank $2$, that is:
\[ G = \Bigl\langle x,y\,\Bigm|\,
x^{p^{a+b}}=y^{p^{a+b}}=[x,y]^{p^a}=[x,y,x]=[x,y,y]=e\Bigr\rangle.\]
Then $G$ is a strong $(p^{a+b},p^a)$-amalgamation base. If $b\leq a+1$,
then it is not a strong $(p^{a+b},p^{a+1})$-amalgamation base.

First, we prove it is a strong $(p^{a+b},p^a)$-amalgamation base. Let
$i>0$ be an integer with $1\leq i\leq a$, and let $g\in G$. Assume
that there exists some overgroup $K\in (p^{a+b},p^a)$ with $g\in
K^{p^i}K'$. We want to prove that $g\in G^{p^i}G'$. Since $g$ is a
$p^i$-th power times a commutator in $K$, $g^{p^{a-i}}$ must be
central in $G$. If $g=x^{\alpha}y^{\beta}[x,y]^{\gamma}$, then we must
have
\begin{eqnarray*}
e & = & [g^{p^{a-i}},x]\\
  & = & [x^{\alpha p^{a-i}}y^{\beta p^{a-i}},x]\\
  & = & [y,x]^{\beta p^{a-i}}.
\end{eqnarray*}
So we must have $p^a|\beta p^{a-i}$, or equivalently $p^{i}|\beta$. A
symmetric argument holds for $\alpha$, so $g\in G^{p^i}G'$, as
desired. Thus $G$ satisfies condition (ii)(b) of
Theorem~\ref{weakstrongbases}, which we know is sufficient for $G$ to
be a strong amalgamation base in $(p^{a+b},p^a)$.

Next, assume that $b\leq a+1$.  We claim $G$ is not a strong
$(p^{a+b},p^{a+1})$-base in this case. Indeed, consider $x^{p^a}$,
which does not have a $p^{a+1}$-th root modulo $G^{p^{a+1}}G'$. The
element $x^{p^a}$ is central in $G$, and in addition, if $\zeta={\rm
lcm}(p^{a+b-(a+1)},p^{a+1})=p^{a+1}$, then
\[(x^{p^a})^{\zeta}=(x^{p^a})^{p^{a+1}} = x^{p^{2a+1}} =e\]
since $b\leq a+1$, so $a+b\leq 2a+1$. Thus, $x^{p^a}$ has a
$p^{a+1}$-th root in some $(p^{a+b},p^{a+1})$-overgroup of $G$, but it
does not lie in $G^{p^{a+1}}G'$, so $G$ does not satisfy condition
(ii)(b) in Theorem~\ref{weakstrongbases} relative to
$(p^{a+b},p^{a+1})$.

\end{example}

\begin{example}
\label{bbig}
Let $p$ be a prime, $a$ and $b$ integers with $b>a+1$,
$a>0$. Let
\[ G = \Bigl\langle x,y \,\Bigm|\,
x^{p^{a+b}}=y^{p^{a+b-1}}=[x,y]^{p^a}=[x,y,x]=[x,y,y]=e
\Bigr\rangle.\]
Then $G$ is a strong $(p^{a+b},p^a)$-amalgamation base, but not a
strong amalgamation base for $(p^{a+b},p^{a+1})$.

First, we prove that $G$ satisfies condition (ii)(b) from
Theorem~\ref{weakstrongbases}, relative to $(p^{a+b},p^a)$. Let $i$ be
an integer, $1\leq i\leq a$, and let $g\in G$. If there is a
$(p^{a+b},p^a)$-overgroup $K$ of $G$ with $g\in K^{p^i}K'$, then
$g^{p^{a-i}}$ must be central in $K$, hence central in $G$. It is now
easy to verify that this means that $g$ is of the form $x^{\alpha
p^i}y^{\beta p^i}[y,x]^{\gamma}$ for some integers
$\alpha,\beta,\gamma$, so $g\in G^{p^i}G'$, proving condition
(ii)(b). Thus, as per Remark~\ref{twobisenough}, $G$ is a strong
amalgamation base in $(p^{a+b},p^a)$. 

To prove that it is not a strong amalgamation base in
$(p^{a+b},p^{a+1})$, consider $i=a+1$ and the element
$g=y^{p^a}$. Trivially, $g\notin G^{p^{a+1}}G'$; however, there is a
$(p^{a+b},p^{a+1})$-overgroup $K$ of $G$ with $g\in K^{p^{a+1}}K'$:
since $g$ is central, it suffices to test if it has the right order,
i.e.~if $g^{p^{b-1}}=e$. But this is true, since $y$ is of order $p^{a+b-1}$.
This proves that $G$ cannot be a strong
amalgamation base in $(p^{a+b},p^{a+1})$.
\end{example}

\begin{example}
\label{firstentryall}
Let $p$ be a prime, and $a,b>0$. Let $G\in
(p^{a+b},p^a)$ be given by:
\[ G = \left\langle x,y,z \,\Biggm|\,
\begin{array}{rcl}
x^{p^{a+b}} & = & y^{p^{a+b}} = [y,x]^{p^a} = [z,x] = [z,y] = e;\\
z^{p^b} & = & [y,x];\quad [x,y,x]=[x,y,y]=e.
\end{array}
\right\rangle.\]
Then $G$ is a strong amalgamation base in $(p^{a+b},p^a)$, but not in
$(p^{a+b+1},p^a)$.

First we show it is a strong amalgamation base in $(p^{a+b},p^a)$. Let
$K$ be a $(p^{a+b},p^a)$ overgroup of $G$, and let $g\in
K^{p^i}K'$. We want to show that $g\in G^{p^i}G'$. We may write
$g=x^{\alpha}y^{\beta}z^{\gamma}[y,x]^{\delta}$, where $0\leq
\alpha,\beta<p^{a+b}$, $0\leq \gamma<p^b$, $0\leq \delta<p^a$. We must
have that $g^{p^{a-i}}$ is central in $G$. Since
\[g^{p^{a-i}}  =  x^{\alpha p^{a-i}} y^{\beta p^{a-i}} z^{\gamma
p^{a-i}} [y,x]^{\delta p^{a-i} + \alpha\beta{p^{a-i}\choose 2}}\]
it follows, by taking commutators with $x$ and $y$, that $p^i|\alpha$,
$p^i|\beta$. So we may rewrite $g$ as:
\[ g = x^{\alpha p^i}y^{\beta p^i}z^{\gamma}[y,x]^{\delta}\]
for new integers $\alpha$, $\beta$, $\gamma$, and~$\delta$.

If $i\leq b$, then we must also have that
$g^{p^{a+b-i}}=e$. Therefore
\begin{eqnarray*}
e & = & g^{p^{a+b-i}}\\
& = & x^{\alpha p^{a+b}} y^{\beta p^{a+b}} z^{\gamma p^{a+b-i}}
[y,x]^{\delta p^{a+b-i} + \alpha\beta p^{2i}{p^{a+b-i}\choose 2}}\\
& = & z^{\gamma p^{a+b-i}}.
\end{eqnarray*}
Therefore, $p^i|\gamma$, and so $g\in G^{p^i}G'$, as desired. On the
other hand, if $i>b$, then we must have that $g^{p^a}=e$, and we have:
\begin{eqnarray*}
e & = & g^{p^a}\\
& = & x^{\alpha p^{a+i}} y^{\beta p^{a+i}} z^{\gamma p^a}
[y,x]^{\delta p^a + \alpha\beta p^{2i}{p^{a}\choose 2}}\\
& = & z^{\gamma p^a}.
\end{eqnarray*}
Therefore, $p^b | \gamma$; but then we may write
$z^{\gamma}=z^{\gamma'p^b}=[y,x]^{\gamma'}$, so we have
\[g = x^{\alpha p^i} y^{\beta p^i} [y,x]^{\delta'},\]
and $g\in G^{p^i}G'$ again. This proves that condition (iii)(b) of
Theorem~\ref{weakstrongbases} holds, so $G$ is a strong amalgamation
base in $(p^{a+b},p^a)$.

To prove that $G$ is not a strong base in $(p^{a+b+1},p^a)$, consider
$z$ and $i=1$. Since $z$ is central, the last clause of condition
(iii)(b) does not hold; and since clearly $z\notin G^{p}G'$, the first
clause does not hold either. As for the second clause, since $b>0$ we
must test the $p^{a+b}$-th power of $z$, but that is trivial, so it
does not satisfy the second clause either. Therefore, $z$ and $p$ do
not satisfy clause (iii)(b) relative to $(p^{a+b+1},p^a)$, so $G$
cannot be a strong base in that variety.
\end{example}

\begin{example}
\label{forzero}
Let $n>0$, and let $G\in(0,n)$ be given by
\[ G = \Bigl\langle x,y \Bigm|
[x,y]^n=[x,y,x]=[x,y,y]=e\Bigr\rangle.\]
Then $G$ is a strong $(0,n)$ amalgamation base, but not a strong
$(0,n')$ amalgamation base for any $(0,n')$ properly containing
$(0,n)$.

First, we prove it is a strong $(0,n)$ base. Let $q|n$, $q>0$, and
$g\in G$. If there is some $(0,n)$-overgroup $K$ of $G$ with $g\in
K^qK'$, then $g^{n/q}$ must be central in $G$; since $Z(G)=G^nG'$,
that means that $g^{n/q}\in G^nG'$, and it is easy to verify that this
means that $g\in G^qG'$. Thus, $G$ satisfies (ii)(b) from
Theorem~\ref{weakstrongbases}, and therefore it is a strong
$(0,n)$-amalgamation base.

But if $n'\not=n$, then $Z(G)\not=G^{n'}G'$; but this is condition
(ii)(a) for a group to be a $(0,n')$ amalgamation base; thus, $G$
cannot be a $(0,n')$ base if $(0,n')$ properly contains $(0,n)$.
\end{example}

Using the examples, we have the following result, which states exactly
to what extent the converse of Theorem~\ref{downisokayforstrong}
holds. Recall that by convention we set ${\rm ord}_p(0)=\infty$.

\begin{theorem}
\label{howfewremain}
Let $(m,n)\subseteq(m',n')$ be subvarieties of ${\cal N}_2$. Every
strong $(m,n)$-amalgamation base is also a strong amalgamation base in
$(m',n')$ if and only if for each prime $p$ one of the following holds:
\begin{itemize}
\item[(a)]~${\rm ord}_p(n)={\rm ord}_p(n')=0$; or
\item[(b)]~${\rm ord}_p(n)={\rm ord}_p(m)$; or
\item[(c)]~${\rm ord}_p(n)={\rm ord}_p(n')$ and ${\rm ord}_p(m)={\rm
ord}_p(m')$.
\end{itemize}
\end{theorem}

\begin{proof}
First, assume that for some $p$, all three conditions fail. If ${\rm
ord}_p(n)=0$, then $p|n'$ and $p|m$. Then we take a sum of two cyclic
groups of order $p$, which is an $(m,n)$-base (by
Corollary~\ref{corabsformzero} if $m=0$, or if $m>0$, because the
$p$-parts consist only of abelian groups), but not an $(m',n')$ base,
regardless of $m'$. If, on the other hand, ${\rm ord}_p(n)>0$, then
${\rm ord}_p(n)<{\rm ord}_p(m)$, and either ${\rm ord}_p(n)<{\rm
ord}_p(n')$, or else ${\rm ord}_p(m)<{\rm ord}_p(m')$. 

If both $m$ and $m'$ are zero, then we have $(0,n)\subseteq (0,n')$,
and the $p$ order of $n$ is strictly less than the $p$ order of $n'$,
so $n\not=n'$. We may use Example~\ref{forzero} to find a group $G$
which is a $(0,n)$-base but not a $(0,n')$-base, so not every
$(0,n)$-base is a $(0,n')$-base.

If ${\rm ord}_p(m)={\rm ord}_p(m')$, and both $m$ and $m'$ nonzero,
then we may use either Example~\ref{bsmall} or Example~\ref{bbig}
(depending on whether ${\rm ord}_p(m)\leq 2{\rm ord}_p(n)+1$ or not,
respectively), to construct a group $G$ which is a strong amalgamation
base for $(p^{{\rm ord}_p(m)},p^{{\rm ord}_p(n)})$ (and thus also an
$(m,n)$ base), but it is not a strong amalgamation base in $(p^{{\rm
ord}_p(m)},p^{{\rm ord}_p(n)+1})$, hence also not an $(m',n')$
amalgamation base.

If ${\rm ord}_p(m)<{\rm ord}_p(m')$, then we may use
Example~\ref{firstentryall} to construct a $p$-group which is a strong
$(p^{{\rm ord}_p(m)},p^{{\rm ord}_p(n)})$ base (and hence an $(m,n)$
base, since its $p$-parts are either a strong base or trivial), but
not a strong base in $(p^{{\rm ord}_p(m)+1},p^{{\rm ord}_p(n)})$, and
hence cannot be an $(m',n')$-base either.

This proves that the conditions (a), (b), and (c) are necessary. To
prove sufficiency, we consider the finite exponent cases first: if
both $m$ and $m'$ are nonzero, we may assume they are powers of the
same prime $p$; given any prime $q\not=p$, (c) holds for $q$. If (c)
also holds for $p$ that means that $m=m'$ and $n=n'$, so trivially the
strong bases in both are the same. If (a) holds for $p$, then $(m,n)$
is the class of all abelian groups of exponent $m$, and $(m',n')$ of
all abelian groups of exponent $m'$, and since $m|m'$ and all groups
are strong bases in both, we again have that every $(m,n)$-base is
also an $(m',n')$-base. Finally, if (b) holds for $p$, then we have
$(m,n)=(p^a,p^a)$, and in particular, $p$ must be an odd prime. But if
$p$ is an odd prime and $a>0$, then every $(p^a,p^a)$-strong
amalgamation base is also an ${\cal N}_2$-strong amalgamation base by
Theorem~4.4 in~\cite{closures}, so it must also be an $(m',n')$ base
by Theorem~\ref{downisokayforstrong}.

Next assume that $m=m'=0$. If for any $p$ we have that (b) holds, then
$n=0$, so $(m,n)=(0,0)=(m',n')$, and there is nothing to
prove. Otherwise, either (a) or (c) hold for each prime, so $n=n'$,
and again there is nothing to prove.

Finally, assume that $m>0$ and $m'=0$. Then (c) can never happen, so
for each prime either (a) or (b) hold; and we may assume that both $m$
and $n$ are power of a prime $p$. For primes other than $p$, (b)
always holds. For the prime $p$, if (a) holds, then we have
$(m,n)=(p^a,1)$ and $(m',n')=(0,n')$, with ${\rm gcd}(p,n')=1$. The groups in
$(p^a,1)$ are all abelian, and they are $n'$-divisible, so they are
$(0,n')$-strong bases by Theorem~\ref{absformzero}. And if (b) holds,
then we again have $(m,n)=(p^a,p^a)$, and we apply
Theorem~4.4 in~\cite{closures} to get that every $(p^a,p^a)$-strong
amalgamation base must also by an $(m',n')$-amalgamation base. This
proves the theorem.
\end{proof}

Let $G\in{\cal N}_2$, and consider the $01$-lattice $\mathcal{L}$ of
all subvarieties of ${\cal N}_2$ containing~$G$. Let $\mathcal{I}$ be
the subset (possibly empty) of all subvarieties $(m,n)$ such that $G$
is a strong amalgamation base in $(m,n)$. We have seen that if
$(m',n')\in \mathcal{I}$, and $(m,n)\in\mathcal{L}$ satisfies
$(m,n)\subseteq (m',n')$, then $(m,n)$ is also in $\mathcal{I}$. This
would lead to the dual question to Theorem~\ref{filter}, of whether
$\mathcal{I}$ is an ideal (that is, whether $\mathcal{I}$ is also
closed under joins). However this is false in general:

\begin{example}
\label{notanideal} Let $G$ be a sum of two cyclic groups of order $p^3$,
with $p$ a prime. Then $G$ is a strong $(p^5,1)$-base, because
everything in $(p^5,1)$ is a strong amalgamation base. It is also a
strong $(p^3,p)$-base, by Theorem~\ref{abelianforp}. The join of
$(p^5,1)$ and $(p^3,p)$ is $(p^5,p)$, but again by
Theorem~\ref{abelianforp} we have that $G$ is not a strong
amalgamation base there. So $\mathcal{I}$ is not closed under joins in
general.
\end{example}

All examples I have where $\mathcal{I}$ fails to be closed under joins
involve abelian groups, or cases where some prime dividing $n'$ does
not divide $n$. We might ask if there are conditions on $n$ and $n'$
that would guarantee that, if each lies in $\mathcal{I}$, then so does
their join. For example, I propose the following two question:

\begin{openquestion} Let $G$ be a $p$-group, and assume that for
$a,b,c,d>0$, $G$ is both a strong $(p^{a+b+c+d},p^a)$-amalgamation
base, and also a strong $(p^{a+b+c},p^{a+b})$-amalgamation base. Is it
also a strong $(p^{a+b+c+d},p^{a+b})$-amalgamation base? If not, under
what conditions on $a,b,c,d$ will it hold?
\end{openquestion}

\begin{openquestion} Let $G$ be a group, and let $n,n'\in \mathbb{Z}$
be positive integers such that for every prime $p$, $p|n$ if and only
if $p|n'$. If $G$ is both a strong $(0,n)$ and a strong $(0,n')$ base,
is it also a strong $(0,{\rm lcm}(n,n'))$ amalgamation base?
\end{openquestion}

In the first question we restrict $a>0$ to avoid the kind of
counterexample given in Example~\ref{notanideal}. We set $b>0$ and
$d>0$ because if either one is equal to zero, then one variety
contains the other, and the result is trivially true; and we have set
$c>0$ because otherwise, the second variety is of the form
$(p^{a+b},p^{a+b})$, and we know that every strong base there is a
strong base in $(0,0)$, and the result would again hold. It is not
hard to verify that if $a=b=c=d$, then the join will also lie in
${\mathcal I}$.

\section{Special amalgamation bases}

\subsection*{Characterization}

We pass now to the special amalgamation bases. Looking at
Theorem~\ref{strongamalgs}, we see that there are essentially two
things that can go wrong in an amalgam $(A,B;D)$; things that
``should'' be central are not; and certain roots of elements of $D$ in
the two different groups do not interact properly.

The conditions for a group to be a strong amalgamation base take care
of the two situations; condition (ii)(a) in
Theorem~\ref{weakstrongbases} ensures that anything which ``should''
be central will be central in any overgroup; and condition (ii)(b) in
Theorem~\ref{weakstrongbases} ensures that roots of elements of $D$
will interact properly. The proof that (i) implies (ii) in
Theorem~\ref{weakstrongbases} is done by constructing two incompatible
overgroups of $G$ when $G$ fails to satisfy~(ii).

When we deal with special amalgamation bases, it is impossible to
construct incompatible overgroups, since both groups are isomorphic
over the core. Instead, we need to ensure that certain commutators not
only agree, but in fact lie in $D$. Which commutators? Looking at
Theorem~\ref{domsmn}, we see that the situation we must be careful of
is when two elements of $G$ can be expressed as $q$-th powers times a
commutator in some overgroup. This leads to the following characterization:

\begin{theorem}[cf. Theorem~2.9 in \cite{absclosed}, Theorem~2.2 in \cite{closures}]
\label{absclosedgen}
Let $G$ lie in~$(m,n)$.  Then $G$ is a special amalgamation base in $(m,n)$
if and only if for all $x,y\in G$ and all integers $q>0$, $q|n$, at
least one of the following occurs:
\begin{itemize}
\item[(a)]~$x^{\zeta}\not= e$ or $y^{\zeta}\not= e$, where
$\zeta={\rm lcm}(\frac{m}{q},n)$; or
\item[(b)]~There exists $g_1,g_2\in G$, and $\alpha,\beta,\gamma,\delta\in\mathbb{Z}$,
with $\beta\equiv \gamma\pmod{\frac{n}{q}}$, such that
\[\begin{array}{rcl}
g_1^{q} & \equiv & x^{\alpha}y^{\beta}\\
g_2^{q} & \equiv & x^{\gamma}y^{\delta}
\end{array}
\pmod{G^{n}G'}\quad\mbox{and}\quad [g_1,x][g_2,y]\not=e;\mbox{ or}\]
\item[(c)]~There exists $g_1,g_2\in G$,
$\alpha,\beta,\gamma,\delta\in\mathbb{Z}$, with
$(\gamma-\beta)\equiv 1\pmod{\frac{n}{q}}$, such that
\[\begin{array}{rcl}
g_1^{q} & \equiv & x^{\alpha}y^{\beta}\\
g_2^{q} & \equiv & x^{\gamma}y^{\delta}
\end{array}
\pmod{G^{n}G'}.\]
\end{itemize}
\end{theorem}

\begin{remark}
For a group $G\in(m,n)$, $g_1,g_2\in G$, saying that they either
satisfy (a) or (b) above is equivalent to saying that there is no
$(m,n)$-overgroup $K$ of~$G$ such that $g_1,g_2\in K^{q}K'$, by
Theorem~\ref{addtworootsgen}. So the statement is equivalent to
saying that either the congruence in (c) has a solution, or no
$(m,n)$-overgroup $K$ of $G$ has $g_1,g_2\in K^qK'$. Phrased that way,
it parallels more clearly condition (ii)(b) in
Theorem~\ref{weakstrongbases}. 
\end{remark}

\begin{remark} 
\label{whenqequalsn}
If $q=n$ then condition (c) trivially holds, by setting $g_1=g_2=e$,
$\alpha=\beta=\gamma=\delta=0$, since the congruence relation between
$\beta$ and $\gamma$ becomes a congruence modulo 1. We keep
$q|n$, rather than specify proper divisors, so it parallels the
previous results.
\end{remark}

\begin{remark}
Note as well that the conditions depend only on the congruence classes
of $x$ and $y$ modulo $G^nG'$, and are symmetric on $x$ and~$y$.
\end{remark}

\begin{proof}
First we prove sufficiency: let $K$ be an $(m,n)$-overgroup of $G$,
$K\in(m,n)$. We want to prove that ${\rm
dom}_K^{(m,n)}(G)=G$.

Let $r,s\in K$, $r',s'\in K^{n}K'$, $q|n$, such that
$r^{q}r',s^{q}s'\in G$. We want to prove that $[r,s]^{q}\in G$. Let
$x=r^{q}r'$, $y=s^{q}s'$. Clearly, $K$ is an $(m,n)$-overgroup of
$G$ where $x$ and $y$ lie in $K^{q}K'$, so  conditions (a) and (b)
cannot hold for $x$, $y$, and $q$. Therefore, (c) holds. But then, by
Lemma~\ref{basiceqs}, we must have, for some $k\in \mathbb{Z}$,
\[ [r,s]^{q(\gamma-\beta)} = [r,s]^{q(1+k(n/q))} = [r,s]^q = [g_1,x][g_2,y] \in G\]
which proves that $[r,s]^{q}\in G$. Thus, $G$ equals its own
dominion in $K$, giving sufficiency.

For necessity, assume that $x,y\in G$, and a given $q>0$, $q|n$, do
not satisfy conditions (a), (b), nor (c). We construct an overgroup
$K\in(m,n)$, where $G\not={\rm dom}_K^{(m,n)}(G)$.

For ease, we do this one prime at a time in the case $m>0$. If $n=1$,
then condition (c) always holds, as mentioned in
Remark~\ref{whenqequalsn}, since the only possible $q$ is $q=n=1$. If
$n>1$, we replace $m$ with $p^{a+b}$, $n$ with $p^a$ ($a>0$, $b\geq 0$
if $p$ is odd, $b>0$ if $p=2$), and replace $q$ with $p^i$, where $1\leq
i\leq a$.

We proceed as we did in Theorem~\ref{addtworootsp}. Let 
\[ K_0 = G\amalg^{{\cal N}_2}\left(\langle r\rangle
\amalg^{(0,p^a)}\langle s\rangle\right),\]
with $\langle r\rangle$, $\langle s \rangle$ infinite cyclic, and let
$N$ be the least normal subgroup containing $xr^{-p^i}$ and
$ys^{-p^i}$. Since conditions (a) and (b) fail, we will conclude, as
we did in Theorem~\ref{addtworootsp}, that $N\cap G=\{e\}$. We also
show that $g[r,s]^{-p^i}\notin N$ for all $g\in G$. This will show that
in $K_0/N$, $[r,s]^{p^i}\notin G$. 

Proceeding as we did in Theorem~\ref{addtworootsp}, we set a general
element of~$N$ equal to $\mathbf{abc}$, and we have $\mathbf{a}=g$,
$\mathbf{b}=[r,s]^{-p^i}$, and $\mathbf{c}=e$. We also get, using equations
(\ref{valueg}), (\ref{valuebeta}), and~(\ref{valuegamma}) in the proof
of Theorem~\ref{addtworootsp}, that they
are equal to:
\begin{eqnarray*}
\mathbf{a} & = & [g_1,x][g_2,y]\\
\mathbf{b} & = & [r,s]^{p^i(c_{21}-c_{12})}\\
\mathbf{c} & = & \left[g_1^{-p^i}x^{c_{11}}y^{c_{12}},r\right]
\left[g_2^{-p^i}x^{c_{21}}y^{c_{22}},s\right]
\end{eqnarray*}
Since we also know that $\mathbf{b} = [r,s]^{-p^i}$, we must have that
\[p^i(c_{21}-c_{12})\equiv -p^i\pmod{p^a},\]
or equivalently, that
$c_{12}-c_{21}\equiv 1 \pmod{p^{a-i}}$. However, this we know to be
impossible, since (c) is not satisfied either. Thus, $G$ is properly
contained in its dominion in $K_1=K_0/N$. If $i\leq b$ we are done, by
looking at the dominion of $G$ in $\langle G,r,s\rangle$, which
contains $[r,s]^{p^i}$ by construction.

If $i>b$, we adjoin central $p^a$-th roots $t$ and $s$ to $x^{p^{a-i}}$ and
$y^{p^{a-i}}$, and then convert certain powers of those roots into
commutators, as we did in the proof of
Theorem~\ref{addtworootsp}. Neither process will collapse
$[r,s]^{p^i}$ into $G$, since they can all be performed via central
amalgams. So again, in the group~$K$ that we obtain after all the
adjunctions are done, we have that $[rt^{-1},sv^{-1}]^{p^i}$
lies in the dominion, but its value is the same is $[r,s]^{p^i}$. So
there is an element in the dominion of $G$ which is not in $G$, and so
$G$ cannot be a special amalgamation base in $(p^{a+b},p^a)$.

For the $m=0$ case, we proceed again as above, and we will always be
in the situation analogous to $i\leq b$.
\end{proof}

\begin{remark} Again, it is not hard to prove that we may restrict $q$
to prime powers, using for example an argument similar to that in
Remark~\ref{ppowersfordoms}.
\end{remark}

\subsection*{Reductions and Examples}

Since dominions are equal in ${\cal N}_2$ and in any subvariety
$(m,n)$ with $m>0$, we have the following easy analogue of
Theorem~\ref{downisokayforstrong}:

\begin{theorem}
\label{downisokayforspecial}
Let $G\in (m,n)\subseteq (m',n')$. If $G$ is a special
amalgamation base in $(m',n')$, then it is also a special amalgamation
base in~$(m,n)$.
\end{theorem}
\begin{proof}
Let $K\in(m,n)$ be any overgroup of $G$. If $G$ is a special
amalgamation base in $(m',n')$, then
\[ {\rm dom}_K^{(m,n)}(G) = {\rm
dom}_K^{(m',n')}(G) = G.\]
So $G$ is absolutely closed in~$(m,n)$ as well.
\end{proof}

\begin{corb}
If $G\in(m,n)$ is cyclic, then it is absolutely closed in~$(m,n)$.  If
$G\in(m,n)$ is of squarefree exponent and $Z(G)/G'$ is cyclic or
trivial, then it is absolutely closed in $(m,n)$. If $G\in(m,n)$ is
abelian, and $G/G^p$ is cyclic or trivial for every prime $p$, then $G$
is absolutely closed in~$(m,n)$.
\end{corb}
\begin{proof} The result follows because in each case described, $G$ will be
absolutely closed in ${\cal N}_2$:
see Theorem~3.7 in~\cite{absclosed} (cyclic groups); Theorem~3.13
in~\cite{absclosed} (squarefree exponent); and Theorem~3.17
in~\cite{absclosed} (abelian group, general case).
\end{proof}

Another useful result uses the fact that dominions respect
quotients. We have:

\begin{propb} 
Let $G\in (m,n)$; if $G^{\rm ab}$ is absolutely closed in $(m,n)$,
then so is $G$. If $G/G^n$ is absolutely closed in $(m,n)$, then so is
$G$. If $G/G^nG'$ is absolutely closed in $(m,n)$, then so is $G$.
\end{propb}
\begin{proof}
Let $K$ be any $(m,n)$ overgroup of $G$. Then $G'$, $G^n$, and $G^nG'$
are all normal in $K$, since they are central subgroups. Let $N$ be
any of them. Then:
\[G/N\subseteq \left({\rm dom}_K^{(m,n)}(G)\right)/N = {\rm dom}_{K/N}^{(m,n)}(G/N)=G/N\]
since $G/N$ is absolutely closed. But that means that ${\rm
dom}_K^{(m,n)}(G)=G$, as desired.
\end{proof}

\begin{remark}
It is important to note, however, that the converse to the above
result does not hold. E.g., see Example~3.19 in~\cite{absclosed}.
\end{remark}

We want simplifying lemmas for absolutely closed groups
similar to the ones we proved for weak and strong amalgamation bases,
as well as simpler characterizations for some subclasses.

\begin{lemmab}
\label{ifithasaroot}
Let $G\in{\cal N}_2$, and let $q>0$. if $x\in G^qG'$, then for all $y$
there exists $g_1,g_2\in G$, $a,b,d\in\mathbb{Z}$, such that
\begin{eqnarray*}
g_1^q & \equiv & x^a y^b \pmod{G'}\\
g_2^q & \equiv & x^{b+1} y^d \pmod{G'}.
\end{eqnarray*}
\noindent Analogously, if $y\in G^qG'$, we may find a solution to the
system of congruences.
\end{lemmab}
\begin{proof}
Suppose that $x=r^nr'$, and $y\in G$. Let $a=b=d=0$, $g_1=e$, and
$g_2=r$. On the other hand, if $y=s^ns'$, and $x\in G$, then set
$g_1=s^{-1}$, $g_2=e$, $a=d=0$, and $b=-1$. 
\end{proof}

\begin{corb}[cf. Corollary~2.13 in~\cite{absclosed}]
If $G\in(m,n)$ is such that for every $q>0$, $q|n$, and every $x\in
G$, either $x\in G^qG'$ or else no $(m,n)$-overgroup $K$ of $G$ has
$x\in K^qK'$, then $G$ is absolutely closed in $(m,n)$.
\end{corb}
\begin{proof}
Given $G$, $x,y\in G$, and $q|n$, $q>0$, if either $x$ or $y$ lies in
$G^qG'$, then condition (c) in Theorem~\ref{absclosedgen} is satisfied
by Lemma~\ref{ifithasaroot}. If neither lies in $G^qG'$, then no
$(m,n)$-overgroup of $G$ has either one as a $q$-th power modulo a
commutator, hence it cannot have both of them, proving that either (a)
or (b) holds in Theorem~\ref{absclosedgen}.
\end{proof}

In particular, it follows that any weak amalgamation base is also a
special amalgamation base, which would also imply that every weak base
is also a strong base, by Theorem~\ref{weakstrongdifference}. This
gives an alternative way to prove that condition (i) and condition
(iv) in Theorem~\ref{weakstrongbases} are equivalent.

\begin{corb} Let $G\in(0,n)$, with $n>0$. If $G=G^nG'$, then $G$ is absolutely
closed in $(0,n)$. If $G$ is divisible, then it is absolutely closed in~${\cal N}_2$.
\end{corb}

The following observation is useful when looking at particular
examples:

\begin{lemmab}[Perturbation Argument]
Let $G\in(m,n)$, let $H$ be a subgroup of~$G$, $q>$, $q|n$, and
$x,y\in G$ such that $x^q,y^q\in H(G^nG')$. If $h_1,h_2\in H$, then
\[[x,y]^q\in H \Longleftrightarrow [xh_1,yh_2]^q\in H.\]
\end{lemmab}
\begin{proof}
Note that if $x^q,y^q\in H(G^nG')$, then so do $(xh_1)^q$ and
$(yh_2)^q$. Thus, both commutator brackets lie in the dominion of
$H$. Expanding the bracket bilinearly, we have
\begin{eqnarray*}
[xh_1,yh_2]^q & = & [x,y]^q[h_1,y]^q[x,h_2]^q[h_1,h_2]^q\\
& = & [x,y]^q[h_1,y^q][x^q,h_2][h_1,h_2]^q.
\end{eqnarray*}
Since $x^q,y^q\in H(G^nG')$, and $h_i\in H$, the last three terms on
the right hand side lie in $H$. Thus, the left hand side lies in~$H$
if and only if the $[x,y]^q$ lies in~$H$, as claimed.
\end{proof}

\begin{lemmab}[cf. Theorem~3.2 in~\cite{absclosed}]
Let $G=A\oplus B\in(m,n)$ and let $\pi$ be a set of primes. If $A$ is
$\pi$-divisible and  $B$ is $\pi'$-divisible, then $G$ is absolutely
closed in $(m,n)$ if and only if both $A$ and $B$ are.
\end{lemmab}
\begin{proof}
In general, if $A\oplus B$ is absolutely closed, then so are $A$ and
$B$. Conversely, assume both $A$ and $B$ are absolutely closed in
$(m,n)$, and let $K$ be an $(m,n)$-overgroup of $A\oplus B$. Let
$x,y\in K$, $x',y'\in K^nK'$, and $q>0$, $q|n$, such that
$x^qx',y^qy'\in A\oplus B$. We want to prove that $[x,y]^q\in A\oplus
B$. We may assume that $q=p^r$ is a prime power. Write
$x^nx'=(a_1,b_1)$, $y^ny'=(a_2,b_2)$. 

If $p\in\pi$, $a_1^{-1}$ and $a_2^{-1}$ both have $q$-th roots in $A$, so there exists
$s,t\in A$ such that $s^q=a_1^{-1}$, $t^q=a_2^{-1}$. Then
\begin{eqnarray*}
(sx)^q \equiv s^qx^q & \equiv & (e,b_1)\pmod{K^nK'}\\
(ty)^q \equiv t^qy^q & \equiv & (e,b_2)\pmod{K^nK'},
\end{eqnarray*}
so $[sx,ty]^q\in{\rm dom}_K^{(m,n)}(B)=B$. By the perturbation
argument, $[sx,ty]^q$ lies in $A\oplus B$ if and only $[x,y]^q\in A\oplus B$,
so the latter lies in $A\oplus B$, as desired. A symmetric argument
holds if $p\notin \pi$.
\end{proof}

With this result in hand, we obtain the analogues of
Lemma~\ref{ppartsanywhere} and Corollary~\ref{cor:pparts}:

\begin{lemmab}
\label{ppartsspecial}
Let $A,B\in(m,n)$ be groups of relatively prime exponents. Then
$A\oplus B$ is absolutely closed in $(m,n)$ if and only if both $A$
and $B$ are.
\end{lemmab}
\begin{corb}
\label{cor:ppartsspecial}
Let $G\in(m,n)$ be a torsion group. Then $G$ is absolutely closed in
$(m,n)$ if and only if its $p$-parts are.
\end{corb}

We use these results to help us characterize the abelian groups and
groups of squarefree exponent which are special amalgamation bases. We
begin with the abelian groups, and as we did in the case of strong
amalgamation bases, we deal with $p$-groups first, and handle the
$(0,n)$ case separately.

\begin{theorem}
\label{abelianspecialchar}
Let $G\in(p^{a+b},p^a)$ be an abelian group, where $p$ is a prime and
$a>1$, $b\geq 0$ ($b>0$ if $p=2$). Then $G$ is a special amalgamation base
in $(p^{a+b},p^a)$ if and only if $G=C\oplus B$, where $C$ is a cyclic
group, and
\begin{itemize}
\item[(i)]~if $b\geq a-1$, then $B=\oplus(Z/p^{a+b}Z)$.
\item[(ii)]~if $b<a-1$, then $B$ is trivial.
\end{itemize}
\end{theorem}

\begin{proof}
First, assume that $b\geq a-1$. Let $G=C\oplus\left(\oplus_{j\in
J}(Z/p^{a+b}Z)\right)$. We want to prove that $G$ is absolutely closed in
$(p^{a+b},p^a)$. Denote the generators by $z$ for $C$, and $r_j$ for
the $j$-th cyclic factor of order $p^{a+b}$. Let $K$ be a
$(p^{a+b},p^a)$ overgroup of $G$, and assume that $x,y\in K$,
$x',y'\in K^{p^a}K'$, and some $i$ with $1\leq i\leq a$ has
$x^{p^i}x',y^{p^i}y'\in G$. We want to prove that $[x,y]^{p^i}\in
G$. If $i=a$, there is nothing to do, since the commutator is
trivial, so assume $i<a$. In particular, $b\geq i$.

Write $x^{p^i}x' = (z^{\alpha},\oplus r_j^{s_j})$, $y^{p^i}y' =
(z^{\beta},\oplus r_{j}^{t_j})$. Their $p^{a+b-i}$-th powers are
trivial, so we must have that $p^{a+b}|s_jp^{a+b-i},t_jp^{a+b-i}$,
Therefore, each $s_j$ and each $t_j$ are multiples of $p^i$. Let $v$
and $w$ be $p^i$-th roots of those elements (we may assume that $v$
and $w$ lie
in $\oplus (Z/p^{a+b}Z)$); perturbing $x$ and $y$ by $v^{-1}$
and $w^{-1}$, we may assume that $x^{p^i}x'$ and $y^{p^i}y'$ both lie
in~$C$. But that means that $[x,y]^{p^i}$ is in the dominion of $C$,
which equals $C$ (since cyclic groups are absolutely closed), and this
proves that the original $[x,y]^{p^i}\in G$. Thus $G$ is absolutely
closed.

For the necessity in the case $b\geq a-1$, assume that $G$ has two
cyclic summands of order less than $p^{a+b}$, and let $x$ and $y$ be
the generators. Since we may adjoin central $p$-th roots to $x$ and
$y$, it follows that $x$, $y$, and $p$ do not satisfy conditions (a)
or (b) of Theorem~\ref{absclosedgen}. To see they also do not satisfy
(c), note that if $x^{\alpha}y^{\beta}$ is congruent, modulo
$G^{p^a}G'$, to a $p$-th power, then $p$ must divide both $\alpha$ and
$\beta$. So if we could solve the congruence in condition~(c), we would have
that $p$ divides both $\beta$ and $\gamma$, which makes it impossible
for $\gamma-\beta$ to be congruent to $1$ modulo $p^{a-1}$. This
proves~(i).

For (ii), note that since $G$ is cyclic, it is
absolutely closed in ${\cal N}_2$, and therefore also in
$(p^{a+b},p^a)$. Conversely, assume that $G$ is not cyclic; decompose
$G$ into a sum of cyclic summands, each of order $p^{i_j}$, $1\leq
i_j\leq a+b$, and let $x$ and $y$ be generators of distinct cyclic
summands. Consider $x^{p^{a-2}}$, $y^{p^{a-2}}$, and
$p^i=p^{a-1}$. Since $b<a-1$, for $x^{p^{a-2}}$ and $y^{p^{a-2}}$ to satisfy condition (a)
of Theorem~\ref{absclosedgen} we would need for either of their
$p^a$-th powers to be nontrivial; but since $a+b\leq 2a-2$, this is
not the case. Since the group is abelian, they also do not satisfy
condition (b). As for condition (c), once again, if some $x^{\alpha
p^{a-2}}y^{\beta p^{a-2}}$ is congruent to a $p^{a-1}$-th power modulo
$G^{p^a}$, then we must have that $p$ divides both $\alpha$ and
$\beta$; we conclude, as we did in case~(i), that condition (c) cannot be
satisfied, for we would have that $p|\beta$ and $p|\gamma$, and yet
$\gamma-\beta\equiv 1\pmod{p}$, which is impossible. This proves
necessity.
\end{proof}

\begin{remark} We do not consider the cases $a=0$ and $a=1$, because
in both cases every dominions is trivial, so every group is a special
amalgamation base.
\end{remark}

\begin{theorem}[cf. Theorem~3.17 in~\cite{absclosed}]
\label{abspecialinzero}
Let $G$ be an abelian group lying in $(0,n)$. Then $G$ is absolutely closed in
$(0,n)$ if and only if for every prime $p$ with $p|n$, either
${\rm ord}_p(n)=1$ or $G/G^p$ is cyclic or trivial.
\end{theorem}
\begin{proof}
First, assume that there is some prime $p$ such that $p^2|n$ and $G/G^p$
is not cyclic. Since $G/G^p$ is the sum of cyclic summands, choose $x$, $y$
in $G$ which project to generators of distinct cyclic summands.
Then note that if for some $g\in G$ we have $g^p\equiv
x^ay^b\pmod{G^n}$ then we must have that both $a$ and $b$ are
multiples of $p$. Since $G$ is abelian and $m=0$, $x$ and $y$ cannot
satisfy conditions (a) or (b) in Theorem~\ref{absclosedgen} relative
to $p$. But if they satisfied condition (c), then we would have
integers $\alpha,\beta,\gamma,\delta$, all multiples of $p$, and
$\gamma-\beta\equiv 1 \pmod{n/p}$. But $n/p$ is a multiple of $p$, as
are $\gamma$ and $\beta$, so this is impossible. Therefore, $G$ is not
absolutely closed in $(0,n)$.

For sufficiency, let $q=p^a$ be a prime power, $q|n$, and $p^2|n$. Let
$K$ be a $(0,n)$ overgroup of $G$, $x,y\in K$, $x',y'\in K^nK'$, with
$x^qx',y^qy'\in G$. We want to prove that $[x,y]^q\in G$. Write
$\mathbf{x}=x^qx'$, $\mathbf{y}=y^qy'$.

If $G/G^p$ is cyclic, then so is $G/G^{p^a}$ (since the latter has the
same number of cyclic summands as $G/G^p$). Let $z\in G$ project to a
generator of $G/p^aG$. Then we may perturb $x$ and $y$ by elements of
$g_1$, $g_2$ of $G$ so that their $q$-th powers lie, modulo $K'$, in
$\langle z\rangle$: write $x^qx'$ as a $z^b w^{p^a}$ for some
$w\in G$, and set $g_1=w^{-1}$, and similarly for $y^qy'$.  Then
$[xg_1,yg_2]^q$ lies in the dominion of $\langle z\rangle$, which is
absolutely closed; by the Perturbation Argument, 
$[x,y]^q$ also lies in $G$, and we are done.

Finally, assume that $q=p$, $p|n$, but $p^2$ does not divide $n$. Then
${\rm gcd}(p,n/p)=1$, so there exists $k\in\mathbb{Z}$ such that $pk\equiv
1\pmod{\frac{n}{q}}$. In particular, we have
\begin{eqnarray*}
\mathbf{y}^p & \equiv & \mathbf{x}^0\mathbf{y}^p\pmod{K^nK'}\\
(\mathbf{x}^{k+1})^p & \equiv & \mathbf{x}^{p(k+1)}\mathbf{y}^0
\pmod{K^nK'}.
\end{eqnarray*}
By Theorem~\ref{basiceqs}, $[x,y]^{p(p(k+1)-p)}\in G$. However,
$p(p(k+1)-p)=p(pk)$. Since $pk\equiv 1\pmod{\frac{n}{p}}$, $p^2k\equiv
p\pmod{n}$, so $[x,y]^p=[x,y]^{p^2k}$, and since the latter lies in
$G$, we are done.
\end{proof}

For the groups of squarefree exponent, we have:

\begin{theorem}[cf. Theorem~3.13 in~\cite{absclosed}]
Let $p$ be a prime, $a>1$, $b\geq 0$ ($b>0$ if $p=2$). A group $G\in
(p^{a+b},p^a)$ of exponent $p$ is absolutely closed in $(p^{a+b},p^a)$
if and only if $Z(G)/G'$ is cyclic or trivial.
\end{theorem}

\begin{proof}
Since $G/G'$ is a vector space over $Z/pZ$, let $\{z_i\}_{i\in I}$ be
elements of $Z(G)$ which project onto a basis for $Z(G)/G'$. Then pick
elements $\{b_j\}_{j\in J}$ whose projections extend
$\{\overline{z_i}\}$ into a basis for $G/G'$. Since $G$ is of
exponent~$p$, $\langle z_i\rangle$ is a direct summand for $G$. If
$|I|>1$, then we have a sum of two cyclic groups, both of order less
than $p^{a+b}$, so the sum is not absolutely closed, and since $G$ has
a direct summand which is not absolutely closed, it cannot itself be
absolutely closed. This proves necessity.

Now assume that $|I|\leq 1$. Let $K$ be an overgroup of $G$, $K\in
(p^{a+b},p^a)$. If some element $g$ of $G$ lies in $K^{p^i}K'$,
$i\geq 1$, then it must be central in $G$, for given any $h\in G$,
\[ [g,h] = [r^{p^i},h] = [r,h^{p^i}] = [r,e] = e.\]

So assume that $x,y\in K$, $x',y'\in K^{p^a}K'$, and
$x^{p^i}x',y^{p^i}y'\in G$. Then these two elements are central in
$G$, so we can perturb $x'$ and $y'$ by elements of $G'$ so that both
$x^{p^i}x''$ and $y^{p^i}y''$ are
powers of $z$, whose projection generates $Z(G)/G'$. Thus, $[x,y]^{p^i}$
lies in the dominion of $\langle z\rangle$, which is 
absolutely closed. Therefore, $[x,y]^{p^i}\in G$, proving sufficiency.
\end{proof}

\begin{theorem}
Let $G\in (0,n)$, $G$ of exponent $k$. If $k$ is squarefree, then $G$
is absolutely closed in $(0,n)$ if and only if for each prime $p$
either ${\rm ord}_p(n)\leq 1$ or $Z(G)/(G^p\cap Z(G))G'$ is cyclic or
trivial.
\end{theorem}

\begin{proof}
We may assume that $G$ is a $q$-group, with $q$ a prime. If $q$ does
not divide~$n$,
then $G$ is abelian, and we know that $G$ is absolutely closed in $(0,n)$
if and only if for every prime either ${\rm ord}_p(n)\leq 1$
or $G/G^p$ is cyclic or trivial, which is the condition we have
above. So we may assume that $q|n$, in which case $G^n=\{e\}$.

Assume first that $G$ satisfies the condition, and let $p^a$ be a
prime power dividing $n$. If $q\not=p$, then $G$ is $p$-divisible, so
any pair of elements will satisfy condition~(c) of
Theorem~\ref{absclosedgen}. If $q=p$ and ${\rm ord}_p(n)=1$, then
again, given any pair of elements $x,y\in G$, find $k$ with $kp\equiv
1\pmod{\frac{n}{p}}$. Then we have
\begin{eqnarray*}
y^p & \equiv & x^0y^p\pmod{G^nG'}\\
(x^{k+1})^p & \equiv & x^{pk+p}y^0\pmod{G^nG'}
\end{eqnarray*}
and $(pk+p)-p\equiv pk\equiv 1 \pmod{\frac{n}{p}}$, which shows that
$x$ and $y$ again satisfy condition~(c) of Theorem~\ref{absclosedgen}.

If, on the other hand, ${\rm ord}_p(n)>1$, then we know that $Z(G)$
modulo $(G^p\cap Z(G))G'$ equals $Z(G)/G'$, and is cyclic.  So
$Z(G)/G'$ is generated by the image of some element $z\in Z(G)$. Since
an element of $G$ has a $p$-th root modulo a commutator in a
$(0,n)$-overgroup of $G$ if and only if it is central, we either have
that $x$ and $y$ are both central, or else they satisfy condition~(b)
of Theorem~\ref{absclosedgen}. If they are both central, we may assume
(by perturbing them by commutators) that they are both powers of
$z$. But then they satisfy condition~(c) as elements of $\langle
z\rangle$, since a cyclic group is absolutely closed, and therefore
they also satisfy condition~(c) as elements of $G$. So the condition
given is sufficient.

Now assume that there is a prime $p$ with $p^2|n$ and $Z(G)/(G^p\cap Z(G))G'$
neither cyclic nor trivial. First note that if $q\not=p$, then
$G^p=G$, so the quotient would be trivial. Thus, we must have
$q=p$. And since $G$ is a $p$-group, the condition given is that
$Z(G)/G'$ is not cyclic nor trivial. Then $G$ is not absolutely closed
in $(p^2,p^2)$, since it has a direct summand which is the sum of two
cyclic groups, so it cannot be absolutely closed in $(0,n)$ either.
\end{proof}

However, it should be noted that we cannot extend the characterization
of absolutely closed groups of squarefree exponent to those of
exponent twice a squarefree number. Here is a counterexample:

\begin{example} A group $G\in (4,2)$, with $Z(G)/G'$ cyclic, which is
not absolutely closed in $(8,4)$ (and so, not absolutely closed in any
variety containing $(8,4)$).

Let $G$ be the $(4,2)$ group given by:
\[ G = \Bigl\langle x,y,z \,\Bigm|\,
x^4=y^2=z^2=[x,y]^2=[x,z]^2=[y,z]=e\Bigr\rangle.\]
Then $G$ is of exponent four, $G/G'=Z/4Z\oplus Z/2Z\oplus Z/2Z$. The
center of $G$ is generated modulo $G'$ by $x^2$, so $Z(G)/G'$ is
cyclic. However, we can embed $G$ into the group
\[ F = \Bigl\langle a,b,c \,\Bigm|\,
a^4=b^4=c^4=[a,c]^4=[a,c]^4=[b,c]^4=e\Bigr\rangle\]
which lies in $(8,4)$, by identifying $x$ with $a$, $y$ with $b^2$,
and $z$ with $c^2$. Here $[b,c]^2$ lies in the dominion of $G$ in
$F$ in $(8,4)$, but does not lie in $G$. So $G$ is not absolutely
closed in $(8,4)$.

Similar examples can be constructed to show that $Z(G)/G'$ is not
sufficient for finitely generated torsion groups of exponent $p^n$,
with $n>1$, and $p$ an arbitrary prime.
\end{example}

Both as a way to explore the differences between the absolutely closed
groups in the different subvarieties, and as illustrations of how one
uses the results we have obtained to establish that a group is or is
not absolutely closed, we present a series of examples and
results. Recall that if $G$ is a strong $(p^a,p^a)$ base, then it is
also a strong base in ${\cal N}_2$. For special bases, we have a
partial analogue:

\begin{propb}
\label{bequalaspecial}
Let $p$ be a prime, $a\geq 0$. If a group $G\in(p^a,p^a)$
is absolutely closed, then it is also absolutely closed in $(0,p^a)$,
and therefore, in $(p^{a+b},p^a)$ for every $b\geq 0$.
\end{propb}

\begin{proof} Note that conditions (b) and (c) from
Theorem~\ref{absclosedgen} do not depend on the value of $m$; so as
long as we do not change $n$ they remain the same. If $G$ lies in
$(p^a,p^a)$, then condition (a) will always be false for any $x,y\in
G$, and any $p^i|p^a$, because the exponent that appears there is always a
multiple of $p^a$, and $x^{p^a}=y^{p^a}=e$. So if $G$ is absolutely
closed in $(p^a,p^a)$, for any $p^i|p^a$ and any pair of elements
$x,y\in G$, either (b) or (c) hold; and in that case, they also hold
in the larger variety $(0,p^a)$.
\end{proof}

The next example, however, shows that we cannot expect to extend
Proposition~\ref{bequalaspecial} to a full analogue of the result for
strong bases.

\begin{example}(cf. Theorem~3.4 in~\cite{closures})
\label{advanceinboth}
Let $p$ be a prime, $a>0$, $b\geq 0$ ($b>0$ if $p=2$); then
\[ G = \Bigl\langle x,y \,\Bigm|\,
x^{p^{a+b}}=y^{p^{a+b}}=[x,y]^{p^{a-1}}=[x,y,x]=[x,y,y]=e\Bigr\rangle\]
is a special amalgamation base in $(p^{a+b},p^a)$, but not in
$(p^{a+b+1},p^{a+1})$.

For $a=1$ the result follows because the group is abelian; it is
absolutely closed in $(p^{b+1},p)$ because everything is in that
variety, but not in $(p^{b+2},p^2)$ because $G$ is the sum of two cyclic
groups of order $p^{b+1}$, and it cannot be absolutely closed
regardless of $b$. So assume that $a>1$.

Consider the $(p^{a+b+1},p^{a+1})$-group 
\[K_1=\langle r,s\,\Bigm|\,
r^{p^{a+b+1}}=s^{p^{a+b+1}}=[r,s]^{p^{a+1}}=[r,s,r]=[r,s,s]=e\Bigr\rangle.\]
It contains $G$ as a subgroup, by identifying $x$ with $r^p$, and $y$
with $s^p$. However, $[r,s]^p$ lies in the dominion of~$G$ in $K_1$, and
not in~$G$, so $G$
cannot be absolutely closed in $(p^{a+b+1},p^{a+1})$.

Suppose, on the other hand, that $K$ is an overgroup of $G$, with $K$
a group in $(p^{a+b},p^a)$. Let $r,s\in K$, $r',s'\in K^{p^a}K'$, with
$r^{p^i}r',s^{p^i}s'\in G$ for some $i$, $1\leq i<a$. We want to prove
that $[r,s]^{p^i}\in G$.

Write 
\begin{eqnarray*}
r^{p^i} & \equiv & x^{\alpha}y^{\beta}=\mathbf{x}\pmod{K^{p^a}K'}\\
s^{p^i} & \equiv & x^{\gamma}y^{\beta}=\mathbf{y}\pmod{K^{p^a}K'}.
\end{eqnarray*}
If $i\leq b$, then $\mathbf{x}^{p^{a+b-i}}=e$, and
\[\mathbf{x}^{p^{a+b-i}} = x^{\alpha p^{a+b-i}}y^{\beta
p^{a+b-i}}[y,x]^{\alpha\beta {p^{a+b-i}\choose 2}}\]
so we must have $p^i|\alpha$, $p^i|\beta$, and similarly with
$\mathbf{y}$. Thus $\mathbf{x},\mathbf{y}\in G^{p^i}G'$. But then,
$\mathbf{x}\equiv g^{p^i}$, $\mathbf{y}\equiv h^{p^i} \pmod{G'}$ for some
$g,h\in G$, and we have:
\[ [r,s]^{p^i} =
[r^{p^i},s]=[\mathbf{x},s]=[g^{p^i},s]=[g,s^{p^i}]=[g,\mathbf{y}]\in
G\]
as desired.

If $i>b$, then we still have that
$\mathbf{x}^{p^{a-i}},\mathbf{y}^{p^{a-i}}$ must be central. It is not
hard to verify that the center of $G$ is equal to $G^{p^{a-1}}G'$;
thus, we can conclude that $\alpha$, $\beta$, $\gamma$, and $\delta$
are multiples of $p^{i-1}$. So we may rewrite them as:
\begin{eqnarray*}
\mathbf{x} & = & x^{\zeta p^{i-1}}y^{\eta p^{i-1}}\\
\mathbf{y} & = & x^{\lambda p^{i-1}}y^{\mu p^{i-1}},
\end{eqnarray*}
and perturbing $\mathbf{x}$ and $\mathbf{y}$ with $x^{p^i}$ and
$y^{p^i}$ if necessary, we may assume that
\[ 0 \leq \zeta,\eta,\lambda,\mu < p.\]

Consider the system of congruences
\begin{eqnarray*}
e^{p^i} & \equiv & \mathbf{x}^0\mathbf{y}^0\pmod{K^{p^a}K'}\\
(x^{\zeta p^{a-i-1}}y^{\eta p^{a-i-1}})^{p^{i}} & \equiv & \mathbf{x}^{p^{a-i}}\mathbf{y}^0
\pmod{K^{p^a}K'}.
\end{eqnarray*}
By Lemma~\ref{basiceqs}, we have that
\begin{eqnarray*}
 [r,s]^{p^ip^{a-i}} &=& \left[ x^{\zeta p^{a-i-1}}y^{\eta p^{a-i-1}},
 x^{\lambda p^{i-1}}y^{\mu p^{i-1}}\right]\\
& = & [x,y]^{p^{a-2}(\zeta\mu-\eta\lambda)}.
\end{eqnarray*}
Since $[r,s]^{p^ip^{a-i}}=[r,s]^{p^a}=e$, we conclude that
$p|\zeta\mu-\eta\lambda$.  Considering $(\zeta,\eta)$, $(\lambda,\mu)$
as vectors over $Z/pZ$, this means that they are proportional, so we
may perturb $\mathbf{x}$ and $\mathbf{y}$ so they generate a cyclic
subgroup; let $\mathbf{z}$ be a generator. Then $[r,s]^{p^i}$ lies in
the dominion of $\langle\mathbf{z}\rangle$; but cyclic groups are
absolutely closed, so $[r,s]^{p^i}\in\langle\mathbf{z}\rangle\subseteq
G$, as desired.

Thus $G$ is absolutely closed in $(p^{a+b},p^a)$.
\end{example}

\begin{example}
\label{bbigspecial}
Let $p$ be a prime, $a>1$, $b\geq a-1$. Let 
\[ G = Z/p^{a+b}Z \oplus Z/p^{a+b}Z\]
and denote by $x$ and $y$ the generators of the cyclic summands.  Then
$G$ is absolutely closed in $(p^{a+b},p^a)$, but not in
$(p^{a+b+1},p^a)$.

Since $a>1$ and $b\geq a-1$, $G$ is absolutely closed, as seen before.
However, it cannot be absolutely closed in $(p^{a+b+1},p^a)$ by
Theorem~\ref{abelianspecialchar}.
\end{example}

\begin{example}
\label{advancesecondsp}
Let $p$ be a prime, $a>1$, $b\geq 1$,  and let $G$ be the
$(p^{a+b},p^a)$ group presented by:
\[ G = \Bigl\langle x,y,z\,\Bigm|\,
x^{p^{a+b-1}}=y^{p^{a+b}}=z^{p^{a+b-1}}=[x,y]^{p^a}=[y,z]^{p^a}=[x,z]=e\Bigr\rangle
\]
together with all identities of $(p^{a+b},p^a)$. Then $G$ is a special amalgamation base
in $(p^{a+b},p^a)$, but not in $(p^{a+b},p^{a+1})$.

Indeed, assume first that some element $g\in G$ is a $p^i$-th power
times a commutator in some overgroup $K\in(p^{a+b},p^a)$. Then
$g^{p^{a-i}}$ is central in $G$. Writing 
\[g=x^{\alpha}y^{\beta}z^{\gamma}[y,x]^{\delta}[y,z]^{\eta}\]
it is easy to verify that this means that $p^i|\alpha,\beta,\gamma$,
so that $g\in G^{p^i}G'$. Therefore, given any pair of elements, if they
fail both conditions (a) and (b) in Theorem~\ref{absclosedgen}, there
is an overgroup $K$ where they are $p^i$-th powers times a commutator;
that means they are $p^i$-th powers in $G$, and by
Lemma~\ref{ifithasaroot}, that suffices for condition (c) to hold.

To prove that $G$ is not absolutely closed in $(p^{a+b},p^{a+1})$,
consider the elements $x,z\in G$, and set
$q=p$. They do not satisfy condition~(a) from
Theorem~\ref{absclosedgen}, since $x^{p^{a+b-1}}=e$ and $b\geq 1$.

For condition (c), note that if $h^p\equiv
x^{\alpha}z^{\beta}\pmod{G^{p^{a+1}}G'}$, then in particular both
$\alpha$ and $\beta$ must be multiples of $p$. If we could solve the
congruences in condition (c), we would have that $p$ divides
$\alpha$, $\beta$, $\gamma$, and $\delta$, but $\gamma-\beta$ must be
congruent to $1$ modulo $p^a$, which is clearly impossible.

Finally, for condition (b), if $h^p\equiv
x^{\alpha}z^{\beta}\pmod{G^{p^{a+1}}G'}$, and we have $h\equiv x^ry^sz^t
\pmod{G'}$, then $sp$ is a multiple of $p^{a+1}$, so $s$ must be a
multiple of $p^a$. Therefore, $[h,x]=[h,z]=e$, and condition~(b)
cannot hold either. Thus, $G$ is
not absolutely closed in $(p^{a+b},p^{a+1})$.
\end{example}

\begin{example}
\label{bsmallspecial}
Let $p$ be a prime, $a,b>1$, and let $G$ be the $(p^{a+b},p^a)$-group
generated by five elements $x$, $y$, $z$, $r$, and $s$, and satisfying
the following relations, in addition to all relations implied by the
laws of $(p^{a+b},p^a)$:
\begin{eqnarray*}
x^{p^{a+b}} & = & y^{p^{a+b}} = z^{p^{a+b}} = e;\\
\null[y,x]^{p^a} &=& [z,y]^{p^a} = [x,z] = e;\\
\null[x,r]&=&[x,s]=[y,r]=[y,s]=[z,r]=[z,s]=[r,s]=e;\\
r^{p^b} & = & [y,x];\quad s^{p^b}=[y,z].
\end{eqnarray*}
Then $G$ is
absolutely closed in $(p^{a+b},p^a)$, but not in $(p^{a+b+1},p^a)$.

To prove it is not absolutely closed in $(p^{a+b+1},p^a)$, let $K$ be
the group obtained by adjoining non-central $p$-th roots $v$ and $w$
to $r$ and $s$, respectively. We can do this, since $r$ and $s$ are
central, and of order $p^{a+b}$. The resulting group lies in
$(p^{a+b+1},p^a)$, and $[v,w]^{p^2}=e$; however, $[v,w]^p$ lies in the
dominion of $G$ and not in $G$, so $G$ is not absolutely closed.

To prove it is absolutely closed in $(p^{a+b},p^a)$, assume $K$ is an
overgroup of $G$, $K\in(p^{a+b},p^a)$, and that some $g\in G$ lies in
$K^{p^i}K'$. 

Write
$g=x^{\alpha}y^{\beta}z^{\gamma}r^{\delta}s^{\eta}[y,x]^{\zeta}[y,z]^{\epsilon}$.
Since $g^{p^{a-i}}$ must be central, we get that $p^i|\alpha$,
$p^i|\beta$, and $p^i|\gamma$. 

If $i\leq b$, then $g^{p^{a+b-i}}$ is trivial, so $p^i|\delta$ and
$p^i|\eta$, proving that $g\in G^{p^i}G'$. If $i>b$, then $g^{p^a}$ is
trivial, so $p^b|\delta$ and $p^b|\eta$. But then
$r^{\delta}s^{\eta}\in G'$, so again we have that $g\in G^{p^i}G'$. In
particular, any pair of elements will satisfy condition (c) in
Theorem~\ref{absclosedgen}, so $G$ is absolutely closed in $(p^{a+b},p^a)$.
\end{example}

We pass now to the varieties where $m=0$, and also present limiting
examples in that case:

\begin{example}
\label{newprimesquare}
Let $n>0$, and $p$ be a prime, with ${\rm gcd}(p,n)=1$. The group
\[ G = \frac{Z}{pZ}\oplus \frac{Z}{pZ}\]
is absolutely closed in $(0,n)$, but not in $(0,p^2n)$.

Indeed, $G$ is $n$-divisible, since $p$ does not divide $n$, so it is
absolutely closed in $(0,n)$. However, $G/G^p\cong G$ is isomorphic to
the sum of two cyclic groups of order $p$, so it cannot be absolutely
closed in $(0,p^2n)$ by Theorem~\ref{abspecialinzero}.
\end{example}

\begin{example}
\label{oldprime}
Let $n>0$, $p$ a prime with $p|n$, and let $a={\rm
ord}_p(n)$. Then
\[ G = \Bigl\langle x,y \,\Bigm|\, x^{p^a}=y^{p^a}=[x,y]^{p^{a-1}} =
[x,y,x] = [x,y,y] = e \Bigr\rangle\]
is absolutely closed in $(0,n)$, but not in $(0,pn)$.

To verify it is absolutely closed in $(0,n)$, let $K$ be an overgroup
of $G$, $K\in (0,n)$, let $q$ be a prime power dividing $n$, and let
$r,s\in K$, $r',s'\in K^nK'$ with $r^qr',s^qs'\in G$. We want to prove
that $[r,s]^q\in G$.

If ${\rm gcd}(p,q)=1$, then $G$ is $q$-divisible,
so both $r^qr'$
and $s^qs'$ have $q$-th roots in $G$, and then we apply
Lemma~\ref{ifithasaroot} to conclude that $[r,s]^q\in G$.

If $p|q$, write $q=p^b$, $b\leq a$. We proceed as
we did in Example~\ref{advanceinboth}. Let $\mathbf{x}=r^qr'$,
$\mathbf{y}=s^qs'$. Since their $\frac{n}{q}$-th powers are central, they are
each $p^{b-1}$-th powers modulo commutators, so we may write:
\begin{eqnarray*}
\mathbf{x} & \equiv & x^{\alpha p^{b-1}} y^{\beta p^{b-1}} \pmod{G^nG'}\\
\mathbf{y} & \equiv & x^{\gamma p^{b-1}} y ^{\delta p^{b-1}} \pmod{G^nG'}.
\end{eqnarray*}

We claim that $[r,s]^{p^a}=e$: since $K\in(0,n)$, we have
\[ \left([r,s]^{p^a}\right)^{n/p^a} = [r,s]^n = e.\]
So the order of $[r,s]^{p^a}$ divides $n/p^a$. In addition, we also have that
\[ \left([r,s]^{p^a}\right)^{p^{2b-1}} = [r^{p^b},s^{p^b}]^{p^{a-1}}
[\mathbf{x},\mathbf{y}]^{p^{a-1}}=e,\]
so the order also divides $p^{2b-1}$. Therefore, the order of
$[r,s]^{p^a}$ divides ${\rm gcd}(n/p^a,p^{2b-1})=1$, since $a={\rm
ord}_p(n)$. So $[r,s]^{p^a}$ is in fact trivial, as claimed.

Consider the congruences
\begin{eqnarray*}
e^{p^b} & \equiv & \mathbf{x}^0\mathbf{y}^0\pmod{K^nK'}\\
\left(x^{\alpha p^{a-b-1}}y^{\beta p^{a-b-1}}\right)^{p^b} & \equiv &
\mathbf{x}^{p^{a-b}}\mathbf{y}^0 \pmod{K^nK'}.
\end{eqnarray*}
\noindent By Lemma~\ref{basiceqs}, this gives that
\begin{eqnarray*}
{} [r,s]^{p^a}  = [r,s]^{p^b(p^{a-b})}
& = & [x^{\alpha p^{a-b-1}}y^{\beta
p^{a-b-1}},x^{\gamma p^{a-1}}y^{\delta p^{a-1}}]\\
& = & [x,y]^{p^{a-2}(\alpha\delta-\beta\gamma)}.
\end{eqnarray*}
Since $[r,s]^{p^a}=e$, this means that $p|\alpha\delta-\beta\gamma$.
That means that by perturbing $\mathbf{x}$ and $\mathbf{y}$, we may
assume that $\langle\mathbf{x},\mathbf{y}\rangle$ is cyclic, hence
absolutely closed, so
\[ [r,s]^q\in {\rm
dom}_K^{(0,n)}\left(\langle\mathbf{x},\mathbf{y}\rangle \right) =
\langle \mathbf{x},\mathbf{y}\rangle \subseteq G.\]
So $G$ is absolutely closed in $(0,n)$.

But as we also showed in Example~\ref{advanceinboth}, $G$ is not
absolutely closed in $(p^{a+1},p^{a+1})\subseteq (0,np)$, so it cannot
be absolutely closed in $(0,np)$ either.
\end{example}

The final case to consider is covered by the following result:

\begin{theorem}
\label{newprimesingle}
Let $n>0$, and $p$ be a prime with ${\rm gcd}(n,p)=1$. If $G$ is
absolutely closed in $(0,n)$, then it is also absolutely closed in
$(0,np)$.
\end{theorem}

\begin{proof} Let $K$ be any $(0,np)$-overgroup of $G$; let $r,s\in
K$, $r',s'\in K^{np}K'$, $q>0$ be a prime power with $q|np$, and assume that
$r^qr', s^qs'\in G$. We want to prove that $[r,s]^q\in G$ as well.

If $q=p$, then note that $[r,s]^{p^2}\in G$. Since $p$ does not divide
$n$, there exist an integer $k$ such that $kp\equiv
1\pmod{n}$. Therefore, $kp^2\equiv p\pmod{np}$. If we raise
$[r,s]^{p^2}$ to the $k$-th power, we are still in $G$; but this
equals $[r,s]^{p^2k}$, which equals $[r,s]^p$, because $[r,s]$ is of
exponent $np$. Therefore, $[r,s]^p\in G$, as desired.

If ${\rm gcd}(q,p)=1$, things are somewhat more complicated. By adding
central $q$-th roots $t$ and $v$ to $r'$ and $s'$, respectively, and
replacing $r$ with $rt$ and $s$ with $sv$, we may assume that
$r'=s'=e$. Consider the subgroup $K_2$ generated by $G$, $r$, and
$s$. Note that $[r,s]^q$ is still in the dominion of $G$ in $K_2$. We
claim that $K_2\in(0,n)$.

To prove this it suffices to show that $[r,s]$ is of exponent $n$, as
are $[g,r]$ and $[g,s]$ for each $g\in G$. To see that $[r,s]$ is of
exponent $n$, note that $[r,s]^{np} = e$, since $K\in (0,np)$; but
also
\[ [r,s]^{q^2n} = [r^q,s^q]^n = e\]
since $r^q,s^q\in G$ and $n$-th powers of elements of $G$ are
central. Therefore, the order of $[r,s]$ divides ${\rm gcd}(q^2n,np) =
n\left({\rm gcd}(q^2,p)\right)=n$.

Now let $g\in G$ be an arbitrary element. Since $K\in(0,np)$,
$[r,g]^{np}=e$. However,
\[ e = [r,g]^{np} = \left( [r,g]^n\right)^p =
\left([r^q,g]^{n/q}\right)^p.\] Thus, $[r^q,g]^{n/q}=[r,g]^n$ is a
commutator in $G$ of exponent $p$. Since ${\rm gcd}(p,n)=1$, and
$G\in(0,n)$, that means that $[r,g]^n=e$. A symmetric argument shows
that $[g,s]^n=e$ as well. Thus, $K_2\in(0,n)$, and it is an overgroup
of $G$. Since $G$ is absolutely closed in $(0,n)$, and $[r,s]^q$ lies
in the dominion of $G$ in $K_2$, it follows that $[r,s]^q\in G$, which
is what we wanted to prove. Therefore, $G$ is absolutely closed in
$(0,np)$, as claimed.
\end{proof}

With these examples in hand, we prove the analogue of
Theorem~\ref{howfewremain}.

\begin{theorem}
\label{howfewremainspecial}
Let $(m,n)\subseteq (m',n')$. Every special amalgamation base in
$(m,n)$ is also a special amalgamation base in $(m',n')$ if and only
if for each prime $p$, one of the following holds:
\begin{itemize}
\item[(a)]~${\rm ord}_p(n),{\rm ord}_p(n')\leq 1$; or
\item[(b)]~${\rm ord}_p(n)={\rm ord}_p(m)={\rm ord}_p(n')$; or
\item[(c)]~${\rm ord}_p(n)={\rm ord}_p(n')$ and ${\rm ord}_p(m)={\rm
ord}_p(m')$.
\end{itemize}
\end{theorem}

\begin{proof}
We prove necessity first. Assume that there is a prime $p$ for which none
of them hold. 

If $m,m'>0$, we may assume both are prime power of the same
prime. Note that for any prime other than that one, condition (c) is
satisfied, so $m$ and $m'$ are powers of $p$ itself. If ${\rm
ord}_p(n)=1$, then we must have ${\rm ord}_p(n')>1$, and we can use
Example~\ref{advanceinboth} to find a $(p,p)$-special base, which is
not a $(p^2,p^2)$-special base, hence not an $(m',n')$-special base;
and the $(p,p)$-base is necessarily an $(m,n)$-special base as
well. If ${\rm ord}_p(n)=0$, then $p^2|n'$, and we may take the sum of
two cyclic groups of order $p$.

If ${\rm ord}_p(n)>1$, then we have that
$p^2|n$, and either ${\rm ord}_p(n)<{\rm ord}_p(n')$, or ${\rm
ord}_p(n)<{\rm ord}_p(m)<{\rm ord}_p(m')$.

If ${\rm ord}_p(n)<{\rm ord}_p(n')$, then we may use
Example~\ref{advanceinboth} if ${\rm ord}_p(m)={\rm ord}_p(n)$, and
Example~\ref{advancesecondsp} if ${\rm ord}_p(n)<{\rm ord}_p(m)$, to see
that not every $(m,n)$-special base is an $(m',n')$-special base.

If ${\rm ord}_p(n)={\rm ord}_p(n')$, and ${\rm ord}_p(n)<{\rm
ord}_p(m)<{\rm ord}_p(m')$, then we use Example~\ref{bbigspecial} or
Example~\ref{bsmallspecial} to show that not every $(m,n)$-absolutely
closed group is absolutely closed in $(m',n')$. This proves necessity
when $m,m'>0$.

If $m=m'=0$, we have $(0,n)\subseteq (0,n')$. Since $p$ fails all
three conditions, ${\rm ord}_p(n)<{\rm ord}_p(n')$, and either $p$
does not divide $n$ and $p^2|n'$, in which case we use
Example~\ref{newprimesquare}, or else ${\rm ord}_p(n)=1$ and $p^2|n'$,
in which case we may use Example~\ref{oldprime} instead.

If $m>0$ and $m'=0$, we may assume that $m$ and $n$ are prime powers;
if they are powers of some prime $q\not=p$, then $q^2|n'$ and we use
Example~\ref{newprimesquare} again. If $q=p$, then we proceed as in
the case where both $m$ and $m'$ are nonzero. This proves necessity.

For sufficiency, let $(m,n)\subseteq(m',n')$, and assume that for
every prime $p$ the conditions are all satisfied.

If $m,m'>0$, we may assume they are prime powers. If (a) holds, then
in both $(m,n)$ and $(m',n')$ all groups are absolutely closed and
there is nothing to do. If (c) holds, then $(m,n)=(m',n')$ and again
the result is trivial. Finally, if (b) holds, then we are going from a
$(p^a,p^a)$ to a $(p^{a+b},p^a)$, and everything which is absolutely
closed in the former is also absolutely closed in the latter by
Proposition~\ref{bequalaspecial}. 

If $m=m'=0$, and (b) holds for any prime, then $n=n'=0$ and the result
is trivial. Otherwise, for each prime either (a) or (c) holds, so we
may write $n'=nk$, with ${\rm gcd}(n,k)=1$, and $k$ squarefree. In
this case, the result follows from Theorem~\ref{newprimesingle}.

Finally, if $m>0$ and $m'=0$, we may assume $m$ and $n$ are powers
of~$p$. Note that condition (c) cannot hold. Given any prime $q$ other
than~$p$, the conditions imply that $q^2$ does not divide $n'$. For
$p$, if condition (b) holds, then we may write $n'=nk$ again, with $k$
squarefree and relatively prime to $n$, and $(m,n)=(p^a,p^a)$, The
result now follows by applying Proposition~\ref{bequalaspecial} and
Theorem~\ref{newprimesingle}. On the other hand, if (a) holds, then
all groups in $(m,n)$ are special amalgamation bases, and since $p^2$
does not divide $n'$, every $p$-group in $(0,n')$ also is a special
base, and we are done.

\end{proof}

\subsection*{Absolute Closures}

We end by considering absolute closures. By a theorem of Isbell
(Corollary~1.8 in~\cite{isbellone}, modified to fit our setting) in a
variety ${\cal V}$, every group $G\in{\cal V}$ can be embedded into a
group $K$ which is absolutely closed in ${\cal V}$. Isbell then gives
the following definition, which we have adapted to our setting:

\begin{defb} Let ${\cal V}$ be a variety of groups, and $G\in {\cal
V}$. We say that $G$ has an \textbf{absolute closure} if there exists
an overgroup $K\in{\cal V}$ of $G$ such that ${\rm dom}_K^{\cal V}(G)$
is absolutely closed in ${\cal V}$.
\end{defb}

That is, an absolute closure for $G$ is an absolutely closed overgroup
of $G$ which is dominated by $G$ in some further overgroup. It is not
hard to give examples of algebras which are not absolutely closed but
have absolute closures in some categories of algebras; for example,
any group is absolutely closed in the category of semigroups, and the
dominion of a subsemigroup of a group in the group itself is the
subgroup it generates. On the other hand, there are examples of
algebras which have no absolute closure.

In the varieties $(m,n)$, it turns out that no group which is not
already absolutely closed can have an absolute closure: 

\begin{theorem}[cf. Theorem~3.6 in~\cite{closures}]
A group $G\in(m,n)$ has an absolute closure in $(m,n)$ if and only if it is
absolutely closed in $(m,n)$.
\end{theorem}

\begin{proof} One implication is trivial. So assume that $G$ is not
absolutely closed. Let $K\in(m,n)$ be an overgroup, and let $D_K={\rm
dom}_K^{(m,n)}(G)$. We want to show that $D_K$ is not absolutely
closed.

If $D_K=G$ there is nothing to do. So assume that $D_K\not=G$. In
particular, there exist $q>0$, $q|n$, $r,s\in K$, $r',s'\in K^nK'$
with $r^qr',s^qs'\in G$, and $[r,s]^q\notin G$.  Let $x=r^qr'$,
$y=s^qs'$. We claim that $x$, $y$, and $q$ do not satisfy conditions
(a), (b), or (c) in Theorem~\ref{absclosedgen}, relative to $D_K$ and
$(m,n)$.

If condition (a) or condition (b) holds, then no overgroup $M$ of
$D_K$ has $x,y\in M^qM'$; however, $x,y\in K^qK'$, so neither (a) nor
(b) can hold. That leaves only (c), and we proceed by
contradiction. Assume that we have $d_1,d_2\in D_K$, and integers
$\alpha,\beta,\gamma,\delta$, with $\gamma-\beta\equiv 1
\pmod{\frac{n}{q}}$, and such that
\begin{eqnarray*}
d_1^q &\equiv& x^{\alpha}y^{\beta}\pmod{D_K^n D_K'}\\
d_2^q &\equiv& x^{\gamma}y^{\delta}\pmod{D_K^n D_K'}.
\end{eqnarray*}
By Lemma~\ref{basiceqs}, it follows that in $K$ we have 
\[ [d_1,x][d_2,y]=[r,s]^{q(\gamma-\beta)} = [r,s]^q.\] 

Note that since $D_K$ is obtained from $G$ by adjoining central
elements, all of which are of exponent $n$, $D_K^n=G^n$ and $D_K'=G'$,
so the congruences above hold modulo $G^nG'$ as well. By the
description of dominions, we may write
\[ d_i = g_i [k_{i11},k_{i12}]^{q_{i1}}\cdots
[k_{is_i1},k_{is_i2}]^{q_{is_i}}\]
for $i=1,2$, where $k_{ijk}\in K$, and
$k_{ij1}^{q_{ij}},k_{ij2}^{q_{ij}}\in G(K^nK')$. Since $d_i\equiv
g_i\pmod{K'}$, we have
\[ [d_1,x][d_2,y]=[g_1,x][g_2,y]\in G\]
which means that $[r,s]^q\in G$, which contradicts our choice of $r$
and $s$. Therefore, condition (c) cannot be satisfied in $D_K$ either,
so $D_K$ is not absolutely closed.
\end{proof}



\bibliography{bibliog}

\providecommand{\bysame}{\leavevmode\hbox to3em{\hrulefill}\thinspace}
\begin{thebibliography}{10}

\bibitem{episandamalgs}
Peter~M. Higgins, \emph{Epimorphisms and amalgams}, Colloq. Math. \textbf{56}
  (1988), no.~1, 1--17, {MR {\bf 89m}:20083}.

\bibitem{higmanpgroups}
Graham Higman, \emph{Amalgams of {$p$}-groups}, J. Algebra \textbf{1} (1964),
  301--305, {MR {\bf 29}:4799}.

\bibitem{isbellone}
J.R. Isbell, \emph{Epimorphisms and dominions}, Proc. of the Conference on
  Categorical Algebra, La Jolla 1965 (S.~Eilenberg et~al, ed.), Lange and
  Springer, New York, 1966, {MR {\bf 35}:105a} ({T}he statement of the {Z}igzag
  {L}emma for {\it rings} in this paper is incorrect. {T}he correct version is
  in \cite{isbellfour}), pp.~232--246.

\bibitem{isbellfour}
J.R. Isbell, \emph{Epimorphisms and dominions {IV}}, J. London Math. Soc. (2)
  \textbf{1} (1969), 265--273, {MR {\bf 41}:1774}.

\bibitem{jonsson}
Bjarni J{\'o}nsson, \emph{Varieties of groups of nilpotency three}, Notices
  Amer. Math. Soc. \textbf{13} (1966), 488.

\bibitem{leinen}
Felix Leinen, \emph{An amalgamation theorem for soluble groups}, Canad. Math.
  Bull. \textbf{30} (1987), no.~1, 9--18, {MR {\bf 88d}:20050}.

\bibitem{machenry}
T.~Mac{H}enry, \emph{The tensor product and the 2nd nilpotent product of
  groups}, Math. Z. \textbf{73} (1960), 134--145, {MR {\bf 22}:11027a}.

\bibitem{closures}
Arturo Magidin, \emph{Amalgamation bases for nil-2 groups of odd exponent},
  {arXiv:math.GR/0006065}.

\bibitem{absclosed}
Arturo Magidin, \emph{Absolutely closed nil-$2$ groups}, Algebra Universalis
  \textbf{42} (1999), no.~1-2, 61--77, {MR {\bf 2001a}:20055}.

\bibitem{fgnil}
Arturo Magidin, \emph{Dominions in finitely generated nilpotent
  groups}, Comm. Algebra \textbf{27} (1999), no.~9, 4545--4559, {MR
{\bf 2001f}:20070}.

\bibitem{nildoms}
Arturo Magidin, \emph{Dominions in varieties of nilpotent groups},
  Comm. Algebra \textbf{28} (2000), no.~3, 1241--1270, {MR {\bf 2000m}:20053}.

\bibitem{amalgone}
Berthold~J. Maier, \emph{Amalgame nilpotenter {G}ruppen der {K}lasse zwei},
  Publ. Math. Debrecen \textbf{31} (1985), 57--70, {MR {\bf 85k}:20117}.

\bibitem{amalgtwo}
Berthold~J. Maier, \emph{Amalgame nilpotenter {G}ruppen der {K}lasse
  zwei {II}}, Publ. Math. Debrecen \textbf{33} (1986), 43--52, {MR
  {\bf 87k}:20050}. 

\bibitem{nilexpp}
Berthold~J. Maier, \emph{On nilpotent groups of exponent {$p$}},
  J. Algebra \textbf{127} (1989), 279--289, {MR {\bf 91b}:20046}.

\bibitem{hneumann}
Hanna Neumann, \emph{Varieties of groups}, Ergebnisse der Mathematik und ihrer
  Grenzgebiete, New Series, vol.~37, Springer Verlag, 1967, {MR {\bf 35}:6734}.

\bibitem{classifthree}
V.N. Remeslennikov, \emph{Two remarks on 3-step nilpotent groups}, Algebra i
  Logika (1965), no.~2, 59--65 (Russian), {MR {\bf 31}:4838}.

\bibitem{saracino}
D.~Saracino, \emph{Amalgamation bases for nil{-$2$} groups}, Algebra
  Universalis \textbf{16} (1983), 47--62, {MR {\bf 84i}:20035}.

\end{thebibliography}

\bibliographystyle{amsplain}



\end{article}
\end{document}